\documentclass{gtart_a}
\pdfoutput=1

\usepackage{graphicx,import,overpic,wrapfig}
\usepackage[import, rotate]{xy}

\title[A random tunnel number one 3--manifold does not fiber over the circle]{A random tunnel number one 3--manifold\\ does not fiber over the circle }

\author{Nathan M Dunfield}
\givenname{Nathan M}
\surname{Dunfield}
\address{Mathematics 253-37\\California Institute of Technology\\\newline
Pasadena, CA 91125\\USA}
\email{dunfield@caltech.edu}
\urladdr{}

\author{Dylan P Thurston}
\givenname{Dylan P}
\surname{Thurston}
\address{Mathematics\\Barnard College\\
Columbia University MC 4436\\\newline
New York, NY 10027\\USA}
\email{dthurston@barnard.edu}
\urladdr{}

\volumenumber{10}
\issuenumber{}
\publicationyear{2006}
\papernumber{55}
\lognumber{0729}
\startpage{2431}
\endpage{2499}

\doi{}
\MR{}
\Zbl{}


\received{8 April 2006}
\revised{}
\accepted{13 November 2006}
\proposed{Cameron Gordon}
\seconded{Rob Kirby, Joan Birman}
\published{15 December 2006}
\publishedonline{15 December 2006}
\corresponding{}
\editor{CPR}
\version{}

\arxivreference{math.GT/0510129}

\let\xysavmatrix\xymatrix
\def\xymatrix{\disablesubscriptcorrection\xysavmatrix}
\AtBeginDocument{\let\bar\wbar\let\tilde\wtilde\let\hat\what}

\newcommand{\G}{{\mathcal G}}
\newcommand{\T}{{\mathcal T}}

\newcommand{\maps}{\colon\thinspace}

\DeclareMathOperator{\rev}{rev}

\newcommand{\abs}[1]{{\left| #1 \right|}}
\newcommand{\pair}[1]{\left\langle #1 \right\rangle}

\newcommand{\spandef}[2]{{  \left\langle  {#1}  \ \left| \   {#2} \right. \right\rangle }}
\newcommand{\setdef}[2]{{  \left\{  {#1}  \ \left| \   {#2} \right. \right\} }}
\newcommand{\mtext}[1]{\quad\mbox{#1}\quad}

\newcommand{\MCG}{\mathcal{MCG}}
\newcommand{\ML}{\mathcal{ML}}
\newcommand{\PML}{\mathcal{PML}}

\newcommand{\W}{\mathcal{W}}
\newcommand{\B}{\mathcal{B}}
\newcommand{\cS}{\mathcal{S}}
\newcommand{\cI}{\mathcal{I}}
\newcommand{\cG}{\mathcal{G}}
\newcommand{\cC}{\mathcal{C}}
\newcommand{\boxmap}{\mathrm{Box}}
\newcommand{\num}[1]{\#\left({#1} \right)}

\newenvironment{xyoverpic*}[3]
{%
\begin{xy}
\xyimport#1{\includegraphics[#2]{#3}}
}{\end{xy}}

\newenvironment{xyoverpic}[2][]
{\begin{xyoverpic*}{(100,100)(0,0)}{#1}{#2}}
{\end{xyoverpic*}}

\makeatletter
\def\cnewtheorem#1[#2]#3{\newtheorem{#1}{#3}[section]
\expandafter\let\csname c@#1\endcsname\c@subsection}

\swapnumbers
\cnewtheorem{theorem}[subsection]{Theorem}
\cnewtheorem{conjecture}[subsection]{Conjecture}
\cnewtheorem{lemma}[subsection]{Lemma}
\cnewtheorem{corollary}[subsection]{Corollary}
\cnewtheorem{proposition}[subsection]{Proposition}
\cnewtheorem{claim}[subsection]{Claim}
\cnewtheorem{question}[subsection]{Question}
\cnewtheorem{criterion}[subsection]{Criterion}

\cnewtheorem{virtual_haken}[subsection]{Virtual Haken Conjecture}
\cnewtheorem{vfconj}[subsection]{Virtual Fibration Conjecture}
\newtheorem*{mainforward}{\fullref{thm-main-nonsep}}
\newtheorem*{mainforwardsep}{\fullref{thm-main-sep}}
\newtheorem*{forwardwalk}{\fullref{thm-random-group}}
\makeautorefname{virtual_haken}{Conjecture}
\makeautorefname{vfconj}{Conjecture}

\theoremstyle{definition}
\cnewtheorem{definition}[subsection]{Definition}
\cnewtheorem{notation}[subsection]{Notation}

\theoremstyle{remark}
\cnewtheorem{remark}[subsection]{Remark}
\cnewtheorem{example}[subsection]{Example}

  \let\c@equation\c@subsection
  \def\theequation{\bf \thesubsection}
\makeatother

\begin{document}

\begin{htmlabstract}
<p class="noindent">
We address the question: how common is it for a 3&ndash;manifold to fiber
over the circle?  One motivation for considering this is to give
insight into the fairly inscrutable Virtual Fibration Conjecture.
For the special class of 3&ndash;manifolds with tunnel number one, we
provide compelling theoretical and experimental evidence that
fibering is a very rare property.  Indeed, in various precise senses
it happens with probability 0.  Our main theorem is that this is
true for a measured lamination model of random tunnel number one
3&ndash;manifolds.
</p>
<p class="noindent">
The first ingredient is an algorithm of
K Brown which can decide if a given tunnel number one 3&ndash;manifold
fibers over the circle.  Following the lead of Agol, Hass and W
Thurston, we implement Brown's algorithm very efficiently by working
in the context of train tracks/interval exchanges.  To analyze the
resulting algorithm, we generalize work of Kerckhoff to understand
the dynamics of splitting sequences of complete genus 2 interval
exchanges.  Combining all of this with a &ldquo;magic splitting
sequence&rdquo; and work of Mirzakhani proves the main theorem.
</p>
<p class="noindent">
The 3&ndash;manifold situation contrasts markedly with random
2&ndash;generator 1&ndash;relator groups; in particular, we show that
such groups &ldquo;fiber&rdquo; with probability strictly between&nbsp;0
and&nbsp;1.
</p>
<h5>Additional material for downloading</h5>
<p class="noindent">
The first two items check certain combinatorial facts asserted in
Section 10.  The other two are the binary and source code for the
genus2fiber program.</p>
<p class="indent">
<a href="p055/all_int.gap">all_int.gap</a><br>
<a href="p055/exchange2.py">exchange2.py</a><br>
<a href="p055/genus2fiber-bin.tar">genus2fiber-bin.tar</a><br>
<a href="p055/genus2fiber-src.tar">genus2fiber-src.tar</a>
</p>
\end{htmlabstract}

\begin{abstract} 
  We address the question: how common is it for a 3--manifold to fiber
  over the circle?  One motivation for considering this is to give
  insight into the fairly inscrutable Virtual Fibration Conjecture.
  For the special class of 3--manifolds with tunnel number one, we
  provide compelling theoretical and experimental evidence that
  fibering is a very rare property.  Indeed, in various precise senses
  it happens with probability 0.  Our main theorem is that this is
  true for a measured lamination model of random tunnel number one
  3--manifolds.
  
  The first ingredient is an algorithm of
  K Brown which can decide if a given tunnel number one 3--manifold
  fibers over the circle.  Following the lead of Agol, Hass and W
  Thurston, we implement Brown's algorithm very efficiently by working
  in the context of train tracks/interval exchanges.  To analyze the
  resulting algorithm, we generalize work of Kerckhoff to understand
  the dynamics of splitting sequences of complete genus 2 interval
  exchanges.  Combining all of this with a ``magic splitting sequence''
  and work of Mirzakhani proves the main theorem.
  
  The 3--manifold situation contrasts markedly with random 2--generator
  1--relator groups; in particular, we show that
  such groups ``fiber'' with probability strictly between~0 and~1.
\end{abstract}

\maketitle

{\small\leftskip=25pt\rightskip25pt\it
We dedicate this paper to the memory of Raoul Bott (1923--2005), a
wise teacher and warm friend, always searching for the simplicity at
the heart of mathematics.\par}

\section{Introduction}
\label{sec:intro}

In this paper we are interested in compact orientable 3--manifolds whose
boundary, if any, is a union of tori.  A nice class of such manifolds
are those that \emph{fiber over the circle}, that is, are fiber
bundles over the circle with fiber a surface $F$:
\[
F \to M \to S^1
\]
Equivalently, $M$ can be constructed by taking $F \times [0,1]$ and
gluing $F \times \{0\}$ to $F \times \{ 1 \}$ by a homeomorphism $\psi$ of $F$.
Manifolds which fiber over the circle are usually easier to understand
than 3--manifolds in general, because many questions can be
reduced to purely 2--dimensional questions about the gluing map $\psi$.  

When $M$ fibers over the circle, the group $H^1(M ; \Z) \cong H_2(M, \partial
M; \Z)$ is nonzero; a nontrivial element is the fibering map to $S^1$,
or dually, the fiber surface $F$.  Our main question is:
\begin{question}\label{vague-question}
  If we suppose $H^1(M; \Z) \neq 0$, how common is it for $M$ to fiber
  over the circle?
\end{question}
\noindent
We will give several reasons for entertaining this question below, but
for now one motivation (beyond its inherent interest) is to try to
estimate how much harder the Virtual Fibration Conjecture is than
other variants of the Virtual Haken Conjecture.  In this
paper we provide evidence, both theoretical and experimental, that
the answer to \fullref{vague-question} is: not very common at
all.  In fact, for the limited category of 3--manifolds that we study
here, the probability of fibering is 0.

The type of 3--manifold we focus on here are those with tunnel number
one, which we now define.  Let $H$ be an orientable handlebody of
genus $2$, and pick an essential simple closed curve $\gamma$ on $\partial H$.
Now build a 3--manifold $M$ by gluing a 2--handle to $\partial H$ along
$\gamma$; that is, $M = H \cup ( D^2 \times I ) $ where $\partial D^2 \times I$ is glued to
$\partial H$ along a regular neighborhood of $\gamma$.  Such a manifold is
said to have \emph{tunnel number one}.  There are two kinds of these
manifolds, depending on whether the curve $\gamma$ is separating or not.
For concreteness, let us focus on those where $\gamma$ is non-separating.
In this case, $M$ has one boundary component, which is a torus.  A
simple example of a 3--manifold with tunnel number one is the exterior
of a 2--bridge knot in $S^3$.

The boundary of a tunnel number one manifold $M$ forces $H^1(M; \Z) \neq
0$, and so it makes sense to consider \fullref{vague-question}
for all manifolds in this class.  To make this question more precise, we
will need a notion of a ``random'' tunnel number one manifold, so that
we can talk about probabilities.  In fact, there are several
reasonable notions for this; here, we focus on two which involve
selecting the attaching curve $\gamma \subset \partial H$ from the point of view of
either measured laminations or the mapping class group.  

For measured laminations, the setup is roughly this.  We fix
Dehn--Thurston coordinates on the set of multicurves (equivalently
integral measured laminations) on the surface~$\partial H$.  Let $\T(r)$
consist of the tunnel number one 3--manifolds whose attaching curve
$\gamma \subset \partial H$ has all coordinates of size less than $r$.  As $\T(r)$ is
finite, it makes sense to formulate a precise version of
\fullref{vague-question} as: what is the proportion of $M \in
\T(r)$ which fiber over the circle when $r$ is large?  The main
theorem of this paper is:
\begin{mainforward}
  The probability that $M \in \T(r)$ fibers over the circle goes to $0$
  as $r \to \infty$.
\end{mainforward}
\noindent
Thus with this notion of a random tunnel number one 3--manifold, being
fibered is very rare indeed.  There is one technical caveat here: the
set $\T(r)$ we consider does not cover all multicurves on $\partial H$,
although we can always change coordinates, preserving~$H$, to put any
curve in $\T(r)$.  See the discussion in
\fullref{subsec:def-of-tunnel-num}.

Another natural model for random tunnel number one 3--manifolds is
to create them using the mapping class group of the surface $\partial H$.
More precisely, fix a finite generating set $S = \{ \psi_1, \psi_2, \ldots,
\psi_n \}$ of the mapping class group $\MCG(\partial H)$; for instance, take
$S$ to be the standard five Dehn twists.  Fix also a non-separating simple
closed curve $\gamma_0$ on $\partial H$.  Now given $r$, create a sequence
$\phi_1, \phi_2, \ldots, \phi_r$ by picking each $\phi_i$ at random from among the
elements of $S$ and their inverses.  Then set
\[
\gamma = \phi_r \circ \phi_{r-1} \circ \cdots \circ \phi_1(\gamma_0),
\]
and consider the corresponding tunnel number one manifold $M$.  That
is, we start with $\gamma_0$ and successively mess it up $r$ times by
randomly chosen generators.  Equivalently, we go for a random walk in
the Caley graph of $\MCG(\partial H)$, and then apply the endpoint of that
walk to $\gamma_0$ to get $\gamma$.  \fullref{vague-question} now
becomes: what is the probability that $M$ fibers over the circle if
$r$ is large?  A priori, the answer could be different from the one
given in \fullref{thm-main-nonsep}.  Because the number of such
manifolds is countably infinite, there is no canonical probability
measure on this set, so our choice of model for a random manifold is
important.  One might hope that all ``reasonable'' models give the
same answer, but it should be emphasized that in some ways our two
notions are fundamentally different.  In any event, we will provide
compelling experimental evidence for the following conjecture.
\begin{conjecture}\label{conj-mcg-nofiber}
  Let $M$ be a tunnel number one 3--manifold created by a random walk
  in $\MCG(\partial H)$ of length $r$.  Then the probability that $M$ fibers
  over the circle goes to $0$ as $r \to \infty$.
\end{conjecture}
\noindent
Thus from this alternate point of view as well, it seems that nearly
all tunnel number one 3--manifolds do not fiber over the circle.

\subsection{Random groups}\label{sub-intro-groups}

One of the fundamental tasks of 3--dimensional topology is to
understand the special properties of their fundamental groups, as
compared to finitely presented groups in general.
From the point of view of this paper there is a surprising
contrast between these two classes of groups.  While the question of
whether a 3--manifold $M$ fibers over the circle might seem
fundamentally geometric, Stallings showed that it can be reduced to an
algebraic question about $\pi_1(M)$ (see \fullref{sec-stallings}).
For a group $G$, let us say that $G$ \emph{fibers} if there is an
automorphism $\rho$ of a free group $B$ so that $G$ is the algebraic
mapping torus:
\[
\spandef{ t, B}{\mbox{$t b t^{-1} = \rho(b)$ for all $b \in B$}} 
\]
If $M$ has tunnel number one, then it fibers over the circle if and
only if $\pi_1(M)$ fibers in this sense
(\fullref{cor:tunnel-one-stallings}).  When $M$ has tunnel
number one, its fundamental group is constructed from the free group
$\pi_1(H)$ by killing the attaching curve $\gamma$ of the 2--handle.
Thus the fundamental group is just 
\[
\pi_1(M) = \spandef{\pi_1(H)}{\gamma = 1} = \spandef{a,b}{R = 1},
\]
that is, a 2--generator, 1--relator group.

In the spirit suggested above, we would like to compare
\fullref{thm-main-nonsep} with the situation for a random group
$G$ of the form $\spandef{a,b}{R=1}$.  A natural meaning for the
latter concept would be to consider the set $\cG(r)$ of all such
groups where the length of the relator $R$ is less than $r$.  This
notion is in fact almost precisely analogous to the setup of $\T(r)$
for manifolds; in particular if $M \in \T(r)$, then the natural
presentation of $\pi_1(M)$ is in $\cG(r)$.  Yet the remarkable thing is
that the probability that $G \in \cG(r)$ fibers experimentally tends to
about $0.94$ as $r \to \infty$.  While we can't prove this, we can at least
show:
\begin{forwardwalk}
  Let $p_r$ be the probability that $G \in \cG(r)$ fibers.  Then
  for all large $r$ one has
  \[
  0.0006 \leq p_r \leq 0.975.
  \]
\end{forwardwalk}
\noindent
In particular, $p_r$ does not limit to $0$ as $r \to \infty$, in marked
contrast to \fullref{thm-main-nonsep}.

\begin{figure}
  \centering
  \includegraphics[scale=0.33]{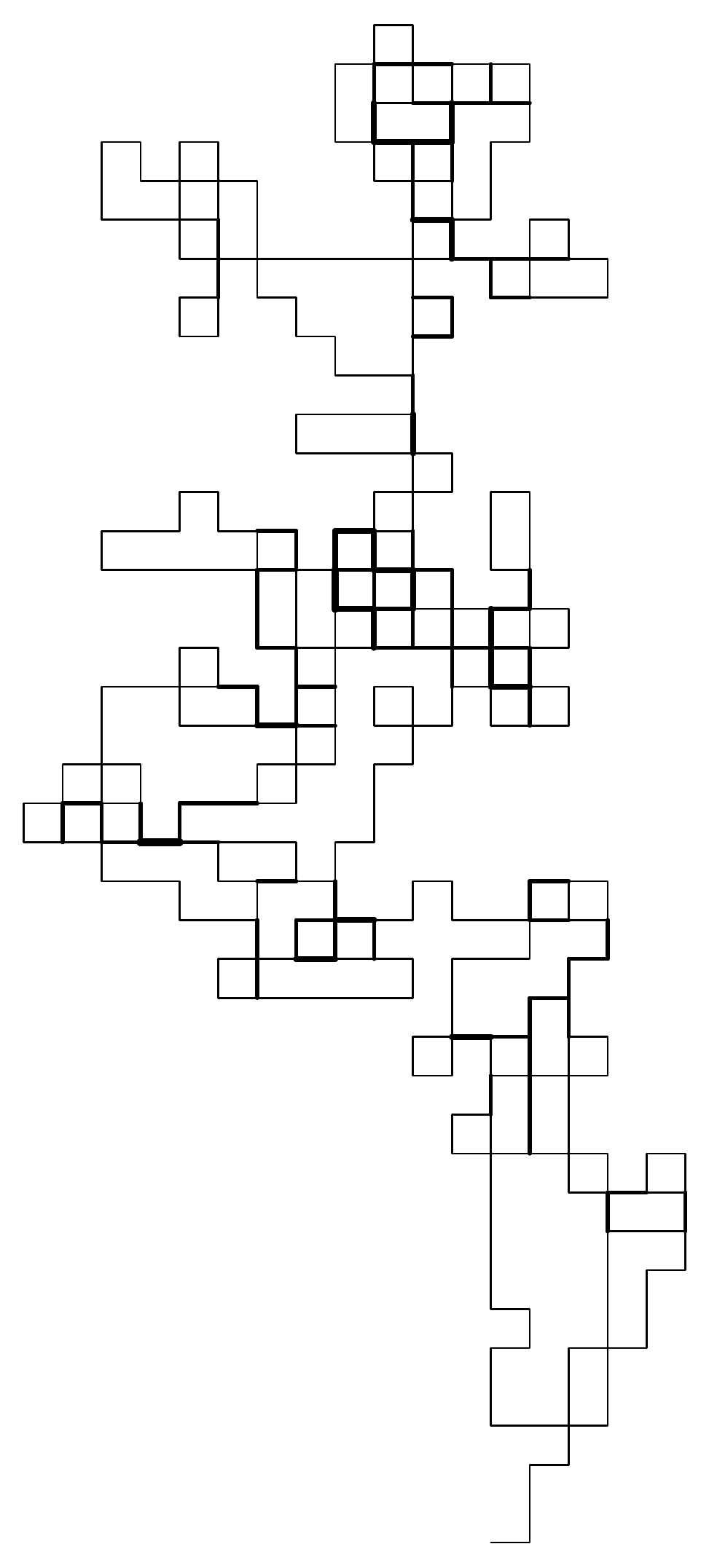}\qquad
  \includegraphics[scale=0.33]
{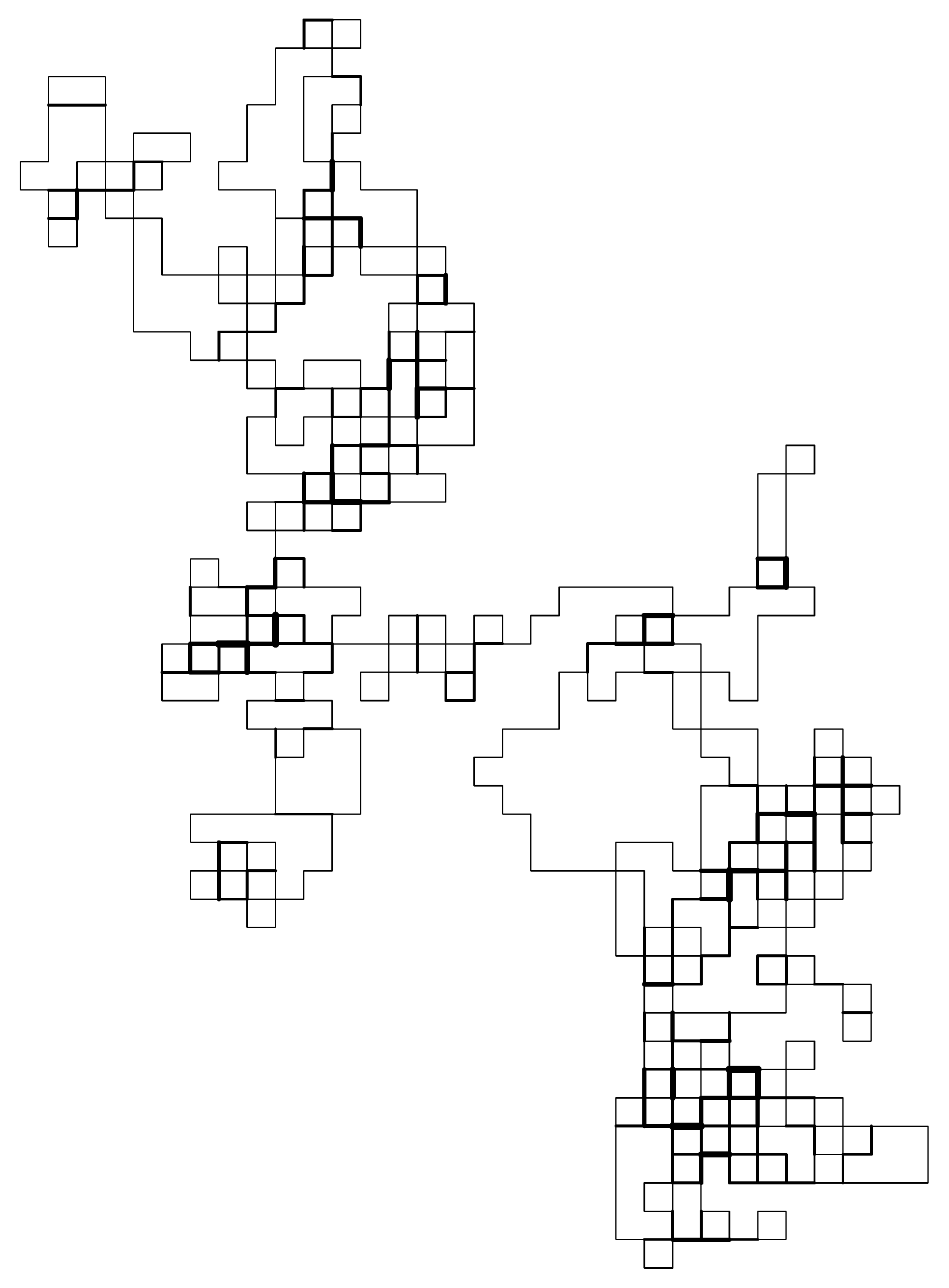}
  \caption{
    Two random words in the free group $F = \pair{a,b}$.  Here a word
    is plotted as a walk in the plane, where $a$ corresponds to a unit
    step in the positive $x$--direction, and $b$ a unit step in the positive
    $y$--direction.  Thicker lines indicate points transversed multiple
    times.  The relevance of these pictures will be made clear in
    \fullref{subsec-Browns-algorithm}.  
 }\label{fig-walk-words}
\end{figure}

\begin{figure}
  \centering
  \includegraphics[scale=0.72]{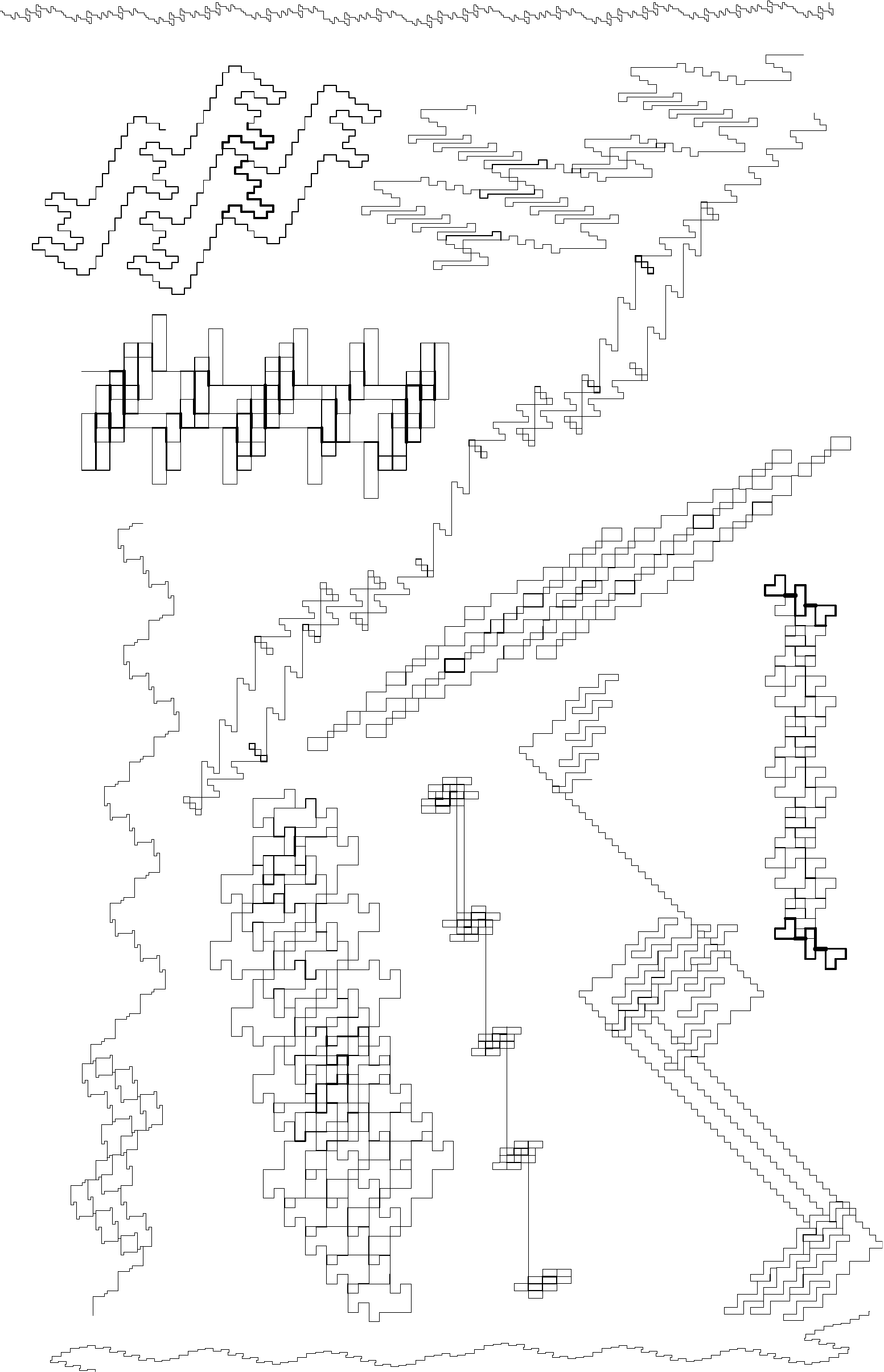}
  \caption{
    Relators of tunnel number one 3--manifolds, plotted in the style of
    \fullref{fig-walk-words}.  To conserve space, they are not all
  drawn to the same scale.}
  \label{fig-mfld-words}
\end{figure}

As we will explain later, whether or not $G = \spandef{a,b}{R=1}$
fibers depends on the combinatorics of the relator $R$ in a certain
geometric sense.  The different behavior for 3--manifold groups comes
down to the fact that the curve $\gamma$ is an \emph{embedded} curve on
the genus 2 surface $\partial H$, and this gives the relator $R$ a recursive
structure where certain ``syllables'' appear repeatedly at varying
scales.  Compare \fullref{fig-walk-words} with \fullref{fig-mfld-words}.

\subsection{Algorithms and experiment}  

The original motivation for \fullref{thm-main-nonsep}, as well as the
basis of \fullref{conj-mcg-nofiber}, was the results of computer
experiments.  While there is an algorithm which decides if a general
3--manifold fibers using normal surface theory (see Schleimer \cite[Section
6]{Saul:thesis}, Tollefson and Wang \cite{TollefsonWang1996} or
Jaco and Tollefson \cite{JacoTollefson1995}), this
is not practical for all but the smallest examples.  However, special
features allow one to rapidly decide if a tunnel number one
3--manifold fibers over the circle.  In particular, we will show that
it is possible to decide if $M \in \T(r)$ fibers in time which is
polynomial in $\log(r)$.  Our algorithm is important not just for the
experimental side of this paper, but also the theoretical; it forms
the basis for \fullref{thm-main-nonsep}.  Indeed, the basic approach
of the proof is to analyze the algorithm and show that it reports
``does not fiber'' with probability tending to $1$ as $r \to \infty$.
 
Because the fundamental group of a tunnel number one manifold is so
simple, one can use a criterion of Ken Brown to determine if $\pi_1(M)
= \spandef{a,b}{R=1}$ fibers in the above algebraic sense.  Brown's
criterion is remarkably elegant and simple to use, and is in terms of
the combinatorics of the relator $R$.  If $R$ is not in the commutator
subgroup and has length $r$, it takes time $O(r\log r\log \log
r)$.   While that may seem quite fast, some of the manifolds we examined
had $r = 10^{3000}$.  It's not even possible to store the
relator $R$ in this case; after all, the number of elementary
particles in the observable universe is well less than $10^{100}$!

However, in the 3--manifold situation it is possible to specify $R$ by
giving the attaching curve $\gamma$, and $\gamma$ can be described
with only $\log(r)$ bits using either Dehn--Thurston coordinates or
weights on a train track.  Agol, Hass and
W~Thurston described~\cite{AgolHassThurston} how to use splittings of train tracks to compute
certain things about $\gamma$ rapidly, eg,~checking whether $\gamma$
is connected, or computing its homology class.  Motivated by their work,
we were able to adapt Brown's algorithm to work in this setting.  The
resulting algorithm uses train tracks which are labeled by ``boxes''
that remember a small amount of information about a segment of $R$.
As mentioned, it can decide if $M \in \T(r)$ fibers in time polynomial in
$\log(r)$.

\subsection{The tyranny of small examples}\label{subsec-tyranny}

Detailed results of our experiments are given in
\fullref{sec-experiment}, and we will just highlight one aspect
here.  With both notions of a random tunnel number one manifold, it
appears that the probability of fibering goes to $0$ as the complexity
increases; however, the rate at which it converges to $0$ is actually quite
slow, excruciatingly so in the $\MCG$ context.  In particular, most
``small'' manifolds fiber for pretty generous definitions of
``small''.

Starting with the measured lamination notion, the first $r$ for which
$M \in \T(r)$ is \emph{less} likely to fiber than not is about $r =
100{,}000$ (recall here that $r$ is essentially the length of the
relator $R$ in the presentation of $\pi_1(M)$).  The probability of
fibering does not drop below 10\% until about $r = 10^{14}$.

For the mapping class group version, we used the standard 5 Dehn
twists as generators for $\MCG(\partial H)$ (Birman~\cite[Theorem
4.8]{Birman74}).  Recall that the notion here is that given $N$, we
apply a random sequence of $N$ of these Dehn twists to a fixed base
curve $\gamma_0$ to get the attaching curve $\gamma$.  To get the
probability of fibering to be less than 50\%, you need to do $N
\approx 10{,}000$ Dehn twists; to get it below 10\% you need to take
$N \approx 40{,}000$.  It's important to emphasize here how the $\MCG$
notion relates to the measured lamination one, as it's the later that
is related to the size of the presentation for $\pi_1(M)$.  For the
$\MCG$ notion, the length $r$ of the relator experimentally increases
\emph{exponentially} in~$N$.  In particular, doing $N = 10{,}000$ Dehn
twists gives manifolds in $\T(10^{500})$, and $N = 40{,}000$ gives
manifold in $\T(10^{1750})$!

The moral here is that typically the manifolds that one can work with
computationally (eg~with SnapPea \cite{SnapPea}) are so small that it
is not possible to discern the generic behavior from experiments on
that scale alone.  For instance, about 90\% of the cusped manifolds in
the census of Callahan, Hildebrand and Weeks
\cite{CallahanHildebrandWeeks} are fibered (Button \cite{Button2005}), and most
of these manifolds have tunnel number one.  Without the naive version
of Brown's criterion, one would not be able to examine enough
manifolds to suggest \fullref{thm-main-nonsep}; without our improved
train track version, we would not have come to the correct version of
\fullref{conj-mcg-nofiber}.  Indeed, initially we did experiments in
the $\MCG$ case using just the naive version of Brown's algorithm, and
it was clear that the probability of fibering was converging to $1$,
not $0$; this provoked much consternation as to why the ``answer''
differed from the measured lamination case.  Thus one must always keep
an open mind as to the possible generic behavior when examining the
data at hand.

\subsection{More general 3--manifolds}

An obvious question that all of this presents is: What about
\fullref{vague-question} for 3--manifolds which do not have tunnel
number one?  While there are certainly analogous notions of random
manifolds for larger Heegaard genus, closed manifolds, etc., we don't
see the way to any results in that direction.  Unfortunately, the
method we use here is based fundamentally on Brown's criterion, which
is very specific to this case.  Without this tool, it seems daunting
to even try to gather enough experimental evidence to overcome the
skepticism bar set by the discussion in \fullref{subsec-tyranny}.
However, our intuition is that for any Heegaard splitting based
notion of random, the answer would remain unchanged: 3--manifolds
should fiber with probability $0$.  For other types of models, such as
random triangulations, the situation is murkier.  

However, there is one generalization of \fullref{thm-main-nonsep} that
we can do.  Recall that we chose to discuss tunnel number one
3--manifolds where the attaching curve $\gamma$ is non-separating.  If
instead we look at those where $\gamma$ is separating, we get
manifolds~$M$ with two torus boundary components.  In this case
$H^1(M; \Z) =
\Z^2$, and this gives us infinitely many homotopically distinct maps
$M \to S^1$, any one of which could be a fibration.  Thus it is perhaps
surprising that the behavior here is no different than the other case:
\begin{mainforwardsep}
  Let $\T^s$ be the set of tunnel number one manifolds with two
  boundary components.  Then the probability that $M \in \T^s(r)$ fibers
  over the circle goes to $0$ as $r \to \infty$.
\end{mainforwardsep}
If we again compare this result to random 1--relator groups, the behavior is
likely even more divergent than in the non-separating case.  In
particular, we conjecture that for groups $\spandef{a,b}{R = 1}$ where
$R$ is in the commutator subgroup, the probability of algebraically fibering is 1.

\subsection{The Virtual Fibration Conjecture}

As we said at the beginning, one motivation for
\fullref{vague-question} is to provide insight into:
\begin{vfconj}[W~Thurston]
  Let $M$ be an irreducible, atoroidal 3--manifold with infinite
  fundamental group.  Then $M$ has a finite cover which fibers over
  the circle.  
\end{vfconj}
Unlike the basic Virtual Haken Conjecture, which just posits the
existence of a cover containing an incompressible surface, there is
much less evidence for this conjecture.  It has proven quite difficult
to find interesting examples of non-fibered manifolds which can be
shown to virtually fiber, though some infinite classes of tunnel
number one 3--manifolds are known to have this property (Leininger
\cite{Leininger2002}, Walsh \cite{Walsh2005}).  Our work here
certainly suggests that the Virtual Fibration Conjecture is likely to
be much more difficult than the Virtual Haken Conjecture.  One pattern
observed above suggests that the following approach is worth
pondering.  As discussed in \fullref{subsec-tyranny}, ``small''
examples are still quite likely to fiber despite
\fullref{thm-main-nonsep}.  Presuming that this pattern persists for
higher Heegaard genus, one strategy would be to try to find covers
which were ``smaller'' than the initial manifold in some sense.  For
example, one measure of smallness might be the maximum length of a
relator in a minimal genus Heegaard splitting.

\subsection[Dynamical ingredients to the proof \ref{thm-main-nonsep}]{Dynamical ingredients to the proof \fullref{thm-main-nonsep}}

In this introduction, we will not say much about the proof of
\fullref{thm-main-nonsep}.  However, let us at least mention the
two other main ingredients besides our adaptation of Brown's algorithm
in the context of splittings of train tracks.  The first is a theorem
of Mirzakhani~\cite{Mirzakhani2004} which says, in particular, that
non-separating simple closed curves have positive density among all
multicurves (see \fullref{thm-maryam} below).  This lets us sample
simple closed curves by sampling multicurves.  The other concerns
splitting sequences of genus 2 interval exchanges.  For technical
reasons, we actually work with interval exchanges rather than train
tracks, as we indicated above.  Given a measured lamination carried by
an interval exchange $\tau$, we can split it to get a sequence of
exchanges carrying the same lamination.  We prove that genus 2
interval exchanges are normal, that is, any splitting sequence that
can occur does occur for almost all choices of initial measured
lamination (\fullref{thm-dyn-ubiq}).  Our proof is a direct
application of a normality criterion of Kerckhoff~\cite{Kerckhoff1985}.

\subsection{Outline of contents}

\fullref{sec-random-mflds} contains a detailed discussion of the
various notions of random tunnel number one 3--manifolds, and gives
the precise setup for Theorems~\ref{thm-main-nonsep} and
\ref{thm-main-sep}.  \fullref{sec-experiment} gives the
experimental data, which in particular justifies
\fullref{conj-mcg-nofiber}.  \fullref{sec-stallings} just
covers Stallings' theorem, which turns fibering into an algebraic
question.  Then \fullref{subsec-Browns-algorithm} discusses
Brown's algorithm in its original form.
\fullref{sec-random-groups} is about random 2--generator 1--relator
groups as mentioned in \fullref{sub-intro-groups} above.
\fullref{sec-eff-brown} is about train tracks, interval exchanges,
and our efficient version of Brown's algorithm in that setting.  The
rest of the paper is devoted to the proof of
Theorems~\ref{thm-main-nonsep} and \ref{thm-main-sep}.  It starts with
an outline of the main idea in \fullref{sec-main-idea}, where a
proof is given for a indicative toy problem.
\fullref{sec-magic-sequence} is devoted to a certain ``magic
splitting sequence'', which is one of the key tools needed.  We then
prove normality for genus 2 interval exchanges in
\fullref{sec-ubiquity}.  Finally, \fullref{sec-main-thm-pf}
completes the proof by a straightforward assembly of the various
elements.

\subsection{Acknowledgments}

Dunfield was partially supported by the US~National Science
Foundation, both by grant \#DMS-0405491 and as a Postdoctoral Fellow.
He was also supported by a Sloan Fellowship, and some of the work was
done while he was at Harvard University.  Thurston was partially
supported by the US~National Science Foundation as a Postdoctoral
Fellow.  Most of the work was done while he was at Harvard University.
The authors also thank Steve Kerckhoff for helpful conversations and
correspondence, as well as the referee for their very careful reading
of this paper and resulting detailed comments.

\section{Random tunnel number one 3--manifolds}\label{sec-random-mflds}

\subsection{Random 3--manifolds}

What is a ``random 3--manifold''?  Since the set of homeomorphism
classes of compact 3--manifolds is countably infinite, it has no
uniform, countably-additive, probability measure.  However, suppose we
filter the set of 3--manifolds by some notion of complexity where
manifolds of bounded complexity are finite in number.  Then we can
consider limiting probabilities as the complexity goes to infinity.
For instance, we could look at all 3--manifolds which are triangulated
with less than $n$ tetrahedra, and consider the proportion $p_n$ which
are hyperbolic.  If the limit of $p_n$ exists as $n \to \infty$, then it is
a reasonable thing to call the limit the probability that a 3--manifold is
hyperbolic.  Of course, unless the property in question is true for
only finitely many, or all but finitely many, 3--manifolds, the answer
depends on the complexity that we choose.  In other words, it depends
on the \emph{model} of random 3--manifolds.  Nonetheless, if we just pick
one of several natural models to look at, it seems worthwhile to
consider these types of questions to get a better global picture of
the topology of 3--manifolds.  For more on different possible models,
and random 3--manifolds in general, see the work of the first author
and W~Thurston in \cite{DunfieldThurston:random}.  Here, we focus on
the special class of tunnel number one 3--manifolds because it is easy
to determine whether they fiber over the circle.  In the next
subsection, we discuss this class of manifolds, and then give some
natural notions of probability on it.

\subsection{Tunnel number one 3--manifolds}
\label{subsec:def-of-tunnel-num}

Look at an orientable handlebody $H$ of genus $2$.  Consider an
essential simple closed curve $\gamma$ on $\partial H$.  Now one can build a
3--manifold $M$ consisting of $H$ and a 2--handle attached along $\gamma$;
that is, $M = H \cup ( D^2 \times I ) $ where $\partial D^2 \times I$ is glued to $\partial H$
along a regular neighborhood of $\gamma$.  A 3--manifold which can be
constructed in this way is said to have \emph{tunnel number one}.
There are two kinds of tunnel number one 3--manifolds, depending on
whether the attaching curve $\gamma$ separates the surface $\partial H$.  If
$\gamma$ is non-separating, then $\partial M$ is a single torus; if it is
separating, then $\partial M$ is the union of two tori.  When we want to
emphasize the dependence of $M$ on $\gamma$, we will denote it by~$M_\gamma$.

There is a dual description of being tunnel number one, which makes
the origin of the name clear.  Consider a compact orientable
3--manifold $M$ whose boundary is a union of tori.  The manifold $M$
has tunnel number one if and only if there exists an arc $\alpha$ embedded
in $M$, with endpoints on $\partial M$, such that the complement of an open
regular neighborhood of $\alpha$ is a handlebody.  While there are clearly
many 3--manifolds with tunnel number one, it's worth mentioning one
class with which the reader may already be familiar: the exterior of a
2--bridge knot or link in $S^3$.  In this case, the arc in question
joins the top of the two bridges.  In general, 3--manifolds with
tunnel number one are a very tractable class to deal with, and much is
known about them.

\subsection{Measured laminations}\label{subsec-ML-intro}
\label{subsec:measured-lam-prob}

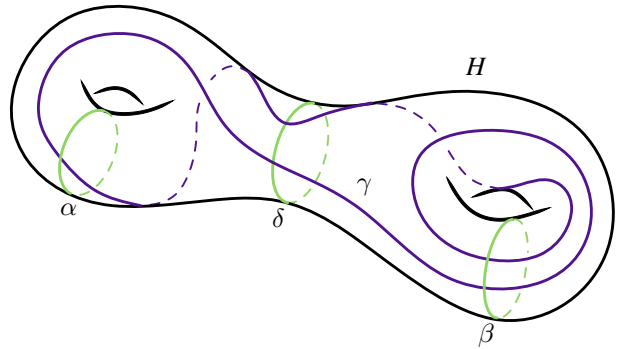
\begin{figure}[htb]\small
\begin{center}
\begin{overpic}[scale=0.72, tics=10]{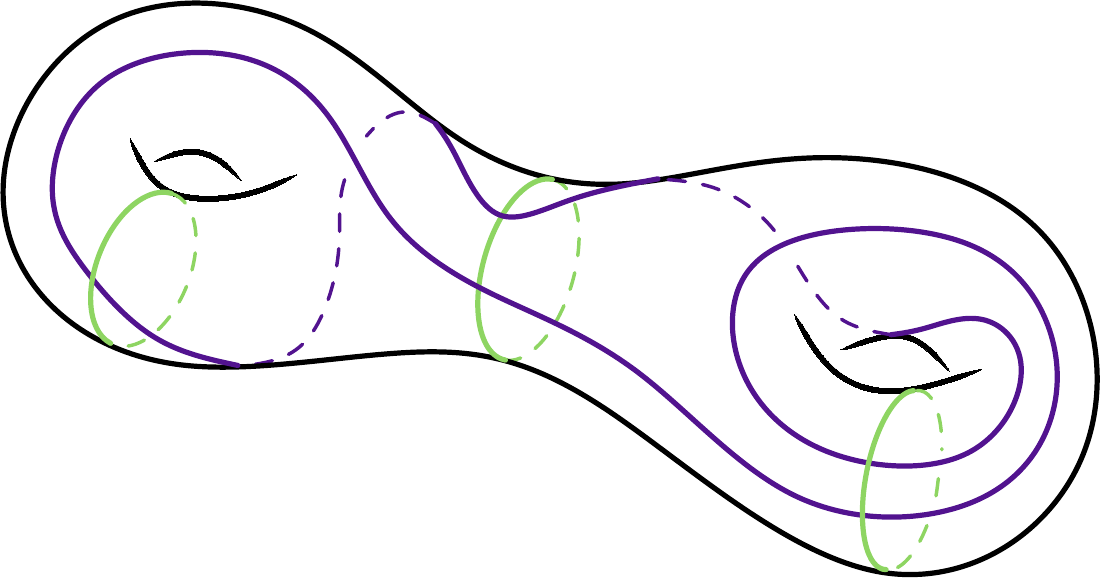}
\put(75,41){$H$}
\put(57,22){$\gamma$}
\put(8,18){$\alpha$}
\put(43, 16){$\delta$}
\put(77,-3){$\beta$}
\end{overpic}
\end{center}
\caption{\
  The curve $\gamma$ has weights $(1, 2, 2)$ and twists $(0, 1, -1)$
  with respect to the Dehn--Thurston coordinates given by the curves
  $(\alpha, \delta, \beta)$.
 }\label{fig-DT-coor}
\end{figure}

Next, we describe our precise parameterization of the tunnel number
one 3--manifolds from the measured laminations point of view.  As
above, let $H$ be a genus 2 handlebody.  Fix a pair of pants
decomposition of $\partial H$ combinatorially equivalent to the curves
$(\alpha, \delta, \beta)$ shown
in \fullref{fig-DT-coor}, so that each of the curves defining this
decomposition bound discs in $H$.

We will use Dehn--Thurston coordinates to parameterize the possible
attaching curves~$\gamma$ for the 2--handle.  A \emph{multicurve} is a
disjoint collection of simple closed curves.  Up to isotopy, any
multicurve $\gamma$ on $\partial H$ is given by weights $w_\alpha$,
$w_\delta$, $w_\beta \in \Z_{\geq0}$ and twist parameters
$\theta_\alpha$, $\theta_\delta$, $\theta_\beta \in \Z$.  Here, the
weights record the (minimal) number of intersections of $\gamma$ with
the curves $(\alpha, \delta, \beta)$, and the twists describe how the
strands of $\gamma$ meet up across the across these curves (with
respect to a certain dual marking).  These coordinates are analogous
to Fenchel--Nielsen coordinates on Teichm\"uller space; see
\fullref{fig-DT-coor} for an example, and Penner and Harer
\cite[Section 1]{PennerHarerBook} or Luo and Stong \cite[Section
2]{LuoStong2004} for details.  When one of the weights is $0$, say
$w_\alpha$, then the twist $\theta_\alpha$ is the number of parallel
copies of $\alpha$ in $\gamma$; thus in this case $\theta_\alpha \geq
0$.  With this convention, Dehn--Thurston coordinates bijectively
parameterize all multicurves, up to isotopy.

Now we let $\T$ be the set of tunnel number one presentations defined by
curves $\gamma$ with the following restrictions with respect to our
choice of Dehn--Thurston coordinates:
\begin{enumerate}
  \item $\gamma$ is a non-separating simple closed curve. \label{en-DT-1}
  \item The weights $w_\alpha, w_\delta, w_\beta$ are $> 0$. \label{en-DT-2}
  \item Each twist~$\theta_i$ satisfies $0 \leq \theta_i < w_i$.  \label{en-DT-3}
  \item $w_\delta \leq \min( 2 w_\alpha , 2 w_\beta )$. \label{en-DT-4}
\end{enumerate}
We now explain why we're making these requirements.  The second
restriction removes some special cases, which are all unknots in lens
spaces and $S^2 \times S^1$; these could be left in without changing the
final theorem as they have asymptotic probability 0.
The third restriction simply accounts for the fact that Dehn twists
along $(\alpha, \delta, \beta)$ extend over $H$; thus any $\gamma$ is equivalent to
one satisfying \eqref{en-DT-3}.

The final restriction serves the following purpose.  If we use the
basis of $\pi_1(H)$ dual to the discs $\alpha, \beta$, then condition
\eqref{en-DT-4} ensures that the word $\gamma$ represents in $\pi_1(H)$
is cyclically reduced.  This is because \eqref{en-DT-4} is the same as
saying that each time $\gamma$ crosses~$\delta$ it then intersects either $\alpha$
or $\beta$ before intersecting $\delta$ again.  It is not immediate that any
$\gamma$ is equivalent, under a homeomorphism of $H$, to one satisfying
\eqref{en-DT-4}; this is the content of
\fullref{lemma-curve-normal-form} below.  The reason that we need to
require (4) is to make the core machinery of the proof work correctly;
as such it is admittedly a tad artificial, and we strongly expect it
is not actually needed.  See \fullref{conj-wo-four} below.

Now, let $\T(r)$ be all elements of $\T$ where $w_\alpha + w_\beta < r$;
equivalently, the ones whose corresponding word in $\pi_1(H)$ has
length $< r$.  As there are only finitely many elements of $\T(r)$, it
makes sense to talk about the probability that they fiber over the
circle.  One form of our main result is:
\begin{theorem}\label{thm-main-nonsep}
  Let $\T$ be the set of tunnel number one manifolds described
  above.  Then the probability that $M \in \T(r)$ fibers over the circle
  goes to $0$ as $r \to \infty$.
\end{theorem}

We can also consider the case of tunnel number one 3--manifolds with
two torus boundary components; these correspond to choosing $\gamma \subset \partial H$
to be a \emph{separating} simple closed curve, replacing
condition~\eqref{en-DT-1}.  We will denote the
corresponding set of manifolds~$\T^s$.  Each $M \in \T^s$ has $H^1(M;
\R) = \R^2$, so they have many chances to fiber over the circle.  Perhaps
surprisingly, the behavior here is no different than the other case:
\begin{theorem}\label{thm-main-sep}
  Let $\T^s$ be the set of tunnel number one manifolds with two
  boundary components.  Then the probability that $M \in \T^s(r)$ fibers
  over the circle goes to $0$ as $r \to \infty$.
\end{theorem}

\subsection{Mapping class group}
\label{subsec:mapping-class-group-prob}

The $\MCG$ model of a random tunnel number one 3--manifold was
completely defined in the introduction.  In this subsection, we
discuss how it differs from the measured lamination version, and what
one would need to leverage \fullref{thm-main-nonsep} into a proof
of \fullref{conj-mcg-nofiber}.  You can skip this section at
first reading, as it's a little technical, and may not make much sense
if you haven't read the proof of \fullref{thm-main-nonsep}.  The
two models for the choice of the attaching curve $\gamma$ can
basically be thought of as different choices of measure on
$\PML(\partial H)$, the space of projectivized measured laminations.  In the
measured lamination model, this measure is just Lebesgue measure on
the sphere $\PML(\partial H)$, whereas for the $\MCG$ model it is a certain
harmonic measure that we describe below.

In proving \fullref{thm-main-nonsep}, we show that there is a
certain open set $U \subset \PML(\partial H)$ so that $\T \cap U$ consists solely
of $\gamma$ so that $M_\gamma$ is not fibered.  The proof then hinges on
showing that $U$ has full Lebesgue measure in the part of $\PML(\partial H)$
defined by requirement \eqref{en-DT-4} above.  The first thing that
one needs to generalize to the $\MCG$ model is to prove a version of
\fullref{thm-main-nonsep} where we drop \eqref{en-DT-4}.  In
particular, for this it suffices to show:
\begin{conjecture}\label{conj-wo-four}
  There exists an open set $V \subset \PML(\partial H)$ of full Lebesgue measure,
  such that for every non-separating curve $\gamma \in V$ the manifold
  $M_\gamma$ is not fibered.  
\end{conjecture}
We are very confident of this conjecture; the proof should be quite
similar to \fullref{thm-main-nonsep} provided certain technical
issues can be overcome.

Assuming this conjecture is true, what one needs to do to prove
\fullref{conj-mcg-nofiber} is show that $V$ has measure 1 with respect
to the following measure.  Fix generators for $\MCG(\partial H)$ and a
base curve $\gamma_0$.  Let $W_n$ be the set of all words in these
generators of length $n$.  Let $\mu_n$ be the probability measure on
$\PML(\partial H)$ which is the average of the point masses supported
at $\phi(\gamma_0)$ for $\phi \in W_n$.  Then we are interested in the
weak limit $\mu$ of these measures, which is called the \emph{harmonic
measure} (Masur and Kaimanovich \cite{Masur_Kaimanovich}).  The
question, then, is whether $\mu(V) \neq 1$.  The relationship between
$\mu$ and Lebesgue messure is not well-understood, and we don't know
how to show that \fullref{conj-mcg-nofiber} follows from
\fullref{conj-wo-four}.  In fact, we suspect that $\mu$ is mutually
discontinuous with Lebesgue measure based in part on
\fullref{sec-experiment}, even though they do agree on $V$.
 
\subsection{Other models}
\label{subsec:triangulation-moves-prob}

There are other ways we could choose $\gamma$ than just the two
detailed above.  For instance, we could start with a one vertex
triangulation of $\partial H$, and then flip edges in a quadrilateral to
obtain a sequence of such triangulations. At the end
of such a sequence of moves, select a edge in the final triangulation
which is a non-separating loop and take that to be~$\gamma$.  Another
approach would be to start with a pair of pants decomposition of
$\partial H$, and then move along a sequence of edges in the pants complex.
Then we would take one of the curves in the final decomposition as
$\gamma$.  We did some haphazard experiments for both these notations as
well, enough to convince us that they also result in a probability of
fibering of 0 and behave generally like the mapping class group
experiments reported above.  For moves in the pair of pants
decomposition, there is a choice of how many Dehn twists to perform on
average before changing the decomposition; as you increase the number
of Dehn twists, the rate at which the manifolds fiber tends to increase.

Strictly speaking, we do not choose our manifolds at random from among
all such manifolds with a given bound on complexity, but rather we
chose from the collection of \emph{descriptions} of bounded
complexity.  These are different as a manifold can have more than one
such description.  Focusing on the measured lamination point of view,
there are two separate issues: first, a manifold can have more than
one unknotting tunnel; second, having fixed an unknotting tunnel,
there may be more than one $\gamma \in \T(r)$ describing it, due to
the action of the mapping class group of the handlebody.  While we do
not prove this here, we believe that, in the measured lamination case,
choosing from descriptions is essentially equivalent to choosing from
among manifolds, as follows.  For the first issue, we strongly believe
that a generic manifold in $\T(r)$ has a unique unknotting tunnel; in
particular, we expect that the distance of the Heegaard splitting
should be very large as $r \to \infty$.  (Another reason why the
number of unknotting tunnels is not a big concern is that this would
only affect our answer if fibered manifolds had many fewer unknotting
tunnels than non-fibered ones.)  About the second issue, namely
multiple descriptions of the same unknotting tunnel, we could further
restrict the conditions (1--4) above on elements of $\T(r)$ to
generically eliminate such multiple descriptions.  As described by
Berge~\cite{BergeDocumentation}, there are simple inequalities in the
weights and twists which ensure that $w_\alpha +w_\beta$ is minimal
among all curves equivalent under the action of the mapping class
group of the handlebody.  This minimal form is typically unique (up to
obvious symmetries, the number of which is independent of the
particular curve at hand). The exception is when there are what
\cite{BergeDocumentation} calls ``level T-transformations''; because
the presence of such transformations is determined by a family of
\emph{equalities}, these occur only in an asymptotically negligible
portion of $\T(r)$.  Thus by supplementing (1--4) we could precisely
parameterize pairs $(M, \mbox{unknotting tunnel})$.  This change would
make no difference in the proof of \fullref{thm-main-nonsep}.

In the case of the mapping class group setup, there is a third issue
which is that there are many random walks in $\MCG$ that end at the
same element.  One could instead work by choosing the elements in
$\MCG$ from larger and larger balls in the Caley graph.
This has two disadvantages.  The first is that in the context of
non-amenable groups such as this one, the study of random walks is
probably more natural than the study of balls; eg, consider the rich
and well-developed theory of the Poisson boundary (Kaimanovich
\cite{KaimanovichSurvey}).  The second is that it is no longer possible
to generate large elements with this alternate distribution, making
experiment impossible (particularly important since experiment is all we
have in this case).  Of course, different elements of $\MCG$ may also
result in the same manifold for the two reasons discussed in the measured
lamination case.  We expect that multiple representatives of the same
curve should be quite rare since the subgroup of $\MCG$ which extends
over the handlebody is very small; in particular, it is of infinite
index.  Indeed, it is easy to show that the probability that a random
walk lies in this subgroup goes to 0 as the length of the random walk
goes to infinity.

\subsection{Curve normal form}  

In this subsection, we justify the claim made in
\fullref{subsec-ML-intro} that given a simple closed curve $\gamma \subset
\partial H$, there is a homeomorphism $\phi$ of the whole handlebody $H$ so
that $\phi(\gamma)$ satisfies condition \eqref{en-DT-4} of
\fullref{subsec-ML-intro}.  Equivalently, we want to find curves
$(\alpha, \delta, \beta)$, arranged as in \fullref{fig-DT-coor}, which satisfy
\begin{equation}\label{eq-normal-form}
  \mbox{Any subarc of $\gamma$ with endpoints in $\delta \cap \gamma$ intersects $\alpha \cup \beta $. } 
\end{equation}
The rest of this section is devoted to:
\begin{lemma}\label{lemma-curve-normal-form}
  Let $\gamma$ be a simple closed curve on the boundary of a genus 2
  handlebody $H$.  Then we can choose $(\alpha, \delta, \beta) \subset \partial H$ bounding
  discs in $H$ as above so that \eqref{eq-normal-form} holds.
\end{lemma}

This lemma is due to Masur~\cite{Masur1986}, and was also described in a
much more general form by Berge~\cite{BergeDocumentation}.  The proof of the
lemma is used in the algorithm for the $\MCG$ case, so as the lemma is
not explicitly set out in \cite{Masur1986}, and
\cite{BergeDocumentation} is unpublished, we include a proof for
completeness.  You can certainly skip it at first reading.

 \begin{figure}[b]\small
    \centering \leavevmode
  \begin{xyoverpic*}{(144,135)}{scale=0.63}{\figdir/n1}
    ,(50,94)*{\alpha}
    ,(126,94)*{\alpha}
    ,(50,17)*{\beta}
    ,(126,17)*{\beta}
    ,(0,0)*{\mbox{(a)}}
  \end{xyoverpic*}
  \hspace{1.0cm}
  \begin{xyoverpic*}{(171,157)}{scale=0.63}{\figdir/n2}
    ,(59,94)*{\alpha}
    ,(135,94)*{\alpha}
    ,(59,17)*{\beta}
    ,(135,17)*{\beta}
    ,(150,137)*+!LD{\delta}
    ,(0,0)*{\mbox{(b)}}
  \end{xyoverpic*}
  \hspace{1cm}
  \caption{}\label{fig-planar-diag}
  \end{figure}

\begin{proof}
  We focus on choosing $\alpha$ and $\beta$ to make the picture as standard as
  possible; the right choice for $\delta$ will then be obvious.  First,
  choose $\alpha$ and $\beta$ to be essential non-separating,
  non-parallel curves that minimize the size of $\gamma \cap (\alpha
  \cup \beta)$.
  Split $H$ open along the discs bounded by $\alpha$ and $\beta$ to get a
  planar diagram as shown in \mbox{\fullref{fig-planar-diag}(a).}
   \begin{figure}\small
    \centering\leavevmode
    \begin{xyoverpic*}{(438,137)}{scale=0.63}{\figdir/n3}
      ,(55,70)*{\alpha}
      ,(131,71)*{\alpha}
      ,(288,70)*{\alpha}
      ,(365,71)*{\beta}
      ,(416,111)*+!LD{\epsilon}
      ,(38,17)*{\mbox{(a)}}
      ,(262,17)*{\mbox{(b)}}
    \end{xyoverpic*}
    \caption{}\label{fig-whitehead-move}
  \end{figure}
  Here the labeled circles, called \emph{vertices}, correspond to the
  discs we cut along, and the arcs are the pieces of $\gamma$.  Note
  vertices with the same label are the endpoints of an equal number of
  arcs, since these endpoints match up when we reglue to get $H$.
  
  We will show that the picture can be made very similar to the one shown
  in \mbox{\fullref{fig-planar-diag}(a);} then the $\delta$ shown in \fullref{fig-planar-diag}(b) works
  to complete the proof.  In particular, it is enough to show:
  \begin{enumerate}
  \item No arc joins a vertex to itself. \label{desiderata-norm-1}
  \item All arcs joining a pair of vertices are isotopic to each
    other in the complement of the other vertices.  \label{desiderata-norm-2}
  \end{enumerate}

  First, suppose we do not have \eqref{desiderata-norm-1}, with $V_0$ being
  the vertex with the bad arc $\gamma_0$.  Consider $V_0 \cup \gamma_0$, which
  separates $S^2$ into two regions.  Both of these regions must
  contain a vertex, or we could isotope $\gamma$ to remove an intersection
  with $\alpha \cup \beta$.  Focus on the component which contains only one
  vertex $V_1$.  If $V_0$ and $V_1$ have the same label as shown in
  \fullref{fig-whitehead-move}(a), we have a contradiction as the
  $V_i$ must be the endpoints of the same number of arcs.
  So we have the situation shown in
  \fullref{fig-whitehead-move}(b).  Replacing $\alpha$ with the
  non-separating curve
  $\epsilon$ indicated reduces $\gamma \cap (\alpha\cup\beta)$, contradicting our initial choice of
  $\alpha$ and $\beta$.
  
  For \eqref{desiderata-norm-2}, there are two basic configurations,
  depending on whether the non-parallel arcs
  join vertices with the same or opposite labels:

\vspace{0.2cm}
{
 \centering\leavevmode\small
 \begin{xyoverpic*}{(163,122)}{scale=0.63}{\figdir/n4}
    ,(15,90)*{\alpha}
    ,(82,90)*{\alpha}
    ,(148,90)*{\beta}
    ,(81,15)*{\beta}
    ,(83,60)*+!U{w}
    ,(83,122)*+!D{z}
    ,(50,93)*+!D{x}
    ,(115,90)*+!D{y}
    ,(31,37)*+!UR{u}
    ,(133,37)*+!UL{v}
  \end{xyoverpic*}
  \hspace{2cm}
 \begin{xyoverpic*}{(163,122)}{scale=0.63}{\figdir/n4}
    ,(15,90)*{\alpha}
    ,(82,90)*{\beta}
    ,(148,90)*{\alpha}
    ,(81,15)*{\beta}
    ,(83,60)*+!U{w}
    ,(83,122)*+!D{z}
    ,(50,93)*+!D{x}
    ,(115,90)*+!D{y}
    ,(31,37)*+!UR{u}
    ,(133,37)*+!UL{v}
  \end{xyoverpic*}

}
\vspace{0.2cm}

  \noindent
  Here, parallel arcs have been drawn as one arc; the label on that arc
  refers to the number of parallel copies (which may be 0).  In the
  case at left, the gluing requirement forces
  \[
  u + w + x + z = x + y \mtext{and} v + w + y + z = u + v,
  \]
  which easily leads to a contradiction.  
  
  In the case at right, we must have $x = y$ and $u = v$ or else we
  can replace $\alpha$ or $\beta$ by a handle slide in the spirit of
  \fullref{fig-whitehead-move}(b) to reduce $\gamma \cap (\alpha \cup  \beta)$.
  Now reglue the $\alpha$ discs to get a solid torus.  Looking at one
  of the $\beta$ vertices, it is joined to $\alpha$ by two families of
  parallel arcs as shown in \fullref{fig-tori-bands}(a).
  \begin{figure}\small
    \centering\leavevmode
    \begin{xyoverpic*}{(163,92)}{scale=0.9}{\figdir/n5}
      ,(60,39)*{\alpha}
      %,(116,34)*[@!90]{\beta}, could rotate beta, but have to choose xydriver
      ,(116,34.5)*{\beta}
      ,(13,5)*{\mbox{(a)}}
    \end{xyoverpic*}
    \hspace{1cm}
    \begin{xyoverpic*}{(163,92)}{scale=0.9}{\figdir/n6}
      ,(13,5)*{\mbox{(b)}}
    \end{xyoverpic*}
    \caption{}\label{fig-tori-bands}
  \end{figure}

  Thinking of this vertex as a bead, slide it along the set of
  parallel $\gamma$ strands past the curve $\alpha$.  Keep sliding past $\alpha$
  in the same direction if possible.  Either:
  \begin{itemize}
  \item This eventually results in an arc joining the pair of $\beta$
    vertices.  In this case, there will also be an arc joining the
    pair of $\alpha$ vertices left over from the final bead slide.
    Because of these two arcs, we can't have non-parallel $\gamma$ arcs
    joining vertices of the same type, which ensures
    \eqref{desiderata-norm-2}.
  \item The bead returns to where it started, so we have something
    like \fullref{fig-tori-bands}(b).  The other $\beta$ vertex must
    be in the same situation, running along a parallel curve on the
    solid torus.  As $\gamma$ is connected, there are no arcs not involved
    in the $\beta$ vertex tracks.  Thus after further sliding, we can
    make the picture completely standard, with the two $\beta$ vertices
    next to each other.  This situation satisfies
    \eqref{desiderata-norm-2} as well.
  \end{itemize}
  Since we have shown $\alpha$ and $\beta$ can be chosen so that
  \eqref{desiderata-norm-1} and \eqref{desiderata-norm-2} hold, we are
  done.
\end{proof}

\section{Experimental results}\label{sec-experiment}

In this section, we give the results of our computer experiments using
the algorithm of \fullref{sec-eff-brown}.  We begin with the
measured lamination notion of random.  For each fixed $r$, we sampled
about $100{,}000$ manifolds $M \in \T(r)$, and used the algorithm to
decide if each one fibers.  Below in \fullref{fig-ML-plots} are
the results for various $r \leq 10^{20}$.

While these results are superseded by \fullref{thm-main-nonsep},
there are still interesting things to notice about the plots.  For
instance, look at the \emph{rate} at which the probability of fibering
approaches $0$; as we already discussed in
\fullref{subsec-tyranny}, it is quite leisurely.  Moreover, the
convergence has a very specific form --- as the log-log plot in the
bottommost part of
\fullref{fig-ML-plots} makes clear, it converges to $0$ like $c_1
e^{-c_2 r}$ for some positive constants $c_i$.  In the proof of
\fullref{thm-main-nonsep}, we will see why this should be the
case.

% The following variables control the size of the pictures
\newcommand{\plotgraph}[1]{\includegraphics[scale=0.72]{#1}}  %%was .8
\newcommand{\plotpicscale}{0.0144bp} % for scale = 1.0, should be 0.020, was .016

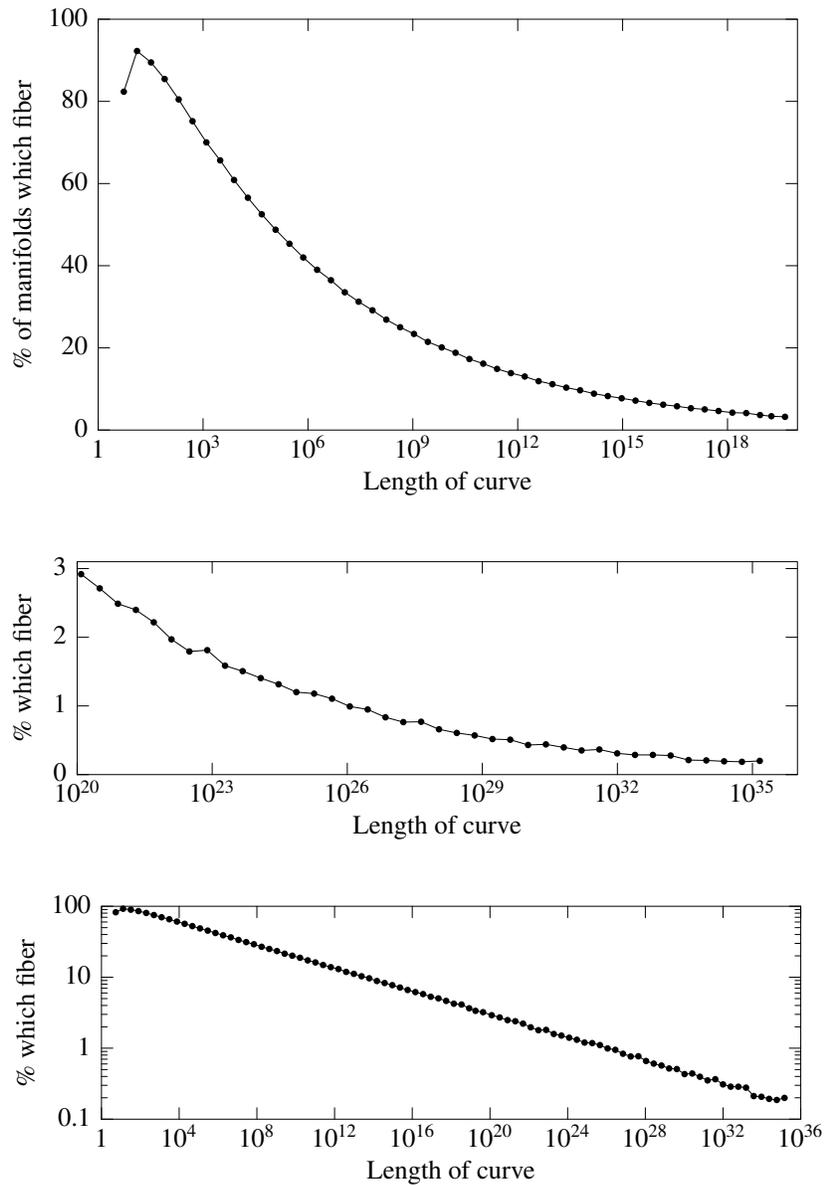
\begin{figure}\small
\begin{center}
  %GNUPLOT: LaTeX picture with Postscript
\begin{picture}(0,0)%
\plotgraph{\figdir/lamination1-pdf}%
\end{picture}%
\begingroup
\setlength{\unitlength}{\plotpicscale}%
\begin{picture}(21600,12960)(0,0)%
\put(2200,1650){\makebox(0,0)[r]{\strut{} 0}}%
\put(2200,3802){\makebox(0,0)[r]{\strut{} 20}}%
\put(2200,5954){\makebox(0,0)[r]{\strut{} 40}}%
\put(2200,8106){\makebox(0,0)[r]{\strut{} 60}}%
\put(2200,10258){\makebox(0,0)[r]{\strut{} 80}}%
\put(2200,12410){\makebox(0,0)[r]{\strut{} 100}}%
\put(18945,1100){\makebox(0,0){\strut{}$10^{18}$}}%
\put(16200,1100){\makebox(0,0){\strut{}$10^{15}$}}%
\put(13455,1100){\makebox(0,0){\strut{}$10^{12}$}}%
\put(10710,1100){\makebox(0,0){\strut{}$10^{9}$}}%
\put(7965,1100){\makebox(0,0){\strut{}$10^{6}$}}%
\put(5220,1100){\makebox(0,0){\strut{}$10^{3}$}}%
\put(2475,1100){\makebox(0,0){\strut{}$1$}}%
\put(550,7030){\rotatebox{90}{\makebox(0,0){\strut{}\% of manifolds which fiber}}}%
\put(11625,275){\makebox(0,0){\strut{}Length of curve}}%
\end{picture}%
\endgroup
 
\vspace{6mm}
  %GNUPLOT: LaTeX picture with Postscript
\begin{picture}(0,0)%
\plotgraph{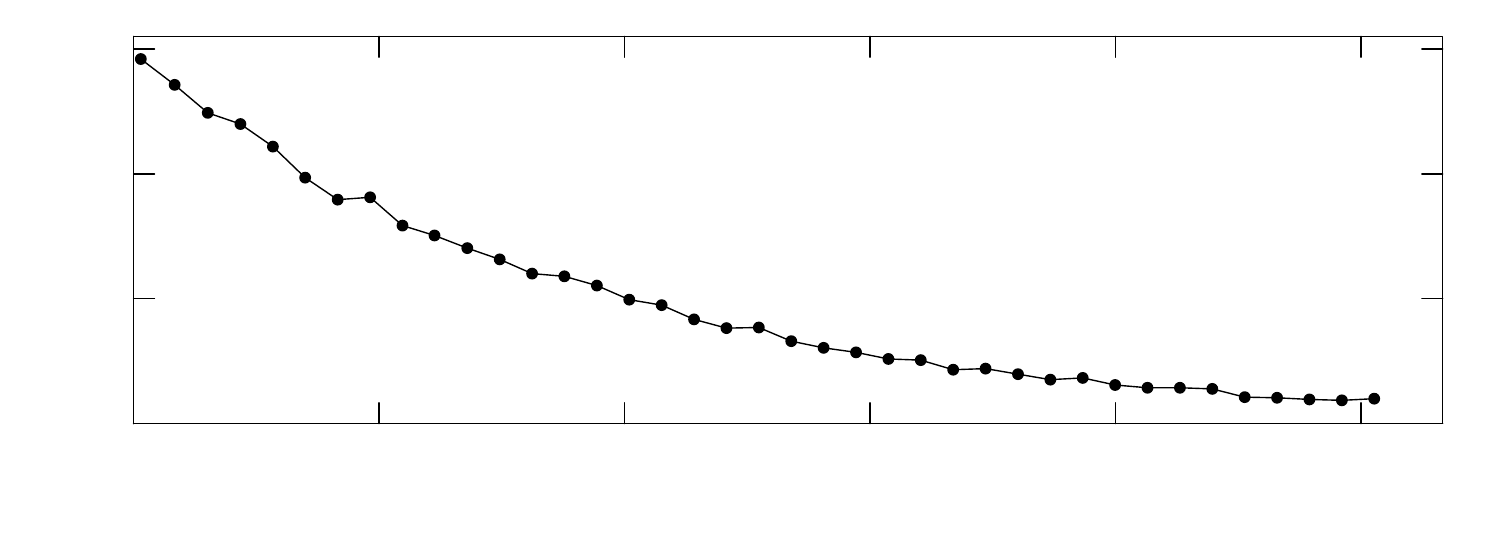}%
\end{picture}%
\begingroup
\setlength{\unitlength}{\plotpicscale}%
\begin{picture}(21600,7776)(0,0)%
\put(1650,1650){\makebox(0,0)[r]{\strut{} 0}}%
\put(1650,3449){\makebox(0,0)[r]{\strut{} 1}}%
\put(1650,5247){\makebox(0,0)[r]{\strut{} 2}}%
\put(1650,7046){\makebox(0,0)[r]{\strut{} 3}}%
\put(19597,1100){\makebox(0,0){\strut{}$10^{35}$}}%
\put(16063,1100){\makebox(0,0){\strut{}$10^{32}$}}%
\put(12528,1100){\makebox(0,0){\strut{}$10^{29}$}}%
\put(8994,1100){\makebox(0,0){\strut{}$10^{26}$}}%
\put(5459,1100){\makebox(0,0){\strut{}$10^{23}$}}%
\put(1925,1100){\makebox(0,0){\strut{}$10^{20}$}}%
\put(550,4438){\rotatebox{90}{\makebox(0,0){\strut{}\% which fiber}}}%
\put(11350,275){\makebox(0,0){\strut{}Length of curve}}%
\end{picture}%
\endgroup
 
\vspace{6mm}
  %GNUPLOT: LaTeX picture with Postscript
\begin{picture}(0,0)%
\plotgraph{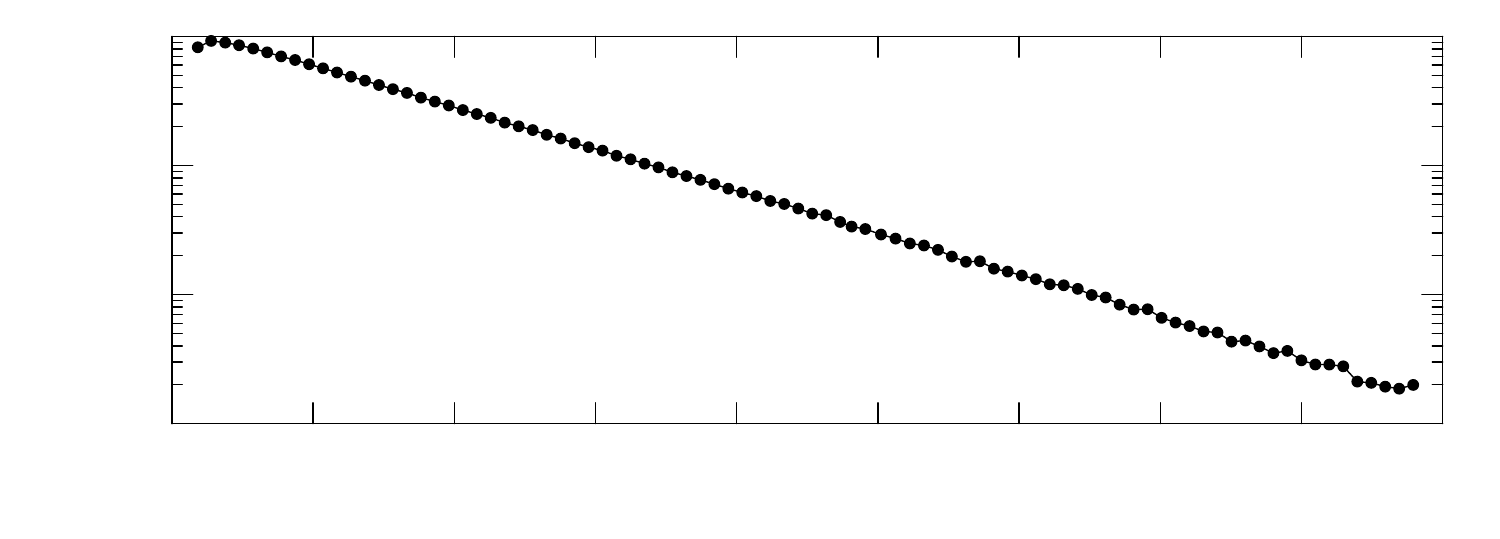}%
\end{picture}%
\begingroup
\setlength{\unitlength}{\plotpicscale}%
\begin{picture}(21600,7776)(0,0)%
\put(2200,1650){\makebox(0,0)[r]{\strut{} 0.1}}%
\put(2200,3509){\makebox(0,0)[r]{\strut{} 1}}%
\put(2200,5367){\makebox(0,0)[r]{\strut{} 10}}%
\put(2200,7226){\makebox(0,0)[r]{\strut{} 100}}%
\put(20775,1100){\makebox(0,0){\strut{}$10^{36}$}}%
\put(18742,1100){\makebox(0,0){\strut{}$10^{32}$}}%
\put(16708,1100){\makebox(0,0){\strut{}$10^{28}$}}%
\put(14675,1100){\makebox(0,0){\strut{}$10^{24}$}}%
\put(12642,1100){\makebox(0,0){\strut{}$10^{20}$}}%
\put(10608,1100){\makebox(0,0){\strut{}$10^{16}$}}%
\put(8575,1100){\makebox(0,0){\strut{}$10^{12}$}}%
\put(6542,1100){\makebox(0,0){\strut{}$10^{8}$}}%
\put(4508,1100){\makebox(0,0){\strut{}$10^{4}$}}%
\put(2475,1100){\makebox(0,0){\strut{}$1$}}%
\put(550,4438){\rotatebox{90}{\makebox(0,0){\strut{}\% which fiber}}}%
\put(11625,275){\makebox(0,0){\strut{}Length of curve}}%
\end{picture}%
\endgroup
 

\caption{
  Data for the probability of fibering from the measured lamination
  point of view.  The horizontal axis is the size $r$ of the curve $\gamma$ in
  Dehn--Thurston coordinates, or equivalently the length of the relator
  in the resulting presentation of $\pi_1$. Each point represents a
  sample of about 100,000 manifolds. 
}\label{fig-ML-plots}

\end{center}
\end{figure}

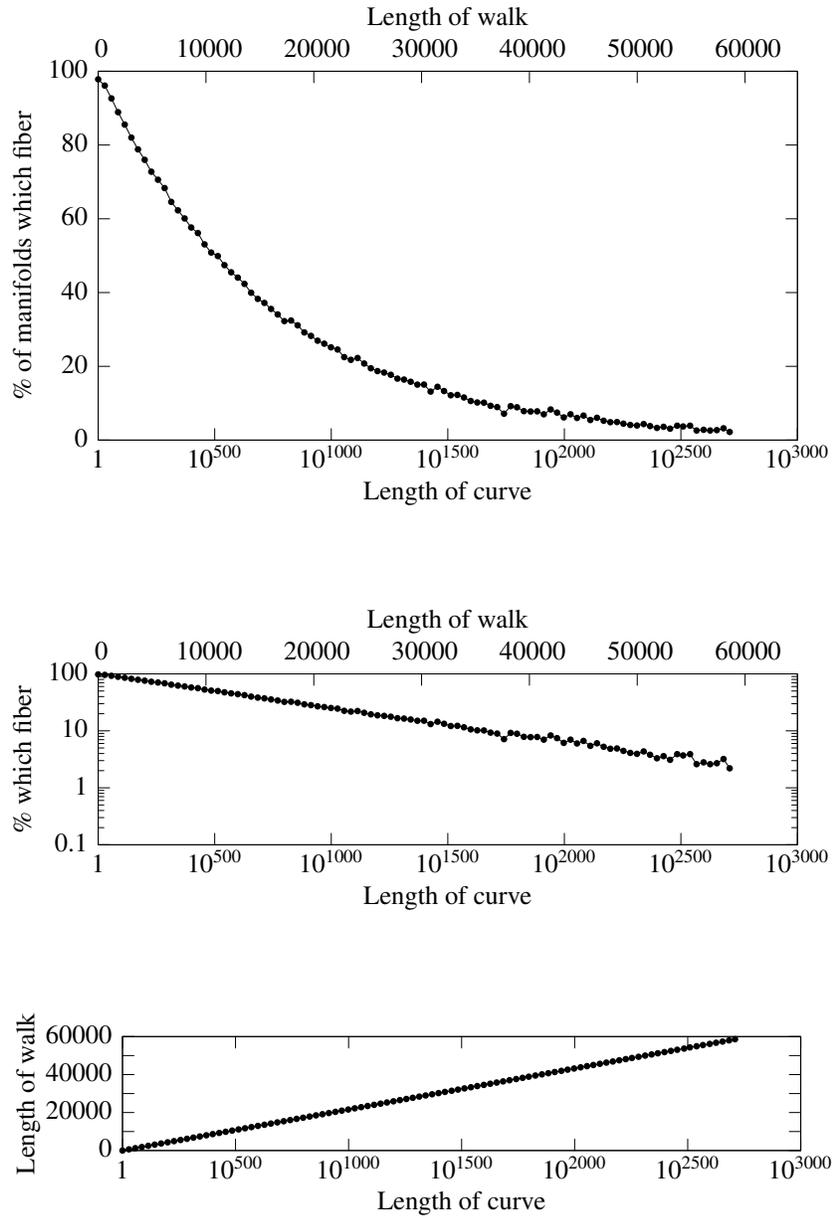
\begin{figure}\small
\begin{center}
%GNUPLOT: LaTeX picture with Postscript
\begin{picture}(0,0)%
\plotgraph{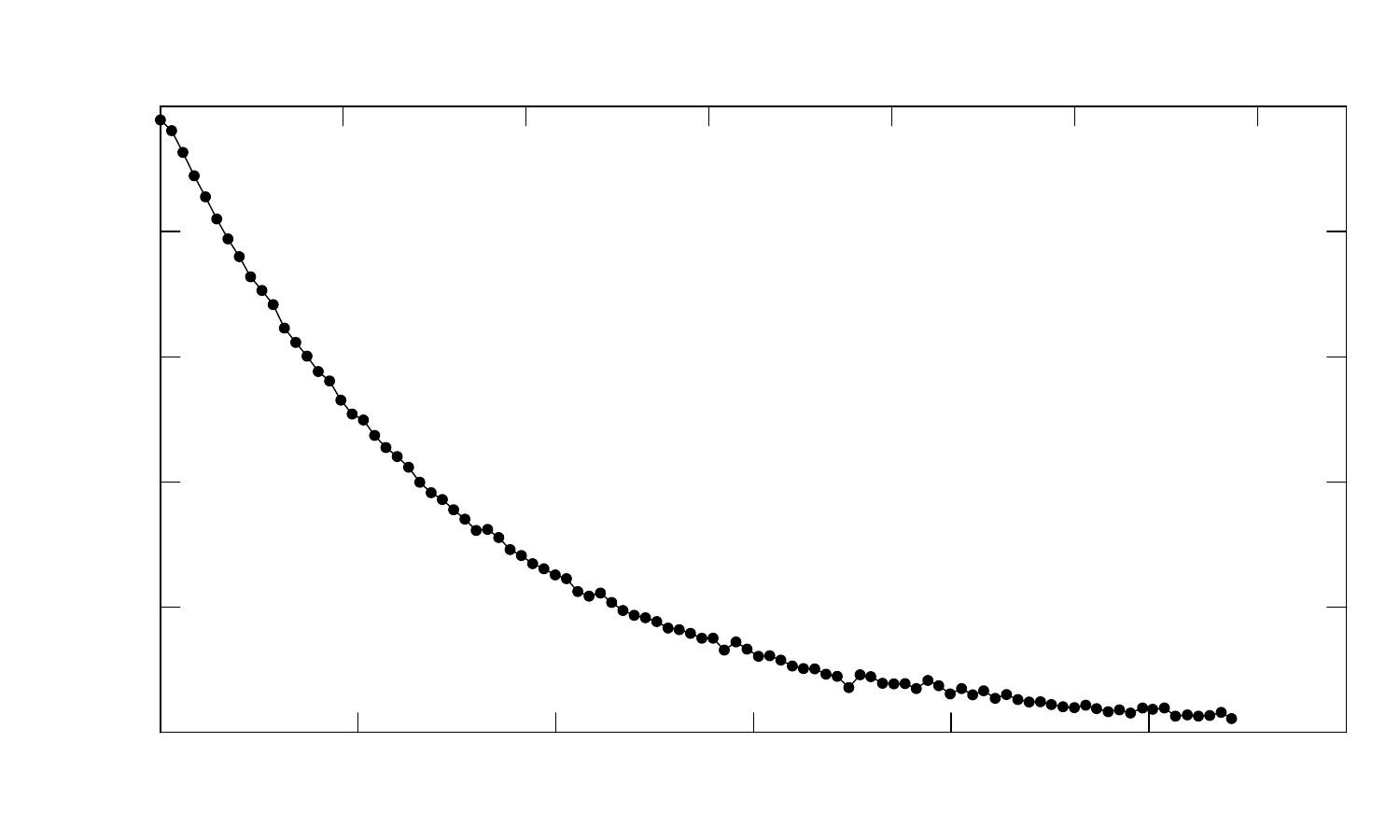}%
\end{picture}%
\begingroup
\setlength{\unitlength}{\plotpicscale}%
\begin{picture}(21600,12960)(0,0)%
\put(2200,1650){\makebox(0,0)[r]{\strut{} 0}}%
\put(2200,3582){\makebox(0,0)[r]{\strut{} 20}}%
\put(2200,5514){\makebox(0,0)[r]{\strut{} 40}}%
\put(2200,7446){\makebox(0,0)[r]{\strut{} 60}}%
\put(2200,9378){\makebox(0,0)[r]{\strut{} 80}}%
\put(2200,11310){\makebox(0,0)[r]{\strut{} 100}}%
\put(20775,1100){\makebox(0,0){\strut{}$10^{3000}$}}%
\put(17725,1100){\makebox(0,0){\strut{}$10^{2500}$}}%
\put(14675,1100){\makebox(0,0){\strut{}$10^{2000}$}}%
\put(11625,1100){\makebox(0,0){\strut{}$10^{1500}$}}%
\put(8575,1100){\makebox(0,0){\strut{}$10^{1000}$}}%
\put(5525,1100){\makebox(0,0){\strut{}$10^{500}$}}%
\put(2475,1100){\makebox(0,0){\strut{}$1$}}%
\put(2475,11860){\makebox(0,0){\strut{} 0}}%
\put(5296,11860){\makebox(0,0){\strut{} 10000}}%
\put(8118,11860){\makebox(0,0){\strut{} 20000}}%
\put(10939,11860){\makebox(0,0){\strut{} 30000}}%
\put(13761,11860){\makebox(0,0){\strut{} 40000}}%
\put(16582,11860){\makebox(0,0){\strut{} 50000}}%
\put(19404,11860){\makebox(0,0){\strut{} 60000}}%
\put(550,6480){\rotatebox{90}{\makebox(0,0){\strut{}\% of manifolds which fiber}}}%
\put(11625,275){\makebox(0,0){\strut{}Length of curve}}%
\put(11625,12684){\makebox(0,0){\strut{}Length of walk}}%
\end{picture}%
\endgroup
 
\vspace{1.4cm}
%GNUPLOT: LaTeX picture with Postscript
\begin{picture}(0,0)%
\plotgraph{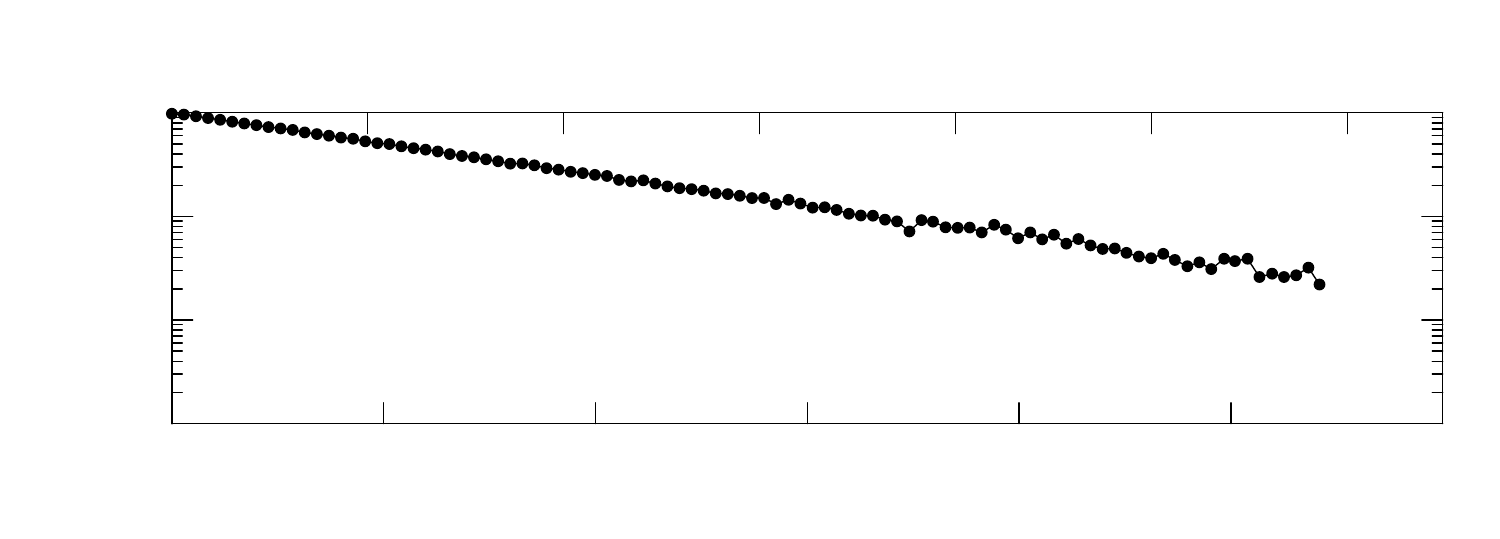}%
\end{picture}%
\begingroup
\setlength{\unitlength}{\plotpicscale}%
\begin{picture}(21600,7776)(0,0)%
\put(2200,1650){\makebox(0,0)[r]{\strut{} 0.1}}%
\put(2200,3142){\makebox(0,0)[r]{\strut{} 1}}%
\put(2200,4634){\makebox(0,0)[r]{\strut{} 10}}%
\put(2200,6126){\makebox(0,0)[r]{\strut{} 100}}%
\put(20775,1100){\makebox(0,0){\strut{}$10^{3000}$}}%
\put(17725,1100){\makebox(0,0){\strut{}$10^{2500}$}}%
\put(14675,1100){\makebox(0,0){\strut{}$10^{2000}$}}%
\put(11625,1100){\makebox(0,0){\strut{}$10^{1500}$}}%
\put(8575,1100){\makebox(0,0){\strut{}$10^{1000}$}}%
\put(5525,1100){\makebox(0,0){\strut{}$10^{500}$}}%
\put(2475,1100){\makebox(0,0){\strut{}$1$}}%
\put(2475,6676){\makebox(0,0){\strut{} 0}}%
\put(5296,6676){\makebox(0,0){\strut{} 10000}}%
\put(8118,6676){\makebox(0,0){\strut{} 20000}}%
\put(10939,6676){\makebox(0,0){\strut{} 30000}}%
\put(13761,6676){\makebox(0,0){\strut{} 40000}}%
\put(16582,6676){\makebox(0,0){\strut{} 50000}}%
\put(19404,6676){\makebox(0,0){\strut{} 60000}}%
\put(550,3888){\rotatebox{90}{\makebox(0,0){\strut{}\% which fiber}}}%
\put(11625,275){\makebox(0,0){\strut{}Length of curve}}%
\put(11625,7500){\makebox(0,0){\strut{}Length of walk}}%
\end{picture}%
\endgroup
 
\vspace{1.4cm}
%GNUPLOT: LaTeX picture with Postscript
\begin{picture}(0,0)%
\plotgraph{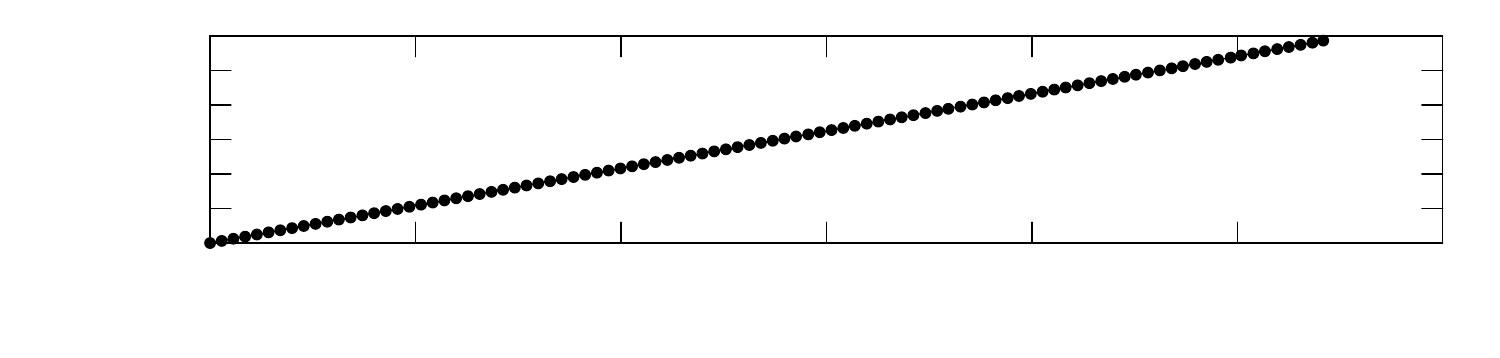}%
\end{picture}%
\begingroup
\setlength{\unitlength}{\plotpicscale}%
\begin{picture}(21600,5184)(0,0)%
\put(2750,1650){\makebox(0,0)[r]{\strut{} 0}}%
%\put(2750,2147){\makebox(0,0)[r]{\strut{} 10000}}%
\put(2750,2645){\makebox(0,0)[r]{\strut{} 20000}}%
%\put(2750,3142){\makebox(0,0)[r]{\strut{} 30000}}%
\put(2750,3639){\makebox(0,0)[r]{\strut{} 40000}}%
%\put(2750,4137){\makebox(0,0)[r]{\strut{} 50000}}%
\put(2750,4634){\makebox(0,0)[r]{\strut{} 60000}}%
\put(20775,1100){\makebox(0,0){\strut{}$10^{3000}$}}%
\put(17817,1100){\makebox(0,0){\strut{}$10^{2500}$}}%
\put(14858,1100){\makebox(0,0){\strut{}$10^{2000}$}}%
\put(11900,1100){\makebox(0,0){\strut{}$10^{1500}$}}%
\put(8942,1100){\makebox(0,0){\strut{}$10^{1000}$}}%
\put(5983,1100){\makebox(0,0){\strut{}$10^{500}$}}%
\put(3025,1100){\makebox(0,0){\strut{}$1$}}%
\put(550,3142){\rotatebox{90}{\makebox(0,0){\strut{}Length of walk}}}%
\put(11900,275){\makebox(0,0){\strut{}Length of curve}}%
\end{picture}%
\endgroup
 

\caption{
  Data for the $\MCG$ case.  All the points represent samples of at
  least 1000 manifolds.  The first half or so of the data represent
  10,000 manifolds. 
}\label{fig-MCG-plots}

\end{center}
\end{figure}

Before moving on to the $\MCG$ case, let us make one quick comment on
why we can easily sample $M \in \T(r)$ uniformly
at random.  While it is easy to pick a random multicurve with $w_\alpha + w_\beta
\approx r$, a priori there is no way to ensure that we sample only
connected non-separating curves.  Fortunately, Mirzakhani has shown
that there is a definite probability, roughly 1/5, that a randomly
chosen multicurve is of this form (see \fullref{thm-maryam}
below).  Thus one simply samples multicurves at random, ignoring all
of those which are not of the desired form.

We turn now to the $\MCG$ case.  For this, we choose the standard five
Dehn twists as our generating set for $\MCG(\partial H)$.  The results are
shown in \fullref{fig-MCG-plots}.  There are two horizontal scales
on each of the upper two plots.  Along the top is the number of Dehn
twists done to create the manifold, ie,~the length of the walk in
$\MCG(\partial H)$.  To give a scale at which to compare it to the previous
figure, along the bottom is the size of the resulting attaching curve
$\gamma$ in terms of the standard Dehn--Thurston coordinates. As the plot
at the bottom shows, the Dehn--Thurston size grows exponentially in
length of the walk, which justifies the use of the two scales on the
upper graphs.

One thing to notice here is just how slowly the probability goes to
zero in terms of the Dehn--Thurston size; in the earlier
\fullref{fig-ML-plots}, the probability of fibering was less than
0.3\% for $r = 10^{35}$, but here the probability is still greater
than 40\% at $r = 10^{500}$.  This reinforces the point made in
\fullref{subsec:mapping-class-group-prob} that the $\gamma$ resulting
from the $\MCG$ process are not generic with respect to Lebesgue
measure on $\PML(\partial H)$.

Because of how large some of these curves are, we had to use much
smaller samples than in the earlier case; this is why the graph looks
so jumpy.  However, if we look at the middle plot, we again see near
perfect exponential decay, just as in the measure lamination case.
Thus we are quite confident that \fullref{conj-mcg-nofiber} is
correct.

\subsection{Fibering in slices of $\PML(\partial H)$}

\begin{figure}[ht!]
\centering
\includegraphics[width=2.9in]{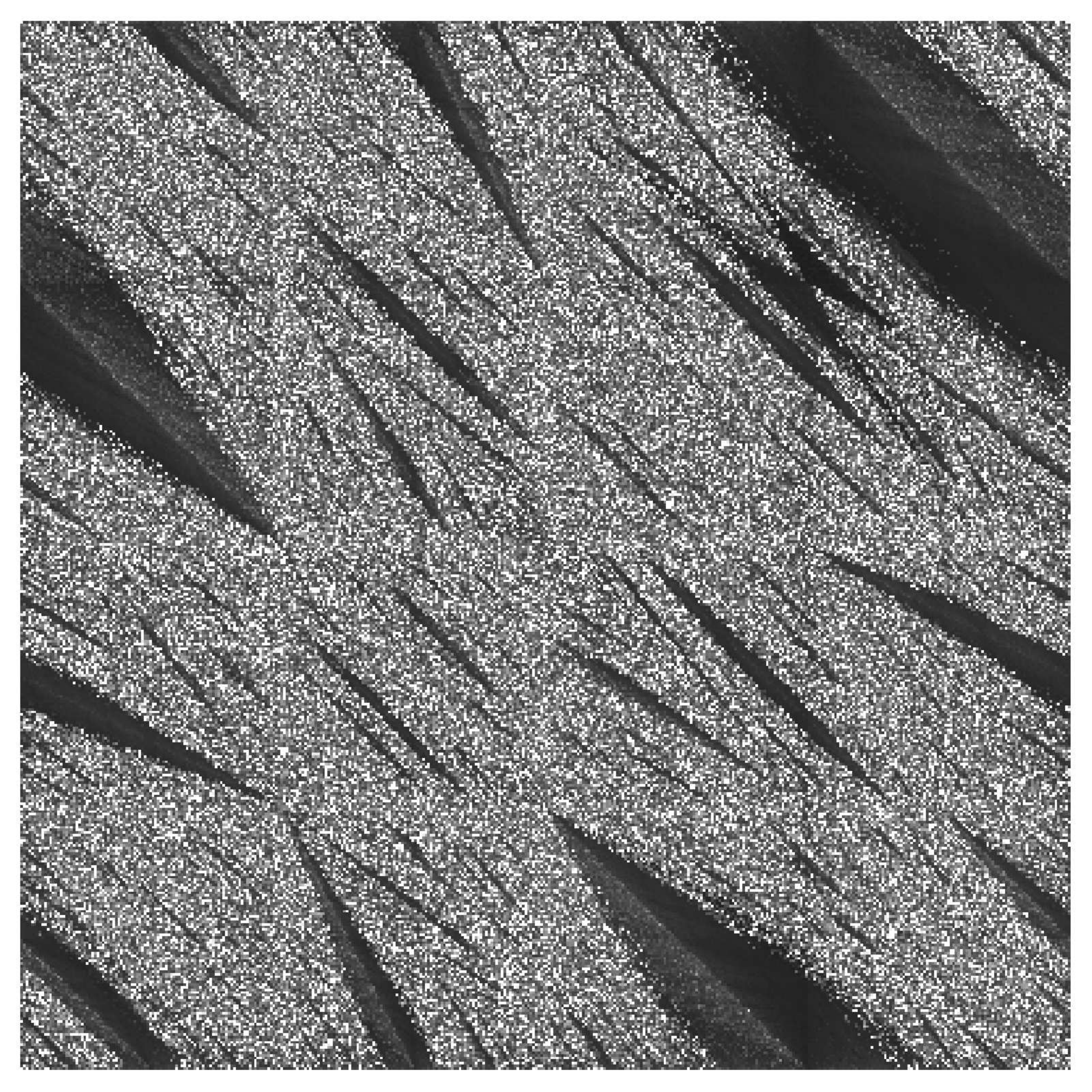}   
\caption{
  A slice of $\PML(\partial H; \R)$, with curves where $M_\gamma$ is fibered
  plotted in white, and where $M_\gamma$ is non-fibered plotted in
  shades of gray; darker grays indicate that the program
  reported ``non-fibered'' in fewer steps.  }\label{fig-slice}
\end{figure}

The parameter space $\T$ of tunnel number one 3--manifolds is a subset
of $\ML(\partial H; \Z)$.  Let us projectivize, and so view $\T$ as a subset
of $\PML(\partial H; \R)$, which is just the 5--sphere.  If we take a two
dimensional projectively linear slice of $\PML(\partial H; \R)$, we can plot
the fibered points of $\T$ in the following sense.  Fix some positive
number $r$.  Divide the slice into little boxes, and in each box pick
a random $\gamma \in \T(r)$ and plot whether or not $M_\gamma$ fibers.  Of
course, as $r \to \infty$ the probability of fibering goes to zero, so it
is much more informative to plot \emph{how many steps} the algorithm
takes before it reports ``not fibered''.  \fullref{fig-slice} shows
the results for one such slice, where we fixed $w_\alpha \approx w_\beta \approx (2/3)
w_\delta \approx 2 \theta_\delta$, and took $r \approx 10^{11}$.  The horizontal and vertical
axes are $\theta_\alpha$ and $\theta_\beta$; since they are well defined modulo $w_a$ and
$w_b$, the figure should be interpreted as living on the torus.

\subsection{Knots in $S^3$}

\begin{figure}[ht!]
\small
\centering
\import{\figdir/}{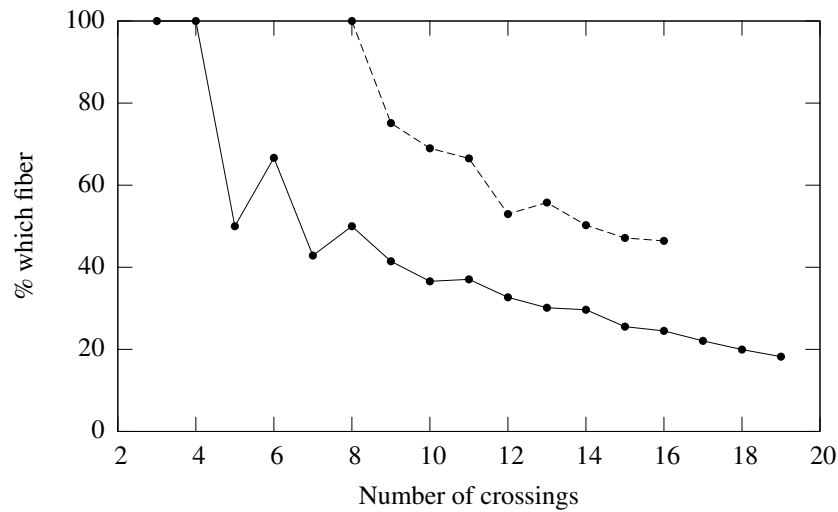}
\caption{
  Proportion of prime knots whose complements fiber.  The lower line
  plots alternating knots, and the upper one non-alternating knots.  The
  vertical axis is really the proportion with non-monic Alexander
  polynomial; as such it is an \emph{upper} bound on the proportion
  that fiber.  }\label{fig-knot-plot}
\end{figure}

As mentioned in the introduction, we suspect that the pattern
exhibited above should persist for any Heegaard splitting based model
of random manifold, regardless of genus.  It is less clear what would
happen for, say, random triangulations.  We did a little experiment
for knots in $S^3$, where one filters (isomorphism classes of) prime
knots by the number of crossings.  Rather than address the difficult
question of whether they fiber, we looked instead at whether the lead
coefficient of the Alexander polynomial is $\pm1$, ie,~the polynomial
is monic.  A monic Alexander polynomial is necessary for fibering, but
not sufficient.  For alternating knots, however, it is sufficient
(Murasugi \cite{Murasugi1963}), and for non-alternating knots with few
crossings there are probably not many non-fibered knots with monic
Alexander polynomials.  The results are shown in
\mbox{\fullref{fig-knot-plot}}.  We used the program Knotscape
\cite{Knotscape} with knot data from Hoste, Thistlethwaite and Weeks
\cite{KnotscapePaper} and Rankin, Flint and Schermann
\cite{RankinFlintSchermannI,RankinFlintSchermannII}.  In light of
\fullref{subsec-tyranny}, we do not wish to draw any conclusions from
this data.  Really, what needs to be done is to figure out how to
generate a random prime knot with, say, 100 crossings with close to
the uniform distribution.  It would be quite interesting to do so even
for alternating knots; for this case, a place to start might be
Poulalhon, Schaeffer and Zinn-Justin
\mbox{\cite{PoulalhonSchaeffer2003,SchaefferZinn-Justin}.}

% The sizing for this figure is controlled in plot/knots.tex:

\subsection{Implementation notes}

The algorithm used in the experiments is the one described in
\fullref{sec-eff-brown}, though some corners were cut in the
implementation of the Dehn twist move; so our resulting program is not
completely efficient in the sense of that section.  In the $\MCG$
case, we had to deal with the fact that the resulting attaching curve
$\gamma$ might not be in $\T$; that is, it fails to satisfy conditions
(2--4) of \fullref{subsec-ML-intro}.  To rectify this, we used the
method of the proof of \fullref{lemma-curve-normal-form} to get
an equivalent curve in $\T$.  The complete source code for our program
can downloaded from the front page for this paper: 
\href{http://dx.doi.org/10.2140/gt.2006.10.2431}{{\tt 
DOI:10.2140/gt.2006.10.2431}}

\section{Algebraic criterion for fibering}\label{sec-stallings}

In the next two sections, we describe how to determine if a tunnel
number one 3--manifold fibers over the circle.  There is an
exponential-time algorithm from normal surface
theory for deciding such questions in general~\cite{Saul:thesis}, but that is
impractical for our purposes.  The criterion we give below is
purely combinatorial in terms of the word in the fundamental group of
the handlebody given by the attaching curve of the 2--handle.  In this
section, we give a theorem of Stallings which reduces this geometric
question to an algebraic one about the fundamental group.  In the next
section, we describe an algorithm of K.~Brown which then completely
solves the algebraic question for the special type of groups coming
from tunnel number one 3--manifolds.

\subsection{Stallings' Theorem} 

Suppose a 3--manifold $M$ fibers over the circle:
\[
F \longrightarrow M \stackrel{\phi}{\longrightarrow} S^1
\]
The map $\phi \maps M \to S^1$ represents an element of $H^1(M ; \Z)$,
which is Poincar\'e dual to the element of $H_2(M ; \Z)$ represented by
the fiber $F$.  Now, take $M$ to be any compact 3--manifold.
Continuing to think of a $\left[ \phi \right] \in H^1(M ; \Z)$ as a
homotopy class of maps $M \to S^1$, it makes sense to ask if $\left[ \phi
\right]$ can be represented by a fibration over $S^1$.  
Associated to $\left[ \phi \right] \in H^1(M ; \Z)$ is the infinite
cyclic cover $\tilde{M}$ of $M$ whose fundamental group is the kernel
of $\phi_* \maps \pi_1(M) \to \pi_1(S^1) = \Z$.  If $\left[ \phi \right]$ can
be represented by a fibration, then $\tilde{M}$ is just
$(\mbox{fiber}) \times \R$.  In particular, the kernel of $\phi_*$ is
$\pi_1(\mbox{fiber})$, and hence finitely generated.  The converse to
this is also true:
\begin{theorem}[Stallings \cite{Stallings62}]\label{thm-stallings}
  Let $M$ be a compact, orientable, irreducible 3--manifold.  Consider
  a $\left[ \phi \right] \neq 0$ in $H^1(M;\Z)$.  Then $\left[ \phi \right]$
  can be represented by a fibration if and only if the kernel of $\phi_*
  \maps \pi_1(M) \to \Z$ is finitely generated.
\end{theorem}
The irreducibility hypothesis here is just to avoid the Poincar\'e
Conjecture; it rules out the possibility that $M$ is the connect sum
of a fibered 3--manifold and a nontrivial homotopy sphere.  When $M$
has tunnel number one, the irreducibility hypothesis can be easily
dropped without presuming the Poincar\'e Conjecture, as follows.
Consider $M$ as a genus 2 handlebody $H$ with a 2--handle attached
along $\gamma \subset \partial H$.  By Jaco's Handle Addition Lemma
(Jaco \cite{Jaco82}, Scharlemann \cite{Scharlemann85}), $M$ is
irreducible if $\partial H \setminus \gamma$ is incompressible in $H$.
If instead $\partial H \setminus \gamma$ compresses, then it is not
hard to see that $M$ is the connected sum of a lens space with the
exterior of the unknot in the 3--ball.  As lens spaces trivially
satisfy the Poincar\'e Conjecture, we have:
\begin{corollary}\label{cor:tunnel-one-stallings}
  Let $M^3$ have tunnel number one.  Then $M$ fibers over the circle
  if and only if there exists a $\left[\phi \right] \neq 0$ in   $H^1(M,
\Z)$ such that the kernel of $\phi_* \maps M \to \Z$ is finitely
generated.
\end{corollary}

In general, if $G$ is a finitely presented group and $G \to \Z$ an
epimorphism, deciding if the kernel is finitely generated is a very
difficult question.  Note that if $H$ is a finitely presented group,
then $H$ is trivial if and only if the obvious epimorphism $H * \Z \to
\Z$ has finitely generated kernel.  Thus our question subsumes the
problem of deciding if a given $H$ is trivial, and hence is
algorithmically undecidable.  Thus, it is not at all clear that
Stallings' Theorem can be leveraged to an  algorithm to decide if a
3--manifold fibers.  However, as we'll see in the next section, the
algebraic problem is solvable in the case of a presentation with two
generators and one relation, giving us a practical algorithm to decide if a
tunnel number one 3--manifold fibers over the circle.

\section{Brown's Algorithm}\label{subsec-Browns-algorithm}

Consider a two-generator, one-relator group $G = \spandef{a,b}{R =
  1}$.  Given an epimorphism $\phi \maps G \to \Z$, Kenneth Brown gave an
elegant algorithm which decides if the kernel of $\phi$ is finitely
generated \cite{Brown87}.  Brown was interested in computing the
Bieri--Neumann--Strebel (BNS) invariant of $G$, which is closely related
to this question.  We will first discuss Brown's algorithm for a fixed
$\phi$, and then move to the BNS context to understand what happens
for all $\phi$ at once.

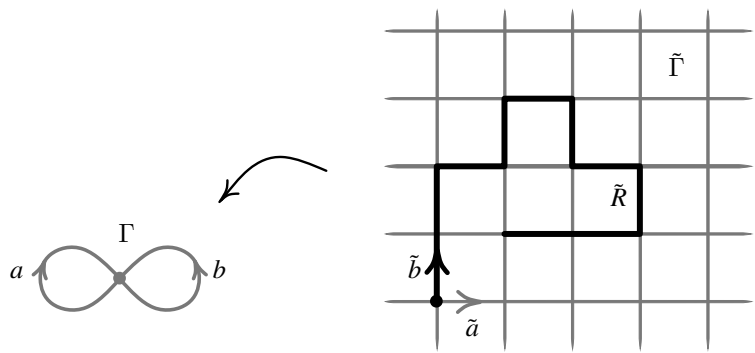
\begin{figure}[b]\small
\centering
\begin{overpic}[scale=0.9, tics=5]{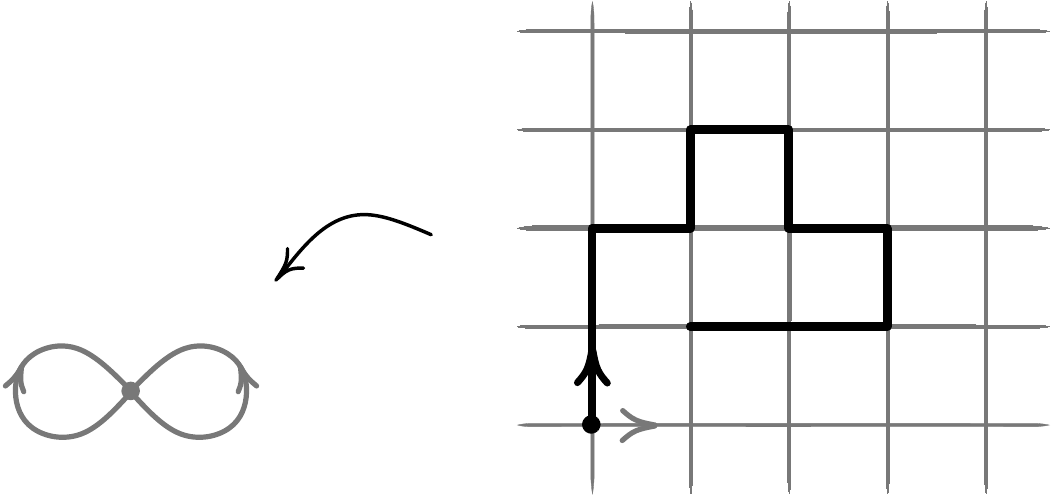}
\put(12, 15){$\Gamma$}
\put(-3, 10){$a$}
\put(25, 10){$b$}
\put(88, 38){$\tilde \Gamma$}
\put(80, 20){$\tilde R$}
\put(60, 2){$\wtilde a$}
\put(52, 10){$\wtilde b$}
\end{overpic}
\caption{
  The lift of the word $R = b^2 a b a b^{-1} a b^{-1} a^{-2}$ to the
  cover $\tilde{\Gamma}$.  }\label{fig-brown1}
\end{figure}

Let us explain Brown's criterion with a geometric
picture.  Regard the group $G = \spandef{a,b}{R = 1}$ as the quotient
of the free group $F$ on $\{a, b\}$.  Think of $F$ as the fundamental
group of a graph $\Gamma$ with one vertex and two loops.  The cover
$\tilde{\Gamma}$ of $\Gamma$ corresponding to the abelianization map $F \to
\Z^2$ can be identified with the integer grid in $\R^2$; the vertices
of $\tilde{\Gamma}$ form the integer lattice $\Z^2 \subset \R^2$ and correspond
to the abelianization of $F$.  A homomorphism $\phi \maps F \to \Z$
can be thought of as a linear functional $\R^2 \to \R$.  Now consider
our relator $R$, which we take to be a cyclically reduced word in $F$.
Let $\tilde{R}$ be the lift of the word $R$ to $\tilde{\Gamma}$, starting
at the origin (see \fullref{fig-brown1}).
An epimorphism $\phi \maps F \to \Z$ descends to $G$ if and only if
$\phi(R) = 1$.  Geometrically, this means that the kernel of $\phi$ is a
line in $\R^2$ joining the terminal point of $\tilde{R}$ to the
origin.  Turing this around, suppose $R$ is not in the commutator
subgroup of $F$ so that the endpoints of $\tilde{R}$ are distinct; in
this case there is essentially only one~$\phi$, namely projection
orthogonal to the line joining the endpoints.  (To be precise, one
should scale this projection so that $\phi$ takes values in $\Z$ rather
than $\R$, and is surjective.)

\begin{figure}\small
\centering
\begin{overpic}[scale=0.81, tics=20]{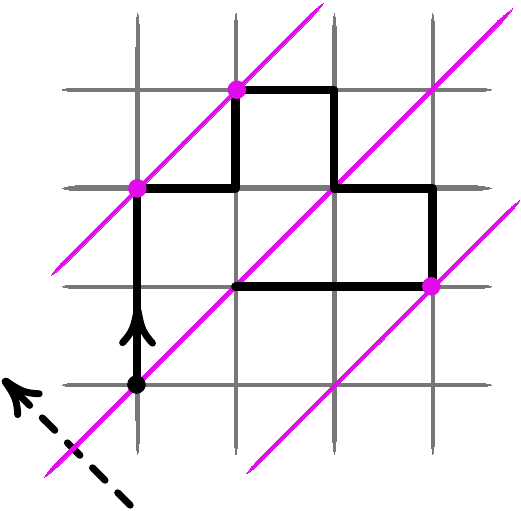}
\put(-5,18){$\phi$}
\end{overpic}
\hspace{1.7cm}
\begin{overpic}[scale=0.81, tics=20]{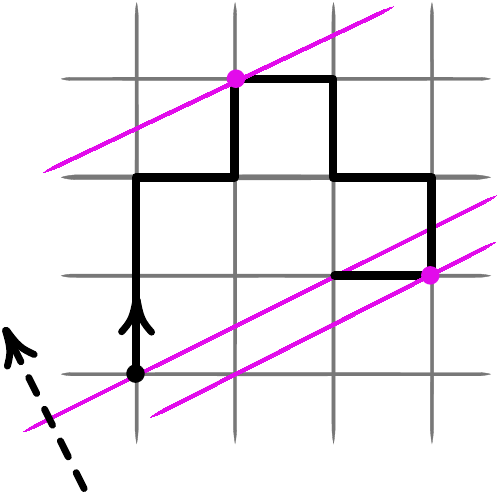}
\put(-8,25){$\phi'$}
\end{overpic}
\caption{At left is the lift of $R = b^2 a b a b^{-1} a b^{-1} a^{-2}$ to
  $\tilde{\Gamma}$.  The essentially unique homomorphism $\phi$ from $G
  = \spandef{a,b}{R = 1}$ to $\Z$ is indicated by the three diagonal
  lines which are among its level sets.  In this case, $\phi$ has two
  global maxes and one global min.  Hence the kernel of $\phi$ is
  \emph{not} finitely generated.  At right is another example with
  $R' = Ra$.  While the words differ only slightly, in this case, the
  kernel of $\phi' \maps \spandef{a,b}{R' = 1} \to \Z$ is finitely
  generated as the global extrema of $\phi$ on $\tilde{R}'$ are unique.
}\label{brown2}
\end{figure}

Now fix a $\phi$ which extends to $G$, and think of $\phi$ as a function
on the lifted path $\tilde{R}$.  Brown's criterion is in terms of the
number of global mins and maxes of $\phi$ along $\tilde{R}$.  Roughly,
$\ker(\phi) \leq G$ is finitely generated if and only if $\phi$
has the fewest extrema possible on $\tilde{R}$; that is, it has only
one global min and one global max.  \fullref{brown2} illustrates
the two possibilities.

To be precise about Brown's criterion, one needs some additional
conventions.  First, extrema are counted with multiplicities: if
$\tilde{R}$ passes through the same point of $\tilde{\Gamma}$ twice and
$\phi \maps \tilde{R} \to \R$ is maximal there, then this counts as two
maxes.  Also, the endpoints of $\tilde{\Gamma}$ can be extrema, and we
include only one of them in our count.  Finally, if the kernel of $\phi$
is horizontal or vertical, then there will be infinitely many global
extrema; in this case we count unit length segments of extrema.  To
ensure that there is no ambiguity, we state Brown's theorem a little
more combinatorially.  For our relator word $R \in F$, let $R_i$ denote
the initial subword consisting of the first $i$ letters of $R$.  The
value that $\phi$ takes on the $i^{\mathrm{th}}$ vertex of $\tilde{R}$
is then $\phi(R_i)$.
\begin{theorem}{\rm {\cite[Theorem 4.3]{Brown87}}}\label{thm-basic-brown}\qua
  Let $G = \spandef{a,b}{R = 1}$, where $R$ is a nontrivial cyclically
  reduced word in the free group on $\{a, b\}$.  Let $R_1, \ldots , R_n$
  be initial subwords of $R$, where $R_n = R$.  Consider an
  epimorphism $\phi \maps G \to \Z$.  
  
  If $\phi(a)$ and $\phi(b)$ are both nonzero, then $\ker(\phi)$ is finitely
  generated if and only if the sequence $\phi(R_1), \ldots, \phi(R_n)$ has a
  unique minimum and maximum.  If one of $\phi(a)$ and $\phi(b)$ is zero,
  then the condition is that there are
  exactly 2 mins and maxes in the sequence and that $R$ is not $a^2$
  or $b^2$.
\end{theorem}

\noindent
The statement above is equivalent to our earlier geometric one; in the
generic case, extrema of $\phi$ on $\tilde{R}$ must occur at vertices.

We'll now briefly outline the proof of Brown's theorem in a way which
elucidates its connections to the classical Alexander polynomial test
for non-fibering of a 3--manifold.  Given a two-generator, one-relator
group $G$ and an epimorphism $\phi \maps G \to \Z$, we can always change
generators in the free group to express $G$ as $\spandef{t,u}{R= 1}$
where $\phi(t) = 1$ and $\phi(u) = 0$.  The kernel of $\phi$ as a map
from the free group $\pair{t,u}$ is (freely) generated by $u_k
= t^k u t^{-k}$ for $k \in \Z$; this is because the cover of $\Gamma$
corresponding to the kernel of $\phi$ is just a line with a loop added
at each integer point.  As $\phi(R) = 0$, we have 
\begin{equation}\label{eq-relator}
R = u_{k_1}^{\epsilon_1}  u_{k_2}^{\epsilon_2} \cdots  u_{k_n}^{\epsilon_n}  \mtext{where each $\epsilon_i = \pm 1$.}
\end{equation}
The geometric condition of \fullref{thm-basic-brown} implies that
the kernel of $\phi \maps G \to \Z$ is finitely generated if and only if
the sequence $k_1, k_2, \ldots ,k_n$ has exactly one max and min.  The
``if'' part is elementary.  For instance, suppose there is a unique
minimum $k_i$, which we can take to be $0$ by replacing $R$ with $t^m
R t^{-m}$.  For any $u_j$, the relation $t^j R t^{-j}$ now implies
that $u_j$ can be expressed as a product of $u_l$'s with $l > j$.
Similarly, a unique maxima allows us to express $u_j$ as a product of
$u_l$'s with strictly smaller indices.  Thus the kernel of $\phi$ is
generated by $u_{\min(k_i)}, \ldots , u_{\max(k_i)}$.  The ``only if''
direction is more subtle, and uses the fact that the relator in a
one-relator group is in a certain sense unique.
  
We can now explain the promised connection to the Alexander
polynomial.  Let $\Delta(t)$ denote the Alexander polynomial associated to
the cyclic cover corresponding to $\phi$.  Recall the classic test in
the 3--manifold context is that if the lead coefficient of $\Delta(t)$
is not $\pm1$ (that is, $\Delta(t)$ is not \emph{monic}), then $\phi$ cannot
be represented by a fibration.  Let us see why this is true for
groups of the form we are looking at here.  First notice that $\Delta(t)$
is just what you get via the formal substitution $u_k \mapsto t^k$ in
\eqref{eq-relator}, where multiplication is turned into addition
(eg~$u_2 u_1^{-2} u_2^{-1} u_0 \mapsto 1 - 2 t$).  Thus in the
``fibered'' case where the kernel of $\phi$ is finitely generated, we
have that the lead coefficient of $\Delta(t)$ is indeed monic, as
expected.  Of course, $\Delta(t)$ can be monic and $G$ still not be
fibered.  Essentially this is because $\Delta(t)$ is only detecting
homological information; geometrically, if we look at the lift of $R$
to the cover corresponding to the kernel of $\phi \co \pair{u,t} \to \Z$,
the issue is that the Alexander polynomial only sees the homology
class of the lift of $R$, whereas Brown's criterion sees the whole
lift.  Thus you can regard Brown's test as a variant of the Alexander
polynomial test that looks at absolute geometric information instead
of homological information, and thereby gives an exact criterion for
fibering instead of only a necessary one.

\begin{remark}\label{remark-mapping-torus}
  One thing that is interesting to note about the proof sketch above
  is that when the kernel of $\phi$ is finitely generated, then in fact
  the group $G$ is the mapping torus (or HNN extension, if you prefer)
  of an automorphism of a free group.  The free group in question here
  is just $u_{\min(k_i)}, \ldots , u_{\max(k_i) - 1}$ (see
  \cite[Section 4]{Brown87} for the details).  When $G$ is the fundamental
  group of a tunnel number one 3--manifold $M$, this makes sense as
  the fiber will be a surface with boundary, whose fundamental group
  is free.
\end{remark}

\subsection{BNS invariants}\label{subsec-BNS}

Let $G = \spandef{a,b}{R = 1}$ be a two-generator one-relator group.
To apply Stallings' \fullref{thm-stallings}, we need to be able to
answer this broader version of our preceding question: does there
exist an epimorphism $\phi \maps G \to \Z$ with finitely generated
kernel?  So far, we just know how to answer this for a particular such
$\phi$.  If the relator $R$ is not in the commutator subgroup of the
free group $F = \pair{a, b}$ then there is, up to sign, a unique such
$\phi$.  So we only need to consider the case where $R \in [F, F]$;
equivalently, the relator $R$ lifts to $\tilde \Gamma$ as a closed loop.
Now there are infinitely many $\phi$ to consider, as every $\phi \maps F
\to \Z$ extends to $G$.  Fortunately, the geometric nature of
\fullref{thm-basic-brown} allows for a clean statement.  It is
natural to give the answer in term of Brown's original context, namely
the Bieri--Neumann--Strebel (BNS) invariant of a group.  This subsection
is devoted to the BNS invariant and giving Brown's full criterion.
The reader may want to skip ahead to \fullref{subsec:boxes}
at first reading; the current subsection will only be referred to in 
\fullref{subsec-rand-com} on random groups of this form.  In
particular, the main theorems about tunnel number one 3--manifolds are
independent of it.

Let $G$ be a finitely-generated group.  Broadening our point of view
to get a continuous object, consider nontrivial homomorphisms $\phi
\maps G \to \R$.  For reasons that will become apparent later, we
will consider such $\phi$ up to \emph{positive} scaling.  Let $S(G)$
denote the set of all such equivalence classes; $S(G) $
is the sphere 
\[
S(G) = \left( H^1(G, \R) \setminus 0 \right) \big/ \R^+.
\]  
The BNS invariant of $G$ is a subset $\Sigma$ of $S(G)$, which captures
information about the kernels of the $\phi$.  Rather than start with the
definition, let us give its key property (see
\cite{BieriNeumannStrebel87,Brown87} for details).
\begin{proposition}
  Let $\phi$ be an epimorphism from $G \to \Z$.  Then the kernel of
  $\phi$ is finitely generated if and only if $\phi$ and $-\phi$ are
  both in $\Sigma$.
\end{proposition}
To define $\Sigma$, first some notation.  For $\left[ \phi \right] \in S(G)$,
let $G_\phi = \setdef{g \in G}{\phi(g) \geq 0}$, which is a submonoid, but
not subgroup, of $G$.  Let $G'$ denote the commutator subgroup of~$G$,
which $G_\phi$ acts on by conjugation.  If $H$ is a submonoid of $G$, we
say that $G'$ is finitely generated over $H$ if there is a finite set
$K \subset G'$ such that $H \cdot K$ generates $G'$.  Then the BNS invariant of
$G$ is
\[
\begin{split}
\Sigma = \left\{ \left[ \phi \right] \in S(G) \ \right| &\ \mbox{$G'$
 is finitely generated over some}\\
   &\left. \mbox{finitely generated submonoid of $G_\phi$} \right\}.
\end{split}
\]
The BNS invariant has some remarkable properties---for instance, it is
always an open subset of $S(G)$.  When $G$ is the fundamental group of
a 3--manifold, $\Sigma$ is symmetric about the origin and has the following
natural description:
\begin{theorem}{\rm{\cite[Theorem E]{BieriNeumannStrebel87}}}\qua
Let $M$ be a compact, orientable, irreducible 3--manifold.  Then
$\Sigma$ is exactly the projection to $S(G)$ of the interiors of the
fibered faces of the Thurston norm ball in $H^1(M; \R)$.
\end{theorem}

In the BNS context, Brown's \fullref{thm-basic-brown} has the
following reformulation:
\begin{theorem}{\rm{\cite[Theorem 4.3]{Brown87}}}\label{thm-bns-brown}\qua
  Let $G = \spandef{a,b}{R = 1}$, where $R$ is nontrivial and
  cyclically reduced.  Let $R_i$ be initial subwords of $R$ and let
  $\left[ \phi \right] \in S(G)$.  If $\phi(a)$ and $\phi(b)$ are
  non-zero, then $\phi$ is in $\Sigma$ if and only if the sequence $\phi(R_1), \ldots,
  \phi(R_n)$ has a unique maximum.  If one of $\phi(a)$ or $\phi(b) = 0$
  vanishes, the condition is that there are exactly 2 maxes.
\end{theorem}

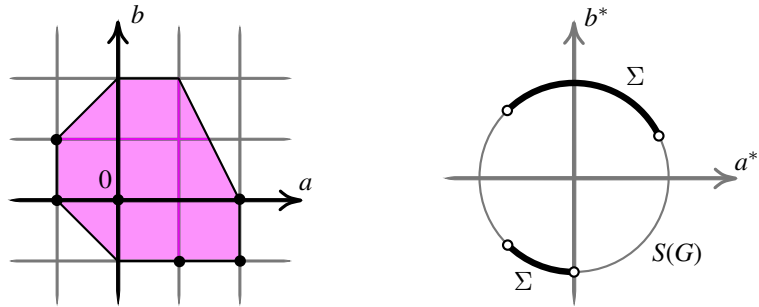
\begin{figure}\small
\centering
\begin{overpic}[scale=0.81, tics=20]{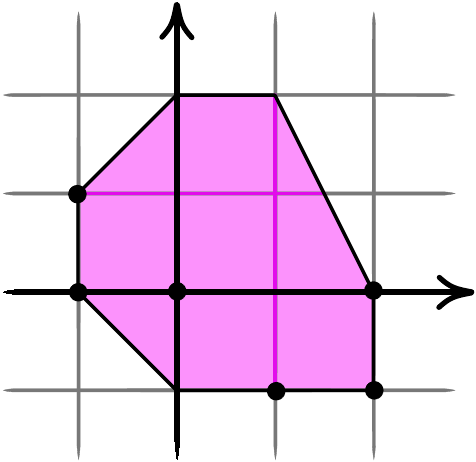}
\put(30,40){$0$}
\put(41,96){$b$}
\put(98,39){$a$}
\end{overpic}
\hspace{1.7cm}
\begin{overpic}[scale=0.81, tics=20]{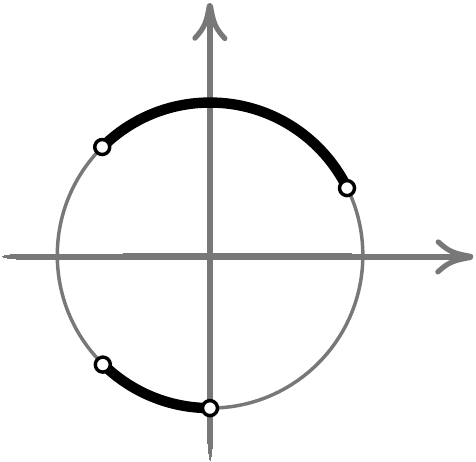}
\put(47,95){$b^*$}
\put(98,46){$a^*$}
\put(70,15){$S(G)$}
\put(23,5){$\Sigma$}
\put(61,75){$\Sigma$}
\end{overpic}
\caption{Here $R = b^{-1} a^2 b a^{-1} b^{-1} a b a^{-1} b^2 a^{-1}b^{-1}a^{-1}b^{-1} a b a^{-1}b^{-1}a $.  At left is the convex hull $C$ of $\tilde{R}$; at right is BNS invariant $\Sigma$, shown as subset of the unit circle $S(G)$ in $H^1(G; \R)$ with respect to the dual basis $\{a^*, b^*\}$.}
\label{fig-bns-example}
\end{figure}

Now consider the case when $R$ is in the commutator subgroup so that
$S(G)$ is a circle.  To describe $\Sigma$, begin by letting $\tilde{R}$ be
the lift of the relator to $\tilde{\Gamma}$ thought of as a subset of
$\R^2 = H_1(G; \R)$.  The focus will be on the convex hull $C$ of
$\tilde{R}$.  For a vertex $v$ of~$C$, let $F_v$ be the open interval
in $S(G)$ consisting of $\phi$ so that the \emph{unique} max of $\phi$ on
$C$ occurs at $v$.  Geometrically, if we pick an inner product on
$H_1(G; \R)$ so we can identify it with its dual $H^1(G;\R)$, then
$F_v$ is the interval of vectors lying between the external
perpendiculars to the sides adjoining $v$.  (Equivalently, we can
think of the dual polytope $D \subset H^1(G;\R)$ to $C$.  Then $F_v$ is
projectivization into $S(G)$ of the interior of the edge of $D$ dual
to $v$.)  We call a vertex of $C$ \emph{marked} if $\tilde R$ passes
through it more than once.  \fullref{thm-bns-brown} easily gives
\begin{theorem}{\rm{\cite[Theorem 4.4]{Brown87}}}\label{thm-bns-brown-full}\qua
  Let $G = \spandef{a,b}{R = 1}$, where $R$ is a nontrivial cyclically
  reduced word which is in the commutator subgroup.  Then the BNS
  invariant $\Sigma$ of $G$ is
  \[
  \bigcup \setdef{ F_v }{\mbox{$v$ is an unmarked vertex of $C$}}
  \]
  together with those $\phi$ whose kernels are horizontal or vertical
  if the edge of $C$ where their maxima occur has length 1 and two
  unmarked vertices.
\end{theorem}
A simple example is shown in \fullref{fig-bns-example}.
The BNS picture can also be connected to the Alexander polynomial, in
particular to the coefficients which occur at the vertices of the
Newton polygon \cite{Dunfield:norms}.

\subsection{Boxes and Brown's Criterion}\label{subsec:boxes}

In this subsection, we show how to apply Brown's criterion by breaking
up the relator~$R$ into several subwords, examining each subword
individually, and then combining the information. This works by
assigning what we call \emph{boxes} to the subwords, together with
rules for multiplying boxes.  This is crucial for adapting Brown's
criterion to efficiently incorporate the topological constraints when
$R$ is the relator for a tunnel number one 3--manifold.  That said,
the contents of this subsection apply indiscriminantly to any relator.

Let $F = \pair{a, b}$ be the free group on two generators.  Let $R \in
F$ be a cyclically
reduced word, and $1 = R_0, R_1, R_2, \ldots, R_n = R$ be the initial
subwords.  (Note that we are including $R_0 = 1$, which differs from
our conventions earlier.)  Suppose
that $\phi \maps F \to \Z$ is an epimorphism with $\phi(R) = 0$.  To apply
Brown's criterion, we are interested in the sequence $\{ \phi(R_i)
\}$, and, in particular, in the number of (global) extreme values.  We
can think of $\{ \phi(R_i)\}$ as a walk in $\Z$.  Thus, we are lead to
consider the set of finite walks $\W$ on $\Z$ which start at $0$,
where steps of any size are allowed, including pausing; that is, an
element of $\W$ is simply a finite sequence of integers whose first
term is $0$.

We now introduce \emph{boxes} to record certain basic features of a
walk $w \in \W$.  In particular, we want to remember:
\begin{itemize}
\item The final position of $w$, which we call the \emph{shift} and
  denote by $s$.
\item The maximum value of $w$, called the \emph{top} and denoted $t$.
\item The minimum value of $w$, called the \emph{bottom} and denoted $b$.  
\item The number of times the top is visited, denoted $n^t$.  We count
  in a funny way: each time the top is visited counts twice, and we
  subtract one if the first integer is the top, and subtract one if
  the last integer is the top.  Unless the walk is just $\{0\}$, this
  amounts to counting a visit at the beginning or end with a
  weight of 1 and all others with a weight of 2. We count in this way
  to make boxes well behaved under operations discussed below.
\item The number of times the bottom is visited, counted in
  the same way, denoted $n^b$.
\end{itemize}

Abstractly, a \emph{box} is simply 5 integers $(s,t,b,n^t,n^b)$ satisfying
$t \geq 0$, $b \leq 0$, $b \leq s \leq t$, and $n^t, n^b \geq 0$.  Graphically,
we denote boxes as

{\centering\leavevmode\small
\begin{xyoverpic*}{(156,89)}{scale=0.63}{\figdir/box1}
      ,(14,88)*+!D{2}
      ,(14,16)*+!U{4}
      ,(56,45)*+!D{1}
      ,(56,2)*+!U{2}
      ,(99,75)*+!D{3}
      ,(99,45)*+!U{1}
      ,(142,65)*+!D{6}
      ,(142,21)*+!U{4}
\end{xyoverpic*}

\vspace{0.1cm}
}
\noindent
and the set of all boxes is denoted $\B$.  We have a natural map
$\boxmap \maps \W \to \B$ implicit in our description above.

Given two walks $w_1, w_2 \in \W$ we can concatenate them into a walk
$w_1 * w_2$ by translating all of $w_2$ so that its initial point
matches the terminal point of $w_1$, dropping the first element of the
translated $w_2$, and joining the two lists.  For example

\vspace{0.2cm}

{\centering\leavevmode\small
\begin{xyoverpic*}{(247,49)}{scale=0.54}{\figdir/boxwalk1}
      ,(57,18)*{*}
      ,(138,18)*{=}
      ,(-28,18)*!R{\{0,-1,0,1\} * \{0,1,0,1\} = \{0,-1,0,1,2,1,2\}}
\end{xyoverpic*}

\vspace{0.2cm}

}
\noindent
where the picture at right is in terms of the graphs of the walks (see
\fullref{fig-walk} if this is unclear).  This operation makes $\W$
into a monoid, with identity element the 1--element walk $\{0\}$.  The
set $\B$ of boxes also has a monoid structure for which $\boxmap \maps
\W \to \B$ is a morphism.  Pictorially, the box multiplication is given
by

{\centering\leavevmode\small
\begin{xyoverpic*}{(369,75)}{scale=0.63}{\figdir/box2}
      ,(15,73)*+!D{2}
      ,(15,17)*+!U{4}
      ,(38,33)*{*}
      ,(64,51)*+!D{1}
      ,(64,2)*+!U{6}
      ,(96,30)*{=}
      ,(131,73)*+!D{2}
      ,(131,3)*+!U{6}
      ,(239,74)*+!D{2}
      ,(239,16)*+!U{4}
      ,(267,32)*{*}
      ,(291,74)*+!D{2}
      ,(291,2)*+!U{6}
      ,(321,29)*{=}
      ,(354,74)*+!D{4}
      ,(354,3)*+!U{6}
\end{xyoverpic*}

}
\noindent
or algebraically by the following rule.  For $i = 1,2$, let $B_i$ be the box
$(s_i, t_i, b_i, n^t_i, n^b_i)$.  Then $B_1* B_2$ is, by definition,
the box which has
\begin{itemize}
\item Shift $s = s_1 + s_2$.
\item Top $t = \max(t_1, s_1 + t_2)$.  
\item Number of top visits
\[
n^t = \begin{cases} 
  n^t_1 & \text{if $t_1 > s_1 + t_2$,} \\ 
  n^t_2 & \text{if $t_1 < s_1 + t_2$,} \\
  n^t_1 + n^t_2 & \text{ if $t_1 = s_1 + t_2$.}
\end{cases}
\]
\end{itemize}
and the corresponding rules for the bottom and $n^b$.
The identity element for this monoid is the box
\[
(0,0,0,0,0).
\]

We can also reverse a walk $w$ into another walk $\rev(w)$ of the same
length by translating by the negative of its last element and reversing
the list.  For example

\vspace{0.2cm}

{\centering\leavevmode\small
\begin{xyoverpic*}{(163,64)}{scale=0.63}{\figdir/boxwalk2}
      ,(75,46)*++!U{\rev}
      ,(-28,18)*!R{\rev\{0,1,2,1,2\} = \{0,-1,0,-1, -2\}}
\end{xyoverpic*}

\vspace{0.1cm}

}
\noindent
This is an anti-automorphism of the monoid structure on $\W$.  There
is a corresponding anti-automorphism of the monoid structure on $\B$
which is compatible with the map $\boxmap$, given pictorially by

{\centering\leavevmode\small
\begin{xyoverpic*}{(95,60)}{scale=0.63}{\figdir/box3}
      ,(14,59)*+!D{2}
      ,(14,2)*+!U{4}
      ,(46,31)*++!D{\rev}
      ,(79,59)*+!D{2}
      ,(79,2)*+!U{4}
\end{xyoverpic*}

}
\noindent
Algebraically, $\rev(s,t,b,n^t,n^b)$ is the box with: 
\begin{itemize}
\item Shift $s' = -s$.
\item Top $t' = t-s$.
\item Number of top visits $n^t{}' = n^t$.
\end{itemize}
and similarly for the bottom and number of bottom visits.

Let us now return to the setting of words in the free group $F$.
Suppose $w$ is a reduced word in $F$; here $w$ need \emph{not} be
cyclically reduced.  Let $1 = w_0, w_1, \ldots, w_n = w$ be the initial
subwords of $w$.  Let $\phi$ be any epimorphism from $F$ to $\Z$ so that
$\phi(a)$ and $\phi(b)$ are non-zero.  (We will deal with the case when
$\phi(a)$ or $\phi(b)$ is zero below.)  We set
\[
\boxmap_\phi(w) = \boxmap(\{ \phi(w_i) \})
\]
Now suppose that $v$ is a reduced word in $F$ so that the
concatenation of $w$ and $v$ is also reduced.  Then we have
\begin{align}\label{eq-box-kinda-morphism}
\boxmap_\phi(w v) &= \boxmap_\phi(w) * \boxmap_\phi (v)\\
\boxmap_\phi(w^{-1}) &= \rev(\boxmap_\phi(w)).
\end{align}
The fact that there is no cancellation when we multiply $w$ with $v$
is important here; $\boxmap_\phi \maps F \to \B$ is not a morphism.  If
we want to think of it as a morphism, we would need to take the domain
to be the monoid of strings in $\{a^{\pm1}, b^{\pm1}\}$.

An alternate way to describe $\boxmap_\phi$ is to give the values on the
generators.  If we assume that $s_1 = \phi(a) > 0$ and $s_2 = \phi(b) > 0$, then
\begin{align*}
  \boxmap(a) &= (s_1, s_1, 0, 1, 1)&   \boxmap(a^{-1}) &= (-s_1, 0, -s_1, 1, 1)\\
  \boxmap(b) &= (s_2, s_2, 0, 1, 1)&   \boxmap(b^{-1}) &= (-s_2, 0, -s_2, 1, 1)
\end{align*}
and $\boxmap$ is multiplicative on reduced words.

In the case when $\phi(a)$ is 0, we instead set
\begin{align*}
  \boxmap(a) &= (0, 0, 0, 2, 2)&   \boxmap(a^{-1}) &= (0, 0, 0, 2, 2)\\
  \boxmap(b) &= (1, 1, 0, 0, 0)&   \boxmap(b^{-1}) &= (-1, 0, -1, 0, 0)
\end{align*}
and extend by multiplicativity on reduced words.  In this case $n^t$
and $n^b$ are twice the number of segments of extrema.

Now, let's restate Brown's criterion in terms of boxes.  We say that
the top (resp. bottom) of a $\boxmap_\phi(w)$ is \emph{marked} if $n^t >
2$ (resp. $n^b > 2$).  Then \fullref{thm-basic-brown} can be
restated as:
\begin{theorem}{\rm{\cite[Theorem 4.3]{Brown87}}}\label{thm-box-brown}\qua
  Let $G = \spandef{a,b}{R = 1}$, where $R$ is a cyclically
  reduced word in the free group $\pair{a, b}$.  Consider an
  epimorphism $\phi \maps G \to \Z$.  Then $\ker(\phi)$ is finitely
  generated if and only if neither the top nor the bottom of
  $\boxmap_\phi(R)$ are marked.
\end{theorem}
Now we relate this to our original question of how to apply Brown's
criterion by breaking $R$ up into pieces.  Suppose 
\[
R = w_1 w_2 w_3 \cdots w_k
\]
where each $w_i$ is reduced and the above product involves no
cancellation to get $R$ in reduced form.  Then we have 
\begin{equation}\label{eq-rel-box-product}
\boxmap_\phi(R) = \boxmap_\phi(w_1)* \boxmap_\phi(w_2)* \cdots * \boxmap_\phi(w_k).
\end{equation}
Notice that if $B_1$ and $B_2$ are two boxes with marked tops, then
$B_1 * B_2$ also has a marked top.  Hence, if it happens that each
$\boxmap_\phi(w_i)$ has a marked top, it follows that $\boxmap_\phi(R)$
does as well without working out the product
\eqref{eq-rel-box-product}.  This yields
\begin{lemma}\label{lemma-marked-tops} 
  Let $G = \spandef{a,b}{R = 1}$, where $R$ is a cyclically
  reduced word, and consider an epimorphism $\phi \maps G \to \Z$.
  Suppose $R = w_1 w_2 w_3 \cdots w_k$ where each $w_i$ is a reduced word,
  and the product has no cancellations.  If each $\boxmap_\phi(w_i)$
  has a marked top, then the kernel of $\phi$ is infinitely generated.
\end{lemma}

In light of the above lemma, it will be useful to have criteria for
when a word has a marked top.  The one we will need is based on the
following simple observation: suppose $B$ is a box with shift $s = 0$
and $n^t \geq 2$.  Then $B * B$ has a marked top.  To apply this,
suppose $w \in F$ is a nontrivial cyclically reduced word with $\phi(w)
= 0$; taking $B = \boxmap_\phi(w)$, we claim that our observation
implies that $\boxmap_\phi(w^2) = B*B$ has a marked top.  If neither
$\phi(a)$ or $\phi(b)$ is $0$, then it is easy to see that $B$ has
$n_t \geq 2$.  If $\phi(a)$ vanishes, then there are words where $B$ has
$n^t = 0$, eg~$w = b^{-1} a b$; however, any such word is
\emph{not} cyclically reduced.  When $w$ is cyclically reduced, $n^t$
must be at least $2$.  This proves the claim that $\boxmap_\phi(w^2) =
B*B$ has a marked top.  More generally
\begin{lemma}\label{lemma-repeat-marked}
  Let $F$ and $\phi \maps F \to \Z$ be as above.  Suppose $w \in F$ is a
  nontrivial cyclically reduced word such that $\phi(w) = 0$.  If $w'$ is
  any subword of $w^n$ of length at least
  twice that of $w$, then $\boxmap_\phi(w')$ has a marked top.
\end{lemma}
\begin{proof}
  By conjugating $w$, we can assume that $w' = w^2 r$ where $r$ is an
  \emph{initial} subword of $w^n$ for some $n \geq 0$.  We have
  $\boxmap_\phi(w^2) = \boxmap_\phi(w)^2$, and since $\phi(w) = 0$ and $w$
  is nontrivial, this implies that $\boxmap_\phi(w^2)$ has a marked top.
  As $r$ is a subword of $w^n$ and $\phi(w) = 0$, the top of $\phi(w^2)$
  forms part of the top of $\phi(w')$; hence $\boxmap_\phi(w')$ has a
  marked top as well.
\end{proof}

\section{Random 1--relator groups}\label{sec-random-groups}

In this section, we consider the following natural notion of a random
2--generator 1--relator group.  Let $\cG(r)$ be the set of presentations
$\spandef{a,b}{R = 1}$ where the relator~$R$ is a cyclically reduced
word of length $r$.  While properly the elements of $\cG(r)$ are
presentations, we will usually refer to them as groups.  A random
2--generator 1--relator group of complexity $r$ is then just an element
of the finite set $\cG(r)$ chosen uniformly at random.  Now given any
property of groups, consider the probability $p_r$ that $G \in \cG(r)$
has this property; we are interested in the behavior of $p_r$ as $r \to
\infty$.  When $p_r$ has a limit $p$, it is reasonable to say that ``a
random 2--generator 1--relator group has this property with probability
$p$''; of course, $p$ is really an asymptotic quantity dependent on
our choice of filtration of these groups, namely word length of the
relator.  An example theorem is that a 2--generator 1--relator group is
word-hyperbolic with probability 1 (Gromov \cite{Gromov87},
Ol'shanskii \cite{Olshanskii92}).

In analogy with the 3--manifold situation, we say that a group $G$
\emph{fibers} if it has an epimorphism to $\Z$ with finitely generated
kernel.  As we noted in \fullref{remark-mapping-torus}, for these
types of groups fibering is equivalent to being the mapping torus of
an automorphism of a free group, which was the definition of fibered
discussed in \fullref{sub-intro-groups}.  This section is devoted
to showing that for 2--generator 1--relator groups the probability of
fibering is strictly between $0$ and $1$.  In particular:
\begin{theorem}\label{thm-random-group}
  Let $p_r$ be the probability that $G \in \cG(r)$ fibers.  Then
  for all large $r$ one has
  \[
  0.0006 < p_r < 0.975.
  \]
\end{theorem}
Experimentally, $p_r$ seems to limit to $0.94$.  It seems quite
remarkable to us that the probability a 2--generator 1--relator group
fibers is neither~$0$ nor~$1$.  In slightly different language, that
most one relator groups fiber was independently discovered
experimentally by Kapovich, Sapir, and Schupp
\cite[Section 1]{BorisovSapir03}; in that context, the proof of
\fullref{thm-random-group} shows that \cite[Theorem
1.2]{BorisovSapir03} does \emph{not} suffice to show that $G \in \cG$
group is residually finite with probability~$1$.

\fullref{thm-random-group} is strikingly different than the
corresponding result (\fullref{thm-main-nonsep}) for tunnel number
one 3--manifolds; these fiber with probability 0.  The setups of the
two theorems are strictly analogous.  Indeed, the parameter space
$\T(r)$ of tunnel number one 3--manifolds is essentially just those $G
\in \cG(r)$ which are \emph{geometric} presentations of the fundamental
group of a tunnel number one 3--manifold.  The differing results can
happen because $\T(r)$ is a vanishingly small proportion of $\cG(r)$
as $r \to \infty$; looking at Dehn--Thurston coordinates, it is clear that
$\# \T(r)$ grows polynomially in~$r$, whereas $\# \cG(r)$ grows
exponentially.  Another example of differing behavior is
word-hyperbolicity --- because of the boundary torus, the groups in
$\T$ are almost never hyperbolic, whereas those in $\cG$ almost always
are.  Still, the different behavior with respect to fibering is
surprising. As the proof of \fullref{thm-main-nonsep} will
eventually make clear, the difference stems from the highly recursive
nature of the relators of $G \in \T(r)$.

For tunnel number one manifolds with two boundary components,
\fullref{thm-main-sep} says that the probability of fibering is
still 0, despite the fact that there are now many epimorphisms to $\Z$.
In contrast, let $\cG'(r)$ be those groups in $\cG(r)$ whose defining
relation is a commutator; then it seems very likely that:
\begin{conjecture}\label{conj-comm-fiber}
  Let $p_r$ be the probability that $G \in \cG'(r)$ fibers.  Then
  $p_r \to 1$ as $r \to \infty$.
\end{conjecture}
\noindent
We will explain our motivation for this conjecture in
\fullref{subsec-rand-com}.

\subsection{A random walk problem}

Let us first reformulate the question answered by
\fullref{thm-random-group} in terms of random walks.  This will
suggest a simplified toy problem whose solution will make it
intuitively clear why \fullref{thm-random-group} is true.  We take
the point of view of \fullref{subsec:boxes}, which runs as
follows.  Start with the free group $F = \left\langle a,b\right\rangle$ and an
epimorphism $\phi \maps F \to \Z$.  A word $R \in F$ gives
1--dimensional random walk $w = \left\{ \phi(R_i) \right\}$ on $\Z$,
where the $R_i$ are the initial subwords of $R$.  Assuming neither
$\phi(a)$ or $\phi(b)$ is $0$, Brown's Criterion is then that
$\spandef{a,b}{R=1}$ fibers if and only if $w$ visits its minimum and
maximum value only once.

Unfortunately, from the point of view of the 1--dimensional walk $w$,
things are a little complicated:
\begin{enumerate}
\item The walk $w$ has two different step sizes, namely $\phi(a)$ and
  $\phi(b)$.  Moreover, the condition that $R$ is reduced means, for instance, that 
    you aren't allowed to follow a $\phi(a)$ step by a $-\phi(a)$ step.
  
\item The walk must end at $0$.  
  
\item Worst of all, the step sizes themselves are determined by the
  relator, as it is the endpoint of $R$ in the plane that determines
  $\phi$ in the first place.  Thus one can't really remove the
  2--dimensional nature of the problem.
\end{enumerate}

To get a more tractable setup, let us consider instead walks on $\Z$
where at each step we move one unit to the left or right with equal
probability.  Let $W(r)$ denote the set of such walks which both start and
end at $0$ (thus $r$ must be even).  For simplicity, let's just focus
on the maxima.   Then:
\begin{proposition}{\rm\cite{Dwass}}\qua
  A walk $w \in W(r)$ visits its maximum value more
  than once with probability $1/2$.
\end{proposition}
So the toy problem at least exhibits the neither 0 or 1 behavior of
\fullref{thm-random-group}.  While the proposition is well known,
we include a proof which, for our limited purposes, is more direct
than those in the literature.  The argument is also similar to what we
will use for \fullref{thm-random-group} itself.
\begin{figure}
  \centering\leavevmode\small
  \begin{xyoverpic}[scale=0.75]{\figdir/walk}
    ,(8,100)*++++!R{\mbox{\small position}}
    ,(3,38)*++!R{0} 
    ,(99,5)*++!U{\mbox{\small time}}
  \end{xyoverpic}
\caption{}\label{fig-walk}
\end{figure}

\begin{proof}
  We focus on the graph of a walk $w \in W(r)$, which we think of as a
  sequence of up and down segments (see \fullref{fig-walk}).  Let $U$ be
  those walks with a unique maximum.  To compute the size of $U$, we
  relate it to the set $D$ of walks which end on a down segment.
  Given $w \in U$ take the down segment immediately after the unique
  maxima, and shift it to the end to produce an element in $D$.  This
  is a bijection; the inverse $D \to U$ is to move the final down to
  immediately after the \emph{leftmost} maximum.  Thus $\# U = \#D = (1/2) \#
  W(r)$, completing the proof.
\end{proof}

If the toy problem was an exact model for
\fullref{thm-random-group}, we would expect the much lower value
of $(1/2)^2 = 1/4$ for the probability of fibering, rather than the
$0.94$ that was experimentally observed.  Next, we consider a slightly
more accurate model, where the probability of a unique maxima rises.
Consider the case where $\phi(a) = \phi(b) = 1$.  Then condition (1) above
becomes a momentum condition --- at each step there is a $2/3$ chance
of continuing in the same direction and a $1/3$ chance of turning and
going the other way.  Intuitively, this increases the chance of a
unique max since it is less likely that a repeat max is created by a
simple up-down-up-down segment.  In this case we have:
\begin{proposition}
  Consider random walks on $\Z$ with momentum as described above.  As
  the length of the walk tends to infinity, the probability of a unique
  maximum limits to $2/3$.
\end{proposition}
\begin{proof}
  We will just sketch the argument, ignoring certain corner cases
  which are why the probability $2/3$ occurs only in the limit.  The
  set of walks of length $r$ is still $W(r)$, unchanged from the
  previous proposition.  What has changed is the probability measure
  $P$ on $W(r)$ --- it is no longer uniform.  While we still have a
  bijection $f\co U \to D$ as above, it is no longer measure preserving.
  Let $U_d$ denote the set of walks with a unique max which end with a
  down, and $U_u$ those that end with a up.  Then $P(U_d) = P
  (f(U_d))$ whereas $P(U_u) = 2 P(f(U_u))$.  Also $f(U_d)$ consists
  of walks in $D$ which end in \emph{two} downs; thus $f(U_d)$
  contributes $2/3$ of the measure of~$D$, whereas $f(U_u)$
  contributes only $1/3$.  Combining gives $P(U) = 4/3 P(D) = 2/3$
  as desired.
\end{proof}

Unfortunately, our approach seems to fail when we allow differing step
sizes as in~(1), even ignoring the momentum issue.  The problem is that 
while the maps between $U$ and $D$ are still defined, they are
no longer bijective.  We turn now to the proof of
\fullref{thm-random-group} which uses similar but cruder methods
which have no hope of being sharp.  

\begin{proof}
  As above, let $\cG(r)$ be our set of 1--relator groups, which we will
  always think of as the set of cyclically reduced words $R$ in $F =
  \left\langle a , b \right\rangle$ of length $r$.  As a first step, we compute
  $\# \cG(r)$.  Counting reduced words, as opposed to cyclically
  reduced words, is easy: there are 4 choices for the first letter and
  3 choices for each successive one, for a total of $4 \cdot 3^{r-1}$.
  What we need to find $\# \cG(r)$ is the probability that a reduced word
  is cyclically reduced.  Thinking of a reduced word $w$ as chosen at
  random, the relationship between the final letter and the
  initial one is governed by a Markov chain whose distribution
  converges rapidly to the uniform one.  Thus the distribution of the
  final letter is (nearly) independent of the first letter, and so the
  odds that $w$ is cyclically reduced is $3/4$.  Thus $\# \cG(r)$ is
  asymptotic to $3^r$.  A more detailed analysis, not needed for what
  we do here, shows that $\# \cG(r) = 3^r + 1$ when $r$ is odd, and
  $3^r + 3$ when $r$ is even.
  
  Let $\cG_0(r)$ denote those $R$ which are not in the commutator
  subgroup, and so that the unique epimorphism $\phi \maps
  \spandef{a,b}{R = 1} \to \Z$ does not vanish on either $a$ or $b$.
  It is not hard to see that the density of $\cG_0(r)$ in $\cG(r)$
  goes to $1$ as $r \to \infty$.  Thus in the remainder of the proof, we
  work to estimate the probability $p_r'$ that $\spandef{a,b}{R = 1}$
  fibers for $R \in \cG_0(r)$.
  \begin{figure}\small
    \centering\leavevmode
    \begin{xyoverpic*}{(184,106)}{scale=0.72}{\figdir/walk2}
      ,(18,93)*+!DL{\phi}
      ,(121,93)*+!DL{\phi}
      ,(80,0)*{\mbox{(a)}}
    \end{xyoverpic*}
    \hspace{1.6cm}
    \begin{xyoverpic*}{(184,106)}{scale=0.72}{\figdir/walk3}
      ,(22,100)*+!DL{\phi}
      ,(127,100)*+!DL{\phi}
      ,(80,0)*{\mbox{(b)}}
    \end{xyoverpic*}
    \caption{}\label{fig-words}
  \end{figure}
  To bound it from above, we construct an injection from $\cG_0(r-4)$
  into the \emph{non-fibered} subset of $\cG_0(r)$.  In particular,
  given $R \in \cG_0(r-4)$, go to the first global maximum and insert a
  commutator as shown in \fullref{fig-words}(a).  As we inserted a
  commutator, $\phi$ is unchanged, but we now have enough maxima to see
  that it is non-fibered. To see that this map is injective, observe
  that there is an inverse process: go to the first global maximum and
  delete the next 4 letters.  Thus
  \[
  1 - p_r' \geq \frac{\# \cG_0(r - 4)}{\# \cG_0(r)} \approx  3^{-4} \quad \mbox{and hence $p_r' < 0.988$ for large $r$}.
  \]
  To improve this, note that we can also insert a commutator at the
  global minimum; the images of these two injections of $\cG_0(r-4)$
  into the non-fibered words have some overlap coming from
  $\cG_0(r-8)$.  Thus, pretending for convenience that $\# \cG_0(n) = 3^n$,
  we have
  \[
  1 - p_r' \geq 2 \cdot 3^{-4} - 3^{-8}  \quad \mbox{and hence $p_r' < 0.975$ for large $r$},
  \]
  as desired.  
  
  To estimate the number of fibered words in $\cG_0(r)$, we inject
  $\cG_0(r-8)$ into them by inserting a commutator at both the first
  global min and the first global max as shown in
  \fullref{fig-words}(b).  In fact, there are four distinct ways of
  doing this, depending on which way we orient the two commutators.
  Thus
  \[
  p_r' \geq 4 \frac{\# \cG_0(r - 8)}{\# \cG_0(r)} \approx  4 \cdot 3^{-8}  > 0.0006 \quad \mbox {for large $r$,}
  \]
 completing the proof.  
\end{proof}

\subsection{The case of $\cG'$} \label{subsec-rand-com}

We end this section by giving the motivation for
\fullref{conj-comm-fiber}.  Consider $G \in \cG'(r)$.  Every $\phi
\maps F \to \Z$ extends to $G$.  If we fix $\phi$, the proof of
\fullref{thm-random-group} shows that the probability that
$\ker(\phi)$ is finitely generated is at least $0.0006$.  The intuition
is that if we fix two such epimorphisms $\phi$ and $\phi'$, then as $r \to
\infty$ the event that $\ker(\phi)$ is finitely generated becomes
independent from the corresponding event for $\phi'$.  This should
extend to any finite collection of $\phi_i$, and as each event has
probability at least $0.006$, independence means that at least one of
the $\phi_i$ will have finitely generated kernel with very high
probability.  Increasing the number of $\phi_i$ would allow one to show
that the probability is at least $1-\epsilon$ for any $\epsilon$, and hence the
probability limits to~$1$.  The different $\phi_i$ should have
independent behavior for the following reason.  As described in
\fullref{subsec-BNS}, whether $G$ fibers depends on the convex
hull of the relator~$R$.  In particular, at the vertices of the convex
hull we care about whether $R$ passes through them multiple times.
Any given $\phi$ picks up its global extrema from some pair of these
vertices.  As long as the extreme vertices associated to $\phi$ and
$\phi'$ are distinct, then whether they are repeated vertices should be
independent.  Thus the key issue is simply that the number of vertices
on the convex hull of $R$ should grow as $r \to \infty$.  In the related
questions that have been studied, the number of vertices grows like
$\log(r)$, and we don't expect the situation here to be any different.
See Steele \cite{Steele2002} and the references therein for details.

\section{Efficient implementation of Brown's algorithm}\label{sec-eff-brown}

In this section, we discuss how to efficiently decide whether a given
curve in our parameter space gives a fibered tunnel number one
manifold.  Here, ``efficiently'' means in time which is polynomial in the
\emph{log} of the Dehn--Thurston coordinates.  The method we present is
also crucial to the proof of \fullref{thm-main-nonsep}.

\subsection{Train tracks}

\begin{figure}
  \centering
  \includegraphics[scale=0.9]{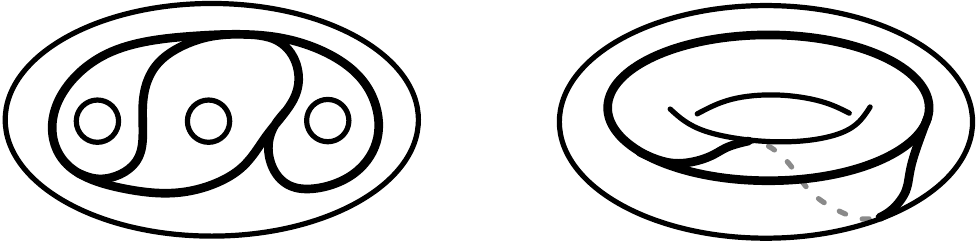}
  \caption{Train tracks on a 4--punctured $S^2$ and on a torus.}
  \label{fig:ex-tracks}
\end{figure}

\begin{minipage}{4in}Our main combinatorial tools for studying curves 
on surfaces are train
tracks and their generalizations.  Roughly, a \emph{train track} in a
surface $\Sigma$ is a 1--dimensional CW complex $\tau$ embedded in $\Sigma$,
which is made up of 1--dimensional branches joined by trivalent
switches.  Here, each switch has one incoming branch and two outgoing
branches.  See \fullref{fig:ex-tracks} for examples, and
\cite{PennerHarerBook} for details.
Associated to a train track $\tau$ is a \emph{space of weights} (or
transverse measures). This consists of assignments of weights $w_e \in
\R_{\geq 0}$\break \end{minipage}
\lower 35pt\hbox to 1.5in{\hglue 15pt \begin{xyoverpic}[scale=0.9]{\figdir/n8}
    ,(0,82)*++!DL{3}
    ,(100,100)*!L{2}
    ,(100,68)*!L{1}
  \end{xyoverpic}}\vglue -20pt

to each branch $e$ of $\tau$, which satisfy the switch
condition: at each switch the sum of the weights on the two outgoing
edges is equal to the weight on the incoming one.  The space of
weights is denoted $\ML(\tau,\R)$, and $\ML(\tau, \Z)$ denotes those where
each $w_e \in \Z$.  As shown at right, an integral measure $w \in
\ML(\tau, \Z)$ naturally specifies a multicurve which lies in a small
neighborhood of $\tau$, that is, is \emph{carried by $\tau$.}  More
generally, $\ML(\tau,\R)$ parameterizes measured laminations carried by
$\tau$.  For suitable train tracks, called \emph{complete} train tracks,
$\ML(\tau,\R)$ gives a chart on the space of measured lamination on the
underlying surface $\Sigma$.\eject

\subsection{Interval exchanges}

\begin{figure}
\centering
\raisebox{0.6cm}{\includegraphics[scale=1.2]{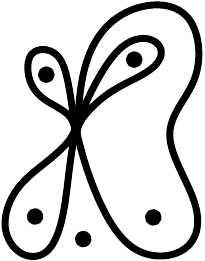}}\hspace{1.7cm} \includegraphics[scale=0.65]{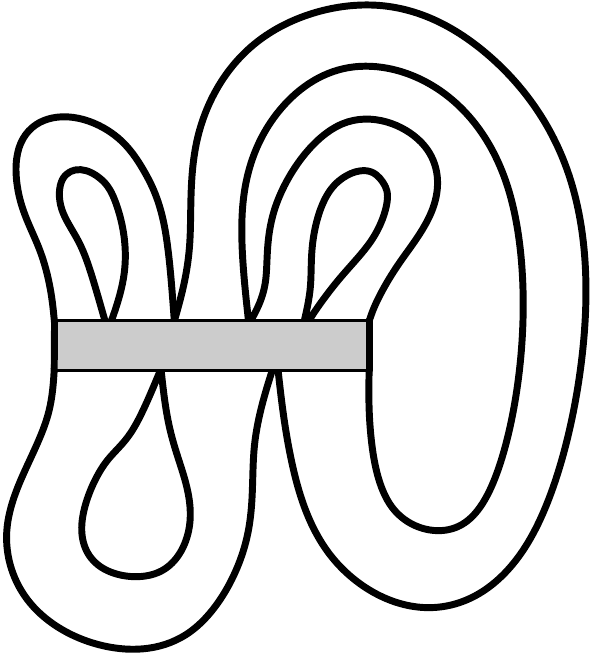}
\caption{At left is an interval exchange $\tau$ on a 5--punctured $S^2$.  At right is a regular neighborhood of $\tau$.}\label{fig:int-exchange}
\end{figure}

A more generalized notion of train tracks is to allow switches where
there are an arbitrary number of incoming and outgoing branches.
Here, we will focus on the class where there is just one switch.
These are called \emph{interval exchanges} for reasons we will see
shortly.  An example on a 5--punctured $S^2$ is shown in
\fullref{fig:int-exchange}.  As you can see from that figure, a
regular neighborhood of such an interval exchange can be decomposed
into a thickened interval (shaded) whose top and bottom are
partitioned into subintervals which are exchanged by means of bands
(the thickened branches).  A $w \in \ML(\tau, \R)$ can be
thought of as assigning widths to the bands so that the total length
of the top and bottom intervals agree.


There are two kinds of bands.  Those that go from the top to the
bottom are termed \emph{orientation preserving} since that is how they
act on their subintervals.  Those joining a side to itself are called
\emph{orientation reversing}.  An interval exchange gives rise to a
natural dynamical system which we describe in the next subsection.  In
that context, they have been studied extensively since the 1970s.
However, generally only orientation preserving bands are allowed in
that literature; we will refer to such exchanges as \emph{classical
  interval exchanges}.

For the rest of this section, one could easily work with train tracks
instead of interval exchanges.  However, the use of interval exchanges
has important technical advantages in the proof of
\fullref{thm-main-nonsep}.

\subsection{Rauzy induction and determining connectivity}\label{subsec-splitting}

\begin{figure}
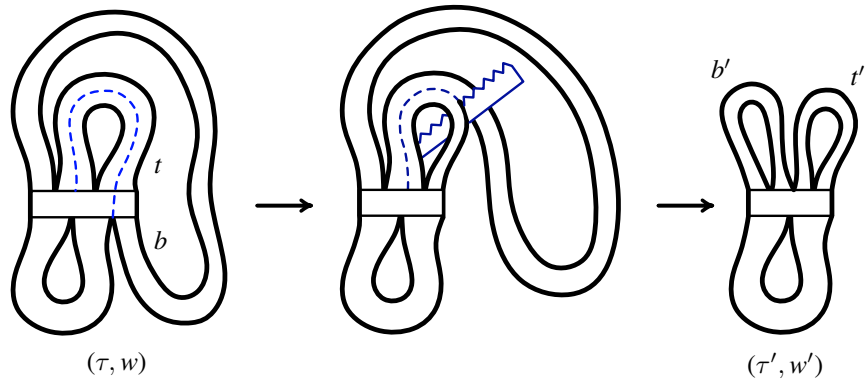
\small
\centering\leavevmode
\begin{xyoverpic*}{(355,140)}{scale=0.9}{\figdir/n11}
  ,(54,62)*++!LD{t}
  ,(54,51)*++!LU{b}
  ,(45,0)*++!U{(\tau,w)}
  ,(305,105)*+!RD{b'}
  ,(348,102)*+!LD{t'}
  ,(324,0)*++!U{(\tau',w')}
\end{xyoverpic*}
\caption{Splitting an interval exchange.}\label{fig:splitting}
\end{figure}

Suppose $\tau$ is an interval exchange in a surface $\Sigma$, and $w \in
\ML(\tau, \Z)$.  We will now describe how to determine the number of
components of the associated multicurve.  This method is also the
basis for our efficient form of Brown's algorithm.  The basic
operation is called Rauzy induction in the context of interval
exchanges, and splitting or sliding in the context of train tracks.
Starting with $\tau$ and $w \in \ML(\tau, \R)$, we will construct a new
pair $(\tau', w')$ realizing the same measured lamination in $\Sigma$.  To
begin, consider the rightmost bands, $t$ and $b$, on the top and
bottom respectively.  First suppose that $w_t > w_b$.  Then we slice
as shown in \fullref{fig:splitting} to construct $(\tau', w')$.  
The new band~$t'$ has weight $w'_{t'} = w_t - w_b$, and the other
modified band $b'$ has weight $w'_{b'} = w_b$.  The other weights are
of course unchanged.  If instead $w_b < w_t$ one does the analogous
operation, flipping the picture about the horizontal axis.

If $w_t = w_b$, then one simply cuts through the middle
interval and amalgamates the bands $b$ and $t$ together.  If $t$ and
$b$ happen to be the same band, this splits off an annular loop.
Hence we enlarge our notion of an interval exchange by allowing the
addition of a finite number of such loops.  We will denote Rauzy
induction, which we usually call \emph{splitting}, by $(\tau, w) \searrow
(\tau', w')$.

Now suppose we start with an integral measure $w \in \ML(\tau, \Z)$ and
want to determine the number of components of the associated
multicurve.  We can split repeatedly to get a sequence $(\tau_i,
w_i)$ carrying the same multicurve.  Here, there is no reason to
remember the (increasingly complicated) embeddings of the $\tau_i$ into
$\Sigma$.  That is, you should think of the $\tau_i$ as \emph{abstract}
interval exchanges, not embedded in any particular surface.  At each
stage, either some weight of $w_i$ is strictly reduced, or we split
off a loop and reduce the number of bands in play.  In the end, we are
reduced to a finite collection of loops labeled by elements of $
\Z_{\geq 0}$; the sum of these labels is the number of components in
the multicurve.

The procedure just described is not always more efficient than the
most naive algorithm; in particular the number of steps can be equal
to $\abs{w} = \max \{ w_b \}$.  For instance, take $\tau$ to be the
exchange on the torus shown in \fullref{fig:ex-tracks} and set the
weights on the bands to be $1$ and $n$.  However, Agol, Hass, and
W.~Thurston have shown that if one adds a ``Dehn twist'' operation,
then the number of steps becomes polynomial in $\log \abs{w}$~\cite{AgolHassThurston}.  

\subsection{Computing other information}

As pointed out in \cite{AgolHassThurston}, one can adapt this
framework to compute additional invariants of the multicurve $\gamma$
given by $w \in \ML(\tau, \Z)$.  The cases of interest for us are derived
from the following setup.  Suppose we know that $\gamma$ is connected,
and fix generators for $\pi_1(\Sigma)$.  Let's see how to find a word in
$\pi_1(\Sigma)$ representing $\gamma$ in terms of splitting.  We can take
the basepoint for $\pi_1(\Sigma)$ to lie in the base interval for $\tau$,
which we think of as very small.  An oriented band of $\tau$ thus gives
rise to an element of $\pi_1(\Sigma)$.  We think of each band as being
\emph{labeled} with this element.  Now suppose we do a splitting
$(\tau, w) \searrow (\tau', w')$.  Resuming the notation of
\fullref{fig:splitting}, we presume $w_t > w_b$.  Orient the bands
$t$ and $b$ so that both orientations point vertically at the
right-hand side of the base interval.  The new bands $b'$ and
$t'$ of $\tau'$ now inherit orientations as well.  If we use $L$ to
denote our $\pi_1(\Sigma)$ labels, then these transform via
\begin{equation}\label{eq:label-mult-rule}
L(b') = L(b) \cdot L(t) \quad \mbox{and} \quad L(t') = L(t),
\end{equation}
with all the other labels remaining unchanged.  Since we are presuming
that $\gamma$ is connected, if we continue splitting in this manner we
eventually arrive at a single loop with weight $1$.  The label on that
loop is then a word representing $\gamma$ in $\pi_1(\Sigma)$.

Since in the end we recover a full word representing $\gamma$, this
splitting algorithm takes time at least proportional to the size of
that word; this can certainly be as large as $\abs{w}$.  The real
payoff is when we want to compute something derived from this word
which carries much less information.  For instance, suppose we want to
know the class of $\gamma$ in $H_1(\Sigma)$.  Then we can use labels which
are the images of the $\pi_1(\Sigma)$ labels under the quotient $\pi_1(\Sigma)
\to H_1(\Sigma)$.  In this way, we can compute the class of $\gamma$ in
$H_1(\Sigma)$ in time polynomial in $\log \abs{w}$
\cite{AgolHassThurston}.

\subsection{The algorithm: boxes on interval exchanges}\label{subsec-boxes-on-IE}

We now turn to the main question at hand.  Suppose $H$ is our genus 2
handlebody, and $\gamma \in \T$ a non-separating simple closed curve.  We
want to (efficiently) decide if the associated tunnel number one
manifold $M_\gamma$ fibers.  Suppose that $\gamma$ is given to us in terms of
weights $w$ on a train track $\tau_0$.  Using the technique of the last
subsection, we can quickly compute the element $\gamma$ represents in
$H_1(H)$.  Let us further suppose that this is not $0$; we now have
determined the essentially unique epimorphism $\phi \maps \pi_1(M_\gamma) \to
\Z$.  In light of \fullref{cor:tunnel-one-stallings}, to decide
if $M_\gamma$ fibers we just need to apply Brown's Criterion to decide if
the kernel of $\phi$ is finitely generated.  In
\fullref{subsec:boxes}, we described how to implement Brown's
Criterion by breaking the defining relation $R$ up into subwords and
using boxes to capture the needed information about these subwords.
Roughly, we initially label $\tau$ by corresponding words of $F =
\pi_1(H) = \pair{a,b}$, and then immediately replace each word $v$ with
$\boxmap_\phi(v)$.  Then at each split, we will combine the boxes via
box multiplication following the rule \eqref{eq:label-mult-rule}.  At
the end we will be left with a single loop labeled by $\boxmap_\phi(R)$
to which we can apply \fullref{thm-box-brown}.  However, in order
for the final box to really be $\boxmap_\phi(R)$, we must restrict
the initial train track $\tau$:
$\boxmap_\phi$ is not a morphism unless we take the domain to be the
monoid of words in $\{a^{\pm 1},b^{\pm 1}\}$, rather than the free group
$F$ itself.  In particular, we must ensure that each time we split the
interval exchange there is no cancellation in the $F$ labels.

While one way of thinking about the final label on $\gamma$ is via the
splitting process, it can be also thought of less dynamically.  Focus
on a neighborhood of $\tau$, and think of each band as having a
vertical dividing line in the middle of its length.  Fix a transverse
orientation for the divider.  Suppose we label each band by a word in
the generators and their inverses.  The label for each band should be
thought of as affixed to
its dividing line.  A connected curve $\gamma$ carried by $\tau$ has a
sequence of intersections with the dividers; reading off the labels as
we go around $\gamma$ (inverting the label if the direction of travel does not
match the transverse orientation of the divider) and taking the
product gives the final word.  The
final word is well-defined up to the choice of starting
point and choice of orientation of $\gamma$.  We say that $\gamma$ is
\emph{tight} if the final word is cyclically reduced.

\begin{definition}
  Let $\tau$ be an interval exchange with bands labeled by elements of
  $F$.  We say that $\tau$ is \emph{tightly labeled} if every $\gamma$ carried by $\tau$ is tight.  
\end{definition}
\begin{figure}[ht!]
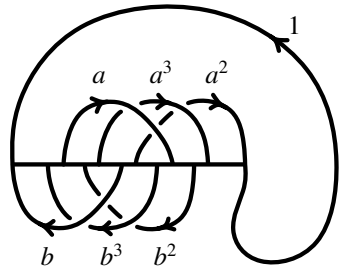
\small
\centering\leavevmode
\begin{xyoverpic*}{(154,120)}{scale=0.81}{\figdir/n12}
,(126,102)*+!DL{1}
,(42,76)*++!D{a}
,(71,76)*++!D{a^3}
,(97,76)*++!D{a^2}
,(18,15.5)*++!U{b}
,(48,17)*++!U{b^3}
,(73,17)*++!U{b^2}
\end{xyoverpic*}
\caption{A tightly labeled interval exchange.}\label{fig:tight-exchange}
\end{figure}

An example of a tightly labeled interval exchange is given in
\fullref{fig:tight-exchange}; the point is that as we
run along $\gamma$, if we cross a label which is a power of $a$,
then the next label other than $1$ that we encounter is a power of $b$.
Concluding the above discussion, we have:
\begin{lemma}\label{lemma:tight-computes-correctly}
  Let $\tau$ be an interval exchange tightly labeled by $F$.  Suppose
  $\phi \maps F \to \Z$ is an epimorphism, and let $\gamma$ be a
  connected simple closed curve carried by $\gamma$.  Split $\tau$ until we
  get a single loop labeled by $R$.  Now start back at the beginning
  and replace the labels on $\tau$ by $L \mapsto \boxmap_\phi(L)$, and again
  split until we get a single loop labeled with a box $B$.  Then $B =
  \boxmap_\phi(R)$.
\end{lemma} 
 
As at the beginning of this subsection, suppose we are given $\gamma \in
\T$, our parameter space of tunnel number one 3--manifolds.  We assume
that $\gamma$ is given to us in terms of Dehn--Thurston coordinates as in
\fullref{subsec:measured-lam-prob}.  As we next describe, the
constraints (2--4) in \fullref{subsec:measured-lam-prob} on the
Dehn--Thurston coordinates of $\gamma \in \T$ allow us to put $\gamma$ on a tight
interval exchange closely related to the one given in
\fullref{fig:tight-exchange}.  In terms of
\mbox{\fullref{fig-DT-coor}}, consider the punctured torus $T$ bounded by
$\delta$ containing $\alpha$.  The intersection of $\gamma$ with $T$ consists of
at most 3 parallel families of arcs.  Thinking homologically, it is
easy to see that we can orient things so that the labels on these
families are $a^i$, $a^j$ and $a^{i + j}$ where $i, j \in \Z$.
Condition (4) of \fullref{subsec:measured-lam-prob} means that
none of $\{i, j, i + j\}$ are zero.  (If there are fewer than 3
families of arcs, we add in empty families to increase the number to
3, making the discussion uniform.  This can be done fairly
arbitrarily, and we can thus ensure that $\{i, j, i + j\}$ are all
nonzero.)  Each family of arcs will contribute one band to our final
$\tau$.  The same picture is true for the other punctured torus.  We can
now make an interval exchange $\tau$ by taking these bands in the
punctured tori and adding one additional band to allow us to effect
the twist around $\delta$.
\begin{figure}[ht!]
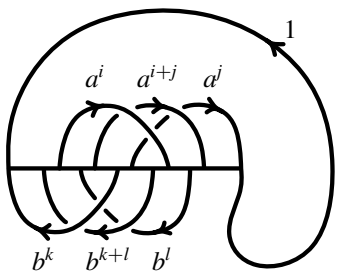
\small
\centering\leavevmode
\begin{xyoverpic*}{(154,120)}{scale=0.81}{\figdir/n12}
,(126,102)*+!DL{1}
,(42,76)*++!D{a^i}
,(71,76)*++!D{a^{i + j}}
,(97,76)*++!D{a^j}
,(18,17)*++!U{b^k}
,(48,17)*++!U{b^{k+l}}
,(73,17)*++!U{b^{l}}
\end{xyoverpic*}
\caption{
  A standard starting interval exchange.  Here the elements $a^i$,
  $a^j$, $a^{i +j}$, $b^k$, $b^l$, and $b^{k + l}$ of $F$ are not the
  identity.  }\label{fig:standard-track}
\end{figure}
The result is shown in \fullref{fig:standard-track}.  (If you are
worried here about whether the final band is consistent with the
orientation of the twisting about $\delta$, note that since a full Dehn
twist about $\delta$ extends over $H$, we can make this twisting have any
sign we like without changing~$M_\gamma$.)  We will call interval
exchanges of this form \emph{standard starting interval exchanges}.
The same reasoning used above shows that $\tau$ is tightly labeled.  For
future reference we record:
\begin{lemma}\label{lemma-std-track}  
  Every $\gamma \in \T$ is carried by one of countably many standard
  starting interval exchanges, each of which is tightly labeled by $F$.
\end{lemma}

To summarize, here is the procedure to efficiently decide if $M_\gamma$
fibers over $S^1$, provided that $\gamma$ is nonzero in $H_1(H)$.  (If
$\gamma$ is zero in $H_1(H)$, there is not a unique $\phi$ to test.  While
Brown's algorithm adapts to work very elegantly to this situation (see
\fullref{subsec-BNS}), it is unclear if it can be implemented
efficiently in this case.  You may have to remember too much in the
appropriate labels.)  First, it is straightforward in terms of the
Dehn--Thurston coordinates to put $\gamma$ on a standard initial exchange
$\tau_0$.  Then run the splitting once using $H_1(H)$ labels to
determine $\phi$.  Once $\phi$ is known, we go back the beginning and
relabel $\tau_0$ with $\boxmap_\phi$ labels.  Now run the splitting to the
bottom again.  In light of \fullref{lemma:tight-computes-correctly}
and \fullref{thm-box-brown}, the label on the final loop
determines if the kernel of $\phi$ is finitely generated, and hence if
$M_\gamma$ fibers over the circle.  Since there isn't much to a box,
really just 5 numbers bounded by the square of the initial weights,
the running time will still be polynomial in
the log of initial weights, or equivalently in the size of the
Dehn--Thurston coordinates.

\section{The idea of the proof of the main theorems}\label{sec-main-idea}

In this section, we explain the outline of the common proof of the
main results of this paper, Theorems~\ref{thm-main-nonsep} and
\ref{thm-main-sep}.  For concreteness, let us focus attention on
\fullref{thm-main-nonsep}.  The basic idea is to analyze the
algorithm given in the last section, and prove that it will report
``non-fibered'' with probability tending to 1 as the input curve $\gamma$
becomes more and more complicated.  Recall the setup is that we are
given a connected curve $\gamma$ on $\partial H$ carried by a tightly labeled
copy of $\tau_0$.  Let $G = \pi_1(M_\gamma) = \spandef{a,b}{R=1}$, and let us
assume we are in the generic case where there is a unique epimorphism
$\phi \maps G \to \Z$.  In the algorithm of
\fullref{subsec-boxes-on-IE}, we start with $\boxmap_\phi$ labels on
$\tau_0$, and then split $(\tau_0, \gamma)$ repeatedly, at each stage replacing 
one of the box labels by its product with another box.  After
many splittings, we are left with a single loop.  What we want to show
is that the top of the remaining box is very likely to be marked, and
thus $\ker(\phi)$ is infinitely generated by
\fullref{thm-box-brown}.

Suppose at some point in the splitting sequence we get to a state
where \emph{every} box label has a marked top.  As
\fullref{lemma-marked-tops} shows, this property will persist as we
continue to split, and so this state alone implies that $\ker(\phi)$
is infinitely generated.  Thus we can stop at this stage, even though
we may still have a million splits to go before we get down to a
loop.  This is in fact what happens when you implement the algorithm
of \fullref{subsec-boxes-on-IE} --- typically, you end up with
marked boxes on all the bands long before you have completed splitting.
The strategy for the proof is therefore to show that having all boxes
marked becomes increasing likely as you do more splits.

We begin by discussing a much simpler problem, whose solution follows
the same strategy, and then explain how the approach must be modified to
account for the constraints in the actual topological situation.

\subsection{A toy problem}

Our toy problem is the following.  Suppose we start with a finite set
of boxes; for concreteness, let us say there are 7 of them to match
the number of bands in one of the standard starting interval
exchanges shown in \fullref{fig:standard-track}.  At each ``split''
we pick two distinct boxes $A$ and $B$ at random, and replace $A$ with
a random selection from
\[
A*B, \quad A*B^{-1}, \quad B*A, \quad \mbox{or} \quad B^{-1} * A,
\]
where here $B^{-1} = \rev(B)$ is the reverse or ``inverse'' box
described in \fullref{subsec:boxes}.  The simplification here is
that any pair of boxes can interact at any stage, whereas for interval
exchanges the pair is fixed by the topology.  However, for interval
exchanges one expects that over time all pairs of labels will
interact, and so it is reasonable to hope that the behavior of the toy
problem will tell us something about the real case of interest.

For any choice of the initial 7 boxes where at least one is
non-trivial, we'll show:
\begin{theorem}\label{thm-toy-splits}
  The probability that all boxes have marked tops after $n$ splits
  goes to $1$ as $n \to \infty$.
\end{theorem}  
\noindent
This result is an easy consequence of the following lemma. 
\begin{lemma}\label{lemma-toy-splits}
  Let $B_1, B_2, \ldots, B_7$ be boxes where at least one is nontrivial.
  Then there is a sequence of 14 splits so that the resulting boxes
  all have marked tops.
\end{lemma}
Assuming the lemma, here's the proof of the theorem.  Start with our
initial boxes, split $14$ times, then $14$ more times, etc.  While we
don't know anything about the state of the boxes at the start of each
chunk of $14$ splits, the lemma tells us that there is at least a
$28^{-14}$ chance that the next $14$ splittings result in all boxes
marked.  Since the splittings in distinct chunks are chosen
independently, the probability that some box is \emph{not} marked
after $14 k$ splits is at most $(1 - 28^{-14})^k$.  Hence the
probability that all of the boxes are marked converges to $1$ as the
number of splits goes to infinity.  To see that we are on the right
track, notice that this exponential decay of the probability of
fibering is consistent with the experimental data in
\fullref{sec-experiment}.

We now head toward the proof of \fullref{lemma-toy-splits}; before
continuing, the reader may want to review the notation of
\fullref{subsec:boxes}.   First, if $X$ is a nontrivial box
with shift 0, and then $X^2$ has a marked top.  Moreover, we can create
a box with shift 0 from any pair of boxes $A, B$ by taking the
commutator $X = A*B *A^{-1}*B^{-1}$.  Now, given our 7 boxes, it is
not easy to create a commutator by the splitting moves; however, one
thing we can do is that if $A, B, C$ are three of the boxes, then we
can replace $C$ with
\begin{equation}\label{eq-Cprime}
C' = C*(A*B *A^{-1}*B^{-1} )^2
\end{equation}
by doing 8 splits.  Roughly, if the box $C$ is shorter than $X =
A*B*A^{-1}*B^{-1}$ then the top of $C'$ should come from the $X^2$
term, and hence $C'$ will have a marked top as well.  We turn now to
the details.

\begin{proof}[Proof of \fullref{lemma-toy-splits}]
  Let $B_1, B_2, \ldots, B_7$ be our initial boxes.  Given the way
  splittings work, we can replace $B_i$ with $B_i^{-1}$ without really
  changing anything, so let's normalize things so that the shifts
  satisfy $s_i \geq 0$.  We will denote the top of $B_i$ by~$t_i$.  Let
  $A$ be the nontrivial $B_i$ with largest top, which we denote $t_a$.
  Pick two of the remaining $B_i$, and denote them $B$ and $C$.  The
  splitting sequence we will use is: first do the sequence of $8$
  splits replacing $C$ by $C'$ in \eqref{eq-Cprime}; then for each
  $B_i$ which is not $C'$, do a split to replace $B_i$ with $B_i *
  C'$.
  
  To see that all boxes are marked at the end of this, first consider
  $X = A*B*A^{-1}*B^{-1}$.  By choice of $A$, the top of $X$ satisfies
  $t_x \geq t_a \geq t_i$ for all $i$.  Then $C'$ above has a
  \emph{marked} top coming from $X^2$ at height $t_x + s_c$, where
  $s_c$ is the shift of $C$.  Thus for any $B_i$ we have that $B_i *
  C'$ has a marked top at height $s_i + s_c + t_x$ since $t_x \geq t_i$.
  Thus we can make all the boxes marked in only 14 splits.
\end{proof}

\subsection{Outline of the proof of the main theorems}

With \fullref{thm-toy-splits} in hand, we now explain how the
approach generalizes to \fullref{thm-main-nonsep}.  As mentioned
above, the difficulty we need to incorporate is that with interval
exchanges, we have much less freedom in how the boxes are changed at
each step.  Despite this, the analog of \fullref{lemma-toy-splits}
is still true.  In \fullref{sec-magic-sequence}, we prove that
there is a single ``magic'' splitting sequence which always gives
marked boxes.  Interestingly, unlike \fullref{lemma-toy-splits},
some assumptions must be made on the starting boxes for this to be
true; however, these always hold for those boxes arising in the
algorithm of \fullref{subsec-boxes-on-IE} (see
\fullref{rmk-why-hard} for more).  In the toy problem, going from
\fullref{lemma-toy-splits} to \fullref{thm-toy-splits} was
essentially immediate.  The key features were that each block of 14
splits is independent of the others, and the desired splitting
sequence always has a definite probability of occurring.  For interval
exchanges, these things are more subtle; essentially, what we need is
that splitting complete genus 2 interval exchanges is ``normal''.
This is shown in \fullref{sec-ubiquity}, relying on work of
Kerckhoff \cite{Kerckhoff1985}.  Finally, in
\fullref{sec-main-thm-pf} we assemble the pieces just discussed
with work of Mirzakhani \cite{Mirzakhani2004} to complete the proof of
\fullref{thm-main-nonsep}.

\section{The magic splitting sequence}\label{sec-magic-sequence}

This section is devoted to a lemma which is one of the central
ingredients in the proof of the main theorems.  Suppose $\gamma$ is a
connected curve on the boundary of our genus~2 handlebody $H$.
Suppose $\gamma$ is carried by some interval exchange $\tau_0$.  Roughly, we
show that if the splitting sequence of $(\tau_0, \gamma)$ has a certain
topological form, then the tunnel number one 3--manifold $M_\gamma$
does not fiber over the circle.  Before stating the lemma, we discuss
its precise context.

\subsection{Complete interval exchanges}

We will work with interval exchanges $\tau$ in a surface $\Sigma$ for which
$\ML(\tau, \R)$ gives a chart for $\ML(\Sigma)$; in particular, we work with
complete interval exchanges, which we now define.  Let $\Sigma$ be a
closed surface of genus at least $2$.  An interval exchange~$\tau$ in
$\Sigma$ is called \emph{recurrent} if there is a $w \in \ML(\tau, \R)$ where
every band has positive weight.  For an interval exchange, the switch
condition is that the sum of the weights of the orientation
reversing bands on the top is equal to the corresponding quantity for the
bottom.  Thus, $\tau$ is recurrent if and only if there are orientation
reversing bands on both sides, or no such bands at all.  The exchange
$\tau$ is \emph{complete} if it is recurrent, and every complementary
region is an ideal triangle.  When $\tau$ is complete, the natural map
$\ML(\tau, \R) \to \ML(\Sigma)$ is a homeomorphism onto its image; if we
restrict the domain to $w$ which are nowhere zero, then we get a
homeomorphism onto an open subset of $\ML(\Sigma)$.  (When working
with train tracks, one also requires transverse recurrence in the
definition of completeness.  However, an easy application of
\cite[Corollary 1.3.5]{PennerHarerBook} shows that any interval exchange is
transversely recurrent.)

In what follows, the embedding of $\tau$ into $\Sigma$ is not really
relevant.  Thus, we will tend to think of interval exchanges
abstractly, that is, as not embedded in any particular surface.  If
one presumes that the complementary regions of $\tau$ are ideal
polygons, then one can reconstruct $\Sigma$ from the combinatorics of $\tau$
alone.  Thus, it makes sense to speak of an abstract interval exchange
as being a complete interval exchange on a genus 2 surface.

\subsection{Statement of the lemma}

Now let's give the setup for the main result of this section.  We are
interested in splitting sequences of interval exchanges.  Suppose $\tau$
is an interval exchange.  Given $w \in \ML(\tau, \R)$, as described in
\fullref{subsec-splitting} we can split $(\tau, w)$ to $(\tau', w')$,
which is denoted by $(\tau, w) \searrow (\tau', w')$.  Independent of the choice
for $w$, there are (at most) 3 distinct possibilities for $\tau'$; in
the notation of \fullref{fig:splitting}, the 3 possibilities
correspond to $w_t > w_b$, $w_t < w_b$, and $w_t = w_b$.  We will use
the notation $\tau \searrow \tau'$ to indicate that $(\tau,w)$ splits to
$(\tau',w')$ for some $w \in \ML(\tau, \R)$.

From now on we will look at complete interval exchanges on a genus 2
surface.  Of special importance is the exchange shown in
\fullref{fig:standard-track}, which we will call $\tau_0$.  Suppose
we have a splitting sequence
\[
S \maps \tau_0 \searrow \sigma_1 \searrow \sigma_2 \searrow \cdots \searrow \sigma_n
\]
where the $\sigma_i$ are also complete genus 2 interval exchanges.  For a
multicurve $\gamma \in \ML(\tau_0,\Z)$, we say that $\gamma$ \emph{exhibits} $S$
if the initial part of the splitting sequence of $(\tau_0, \gamma)$ is
\[
\tau_0 \searrow \tau_1 \searrow \tau_3 \searrow \cdots \searrow \tau_m
\]
where the tail $ \tau_{m-n} \searrow \tau_{m-n} \searrow \cdots \searrow \tau_m$ is abstractly
isomorphic to $S$.   The point of this section is to prove:
\begin{lemma}\label{lemma-magic}
  There exists a splitting sequence of complete genus 2 interval
  exchanges
  \[
  S \maps  \tau_0 \searrow \tau_1 \searrow \tau_2 \searrow \cdots \searrow \tau_n
  \]
  such that the following holds.  Suppose $\gamma$ is a connected simple
  closed curve on $\partial H$ carried by a tightly labeled copy of $\tau_0$.
  If the splitting sequence for $(\tau_0, \gamma)$ exhibits $S$, then the
  manifold $M_\gamma$ does not fiber over the circle.
\end{lemma}

Note that in the lemma $\gamma$ is allowed to be either separating or
non-separating.  The splitting sequence $S$ will be referred to as the
\emph{magic splitting sequence}.  It is quite complicated, and so even
to describe it, we must first give another point of view on splitting
interval exchanges.

\subsection{Flexible splitting of interval exchanges}

Let $\tau$ be an interval exchange with some initial measured lamination
$w \in \ML(\tau, \R)$.  Let $I = [0, L]$ be the base interval of the
exchange.  In the notation of \fullref{subsec-splitting}, during
each splitting we reduce the length of $I$ by $\min(w_b, w_t)$.  We
now describe a way of seeing the result of several splittings at once.
Consider a subinterval $J = [0, L']$ for some $L' < L$.  Now take a
knife to $(\tau, w)$ and begin to slice it starting at some notch
between bands, following the lamination as you go.  Here $J$ should
be viewed as indestructible, and when the knife collides with $J$ you
stop.  Repeat for each of the other notches until no more progress can
be made.  Thus we have created a new pair $(\tau', w')$ which describes
the same lamination.  Moreover $(\tau', w')$ is a stage of the splitting
sequence for $(\tau, w)$, in particular the one right before the base
interval shrinks to a proper subinterval of $J$.  In describing this
cutting process, there is no need to cut each notch down to $J$ in one
go --- we can start somewhere, cut for a bit, and work somewhere else
before coming back to finish the job.

\subsection{The magic sequence}

\begin{figure}[ht!]
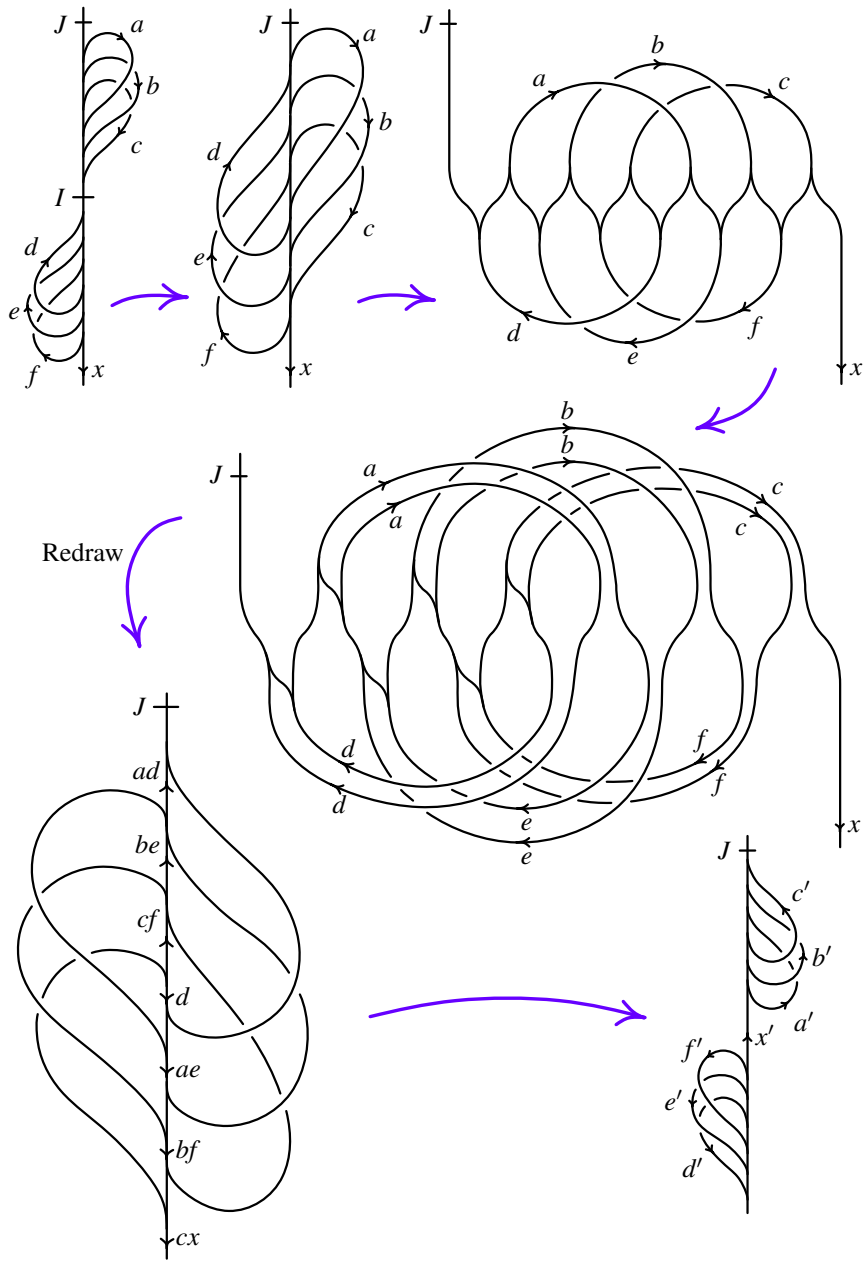
\small
  \centering\leavevmode
  \begin{xyoverpic*}{(368,557)}{scale=0.85}{\figdir/universal}
    ,(26,549)*+!R{J}
    ,(47,542)*+!LD{a}
    ,(54,522)*+!L{b}
    ,(47,501)*+!LU{c}
    ,(25,472)*+!R{I}
    ,(13,443)*+!RD{d}
    ,(5,421)*+!R{e}
    ,(14,401)*+!RU{f}
    ,(30,395)*+!L{x}
    ,(117,549)*+!R{J}
    ,(150,538)*+!LD{a}
    ,(157,506)*+!L{b}
    ,(150,464)*+!LU{c}
    ,(94,485)*+!RD{d}
    ,(87,445)*+!R{e}
    ,(92, 411)*+!RU{f}
    ,(122,396)*+!L{x}
    ,(188,549)*+!R{J}
    ,(238,518)*+!RD{a}
    ,(284,533)*+!D{b}
    ,(334,517)*+!LD{c}
    ,(227,420)*+!RU{d}
    ,(273,407)*+!U{e}
    ,(322, 423)*+!UL{f}
    ,(366,396)*+!L{x}
    ,(51,307)*+!RD{\mbox{Redraw}}
    ,(95,349)*+!R{J}
    ,(163,345)*+!DR{a}
    ,(168,334)*+!U{a}
    ,(244,370)*+!D{b}
    ,(244, 355)*+!D{b}
    ,(331, 337)*+!DL{c}
    ,(327, 331)*+!UR{c}
    ,(147, 220 )*+!D{d}
    ,(143, 210 )*+!U{d}
    ,(226, 200 )*+!U{e}
    ,(226, 185)*+!U{e}
    ,(304, 222)*+!D{f}
    ,(312, 220)*+!U{f}
    ,(365, 193)*+!L{x}
    ,(61, 246 )*+!R{J}
    ,(67, 210)*+!RD{ad}
    ,(67, 177)*+!RD{be}
    ,(67, 142)*+!RD{cf}
    ,(67, 117)*+!L{d}
    ,(67, 84)*+!L{ae}
    ,(67, 48)*+!L{bf}
    ,(67, 9)*+!L{cx}
    ,(320,182)*+!R{J}
    ,(340,156)*+!LD{c'}
    ,(349,135)*+!L{b'}
    ,(341,114)*+!LU{a'}
    ,(324,100)*+!L{x'}
    ,(308,84)*+!DR{f'}
    ,(299,71)*+!R{e'}
    ,(308,49)*+!UR{d'}
  \end{xyoverpic*}
  \caption{
    The first part $S_1$ of the magic splitting sequence, which starts
    in the top right with $\tau_0$, and ends in the bottom left with
    another copy of $\tau_0$.}
  \label{fig:universal-splitting}
\end{figure}

We now describe the first part of the magic sequence, using the setup
just given.  In \fullref{fig:universal-splitting}, we start in the
upper left with $\tau_0$ drawn as a train track; in this picture of
$\tau_0$ the top and the bottom of the vertical segment are identified
and this convention persists throughout the figure.  Also marked on
$\tau_0$ is the base interval $I$ and the smaller initial subinterval
$J$.  \fullref{fig:universal-splitting} describes a splitting
sequence $S_1$ where we split $\tau_0$ down to $J$.  The final train
track is again the interval exchange $\tau_0$.  A choice of $w \in
\ML(\tau_0, \R)$ which induces this splitting is not indicated; it can
be determined a posteriori by choosing nowhere zero weights $w'$ on
the final copy of $\tau_0$, and working backwards up the sequence to
determine $w$.  In proving \fullref{lemma-magic}, it will be very
important to know what the labels (in the sense of
\fullref{subsec-boxes-on-IE}) are on the final copy of $\tau_0$.
The initial $\tau_0$ is labeled by $\{ x, a, b, c, d, e, f \}$ in
$\pi_1(H)$.  When following through what happens in
\fullref{fig:universal-splitting}, you should view each label as
sitting precisely at the indicated arrow.  The final labels are as
follows, where $A = a^{-1}$, etc.~and the vertical bars should simply
be ignored for now.
\begin{align}
\label{eqn-new-labels}
x' &= XCFBEADcfbead \notag \\ 
a' &= DAEBFCdaebf | bead  \notag \\ 
b' &= DAEBFCdae|ad \notag \\
c' &= DAEBFCd| \\ 
d' &= XCFBEAD|cx   \notag \\ 
e' &= XCFBEADcf|bfcx \notag \\ 
f' &= XCFBEADcfbe|aebfcx \notag  
\end{align}

The second part of the magic splitting sequence is much easier to
describe.  Suppose $w \in \ML(\tau_0, \R)$ has larger weight on the $x$
band than the $c$ band.  Then splitting one step gives us $\tau_0$
again.  We refer to this as the \emph{stable splitting} of $\tau_0$.
Let $S_2$ consist of repeating the stable splitting 6 times.  The key
properties of $S_2$ that we will use are as follows.  The sequence
$S_2$ affects the labels by replacing $a$ with $X a x$ and the same
for $b$ and $c$.  Suppose the initial splitting sequence of some $w \in
\ML(\tau_0, \R)$ is $S_2$, what does that tell us about $w$?  Let $w_x$
be the width of the $x$--band, and $w_r = 2 (w_a + w_b + w_c)$ be the
difference between the length of the base interval of $\tau_0$ and the
$x$--band.  It is not hard to see that the splitting sequence for
$(\tau_0,w)$ starts off with $S_2$ if and only if $w_x \geq w_r$.  More
generally, the splitting sequence starts off with $n$ copies of $S_2$
if and only if $w_x \geq n w_r$.
 
Finally, the magic sequence itself is
\[
S = S_1 \searrow S_2 \searrow S_2 \searrow S_2 \searrow S_2 \searrow S_2 \searrow S_2.  
\]

\subsection{Proof of the lemma}

We first outline the approach to proving
\fullref{lemma-magic}.  Recall the setup is that we are given a
connected curve $\gamma$ on $\partial H$ carried by a tightly labeled copy of
$\tau_0$.  Supposing the splitting sequence for $(\tau_0, \gamma)$ exhibits
$S$, then we need to show that $M_\gamma$ does not fiber over the circle.
If $G = \pi_1(M_\gamma) = \spandef{a,b}{R=1}$, for each epimorphism $\phi
\maps G \to \Z$ we need to show that $\ker(\phi)$ is infinitely
generated.  From now on, we view $\phi$ as fixed.  Roughly, the strategy
of the proof is we start with $\boxmap_\phi$ labels on $\tau_0$ as in
\fullref{subsec-boxes-on-IE}, and then split until the sequence
$S$ occurs.  At that point, all the box labels will have their tops
\emph{marked} in the sense of \fullref{subsec:boxes}.  Thus the
relation $R$ for $G$ coming from $\gamma$ is a product where each factor
has $\boxmap_\phi$ with a marked top.  By \fullref{lemma-marked-tops},
the subgroup $\ker(\phi)$ is then infinitely generated.  For technical
reasons, the proof of the lemma deviates slightly from the above sketch,
though we suspect it could be made to work on the nose by (further)
complicating $S$.   We turn now to the details.

\begin{proof}[Proof of \fullref{lemma-magic}]
  Continuing with the notation above, we can assume that the
  \emph{initial} splitting sequence of $(\tau_0, \gamma)$ is $S$.  Beyond
  the tightness restraint, the only thing we need to show about the
  initial labels on $\tau_0$ is that $\phi(x) = 0$, as follows.  Consider
  the curve $\epsilon$ carried by $\tau_0$ where the
  only nonzero weight is $w_x = 1$.  From the way $\epsilon$ divides up
  $\tau_0$, we see that $\epsilon$ is a \emph{separating} curve in $\partial H$.  Thus
  the word $x \in \pi_1(H)$ associated to $\epsilon$ lies in the commutator
  subgroup, and so $\phi(x) = 0$.

  Now split $(\tau_0, \gamma)$ along $S_1$, getting back to $\tau_0$
  with the labels as in \eqref{eqn-new-labels}.  Next do the
  $S_2$ splitting twice, so that the we have new labels
  \[
  a'' = (X')^2a'(x')^2 \quad b'' = (X')^2b'(x')^2  \quad c'' = (X')^2c'(x')^2 
  \]  
  with the others unchanged.  Changing tacks, rather than implement
  the 4 remaining $S_2$ splits, we just use the fact noted above that
  this means that the $x'$ band is much wider than all the other bands
  put together.  In particular, we now split starting from the
  \emph{left} side of the base interval rather than the right, and do
  two analogs of the $S_2$ splittings there. This has the effect of
  changing the labels by
  \[
  c'' = (x')^2c'(X')^2 \quad   e'' = (x')^2e'(X')^2 \quad   f'' = (x')^2f'(X')^2
  \]  
  We claim that the boxes of $a'', b'', c'', d'', e'', f''$ all
  have marked tops.  Consider for instance
  \[
  a'' = (DAEBFCdaebfcx)^2 \cdot DAEBFCdaebf | bead \cdot(XCFBEADcfbead)^2.
  \]
  Now we have $\phi(x') = \phi(X') = 0$, by the same argument that shows
  $\phi(x) = 0$.  If we look at the part of $a''$ lying to the left of
  the vertical line, we see a subword of $(X')^3$ that is long enough
  so \fullref{lemma-repeat-marked} implies that its box has a marked
  top.  Similarly, the right half of $a''$ also has a box with a
  marked top, and thus $\boxmap_\phi (a'')$ has a marked top.  The same
  argument works for all the other labels except $x'$, where the
  division of the word into two parts is indicated in
  \eqref{eqn-new-labels}.
  
  To conclude the proof, it would be enough to know that
  $\boxmap_\phi(x')$ has a marked top.  Rather than show this, first
  note that $\boxmap_\phi((x')^n)$ has a marked top for any $n \geq 2$.
  We still have two $S_2$ splits ``left'' that we haven't used.  This
  means the weight on the $x'$ band is large enough that the places
  where the $x'$ band is attached on the top and bottom have
  considerable overlap.  As a result, anytime the curve $\gamma$ enters
  the $x'$ band it goes around it at least twice before moving to a
  different band.  Thus the defining relation $R$ for $G$ is a
  (non-canceling) product of the words
  \[
  a'', b'', c'', d'', e'', f'', (x')^n \quad \mbox{for $n \geq 2$}
  \]
  and their inverses.  Since all of these have marked boxes, the subgroup
  $\ker(\phi)$ is not finitely generated by
  \fullref{lemma-marked-tops}, completing the proof.
\end{proof}

\begin{remark}\label{rmk-why-hard}
  The subtle issue in the proof is not that marked boxes tend to
  persist, but why any marked boxes must be created in the first
  place.  For the latter, the fact that $\phi(x)=0$ is crucial, as
  otherwise the box of a power like $(x')^2$ need not have a marked
  top.  We found situations where you start with a train track and put
  random box labels on it so that no sequence of splitting operations
  can produce marked boxes on all the bands.
\end{remark}

\section{Ubiquity of splitting sequences}\label{sec-ubiquity}

In this section, we study splitting sequences of complete genus 2
interval exchanges.  We show that certain finite splitting sequences
(including the magic one)
are ubiquitous --- in a suitable sense, they occur in the splitting
sequence of almost every multicurve.  Before giving the precise
statement, it is worth considering the analogous fact for continued
fractions.  Let $s_1, s_2, \ldots , s_n$ be a sequence of positive
integers.  Given a rational number $p/q \in [0,1]$, we can look at the
partial quotients in its continued fraction expansion.  We can ask if
the sequence $\{s_i\}$ occurs as successive terms in these partial
quotients. In fact, the probability that this occurs goes to $1$ as
$\max(p,q) \to \infty$.  The point of this section is to prove the
corresponding result, \fullref{thm-ubiquity}, for genus~2 interval
exchanges.

Among complete interval exchanges on a genus 2 surface, we focus on
the exchange~$\tau_0$ shown in \fullref{fig:standard-track}.  Suppose
we have a splitting sequence
\[
S \maps \tau_0 \searrow \sigma_1 \searrow \sigma_2 \searrow \cdots \searrow \sigma_n
\]
where the $\sigma_i$ are also complete genus 2 interval exchanges.  As in
the last section, for a multicurve $\gamma \in \ML(\tau_0,\Z)$, we say that
$\gamma$ \emph{exhibits} $S$ if its splitting sequence contains a copy of
$S$.  We will show that, asymptotically, the set of $\gamma$ in $\ML(\tau_0,
\Z)$ which exhibit $S$ has density $1$.  More precisely:
\begin{theorem}\label{thm-ubiquity}
  Let $\tau_0$ be the complete genus 2 interval exchange specified
  above, and $S$ a splitting sequence of $\tau_0$ consisting of complete
  interval exchanges.  Set
  \[
  \cS = \setdef{ \gamma \in \ML(\tau_0, \Z)}{ \mbox{$\gamma$ exhibits $S$}}.
  \]
  If $U$ is a bounded open set in $\ML(\tau_0, \R)$, then
  \[
  \frac{  \num{ \cS \cap t U} } { \num{ \ML(\tau_0, \Z) \cap t U}} %
    \to 1 \mtext{as} t \to \infty.
  \]
\end{theorem}

This theorem may strike the reader as excessively narrow: what is so
special about genus 2 and this particular choice of $\tau_0$?  The
choice of $\tau_0$ can be broadened, but not to the point where every
complete genus 2 exchange can play its role.  The restriction to genus
2 is necessary for the proof as given; it is unclear to us if it
is actually needed.  Both of these issues are discussed at length
below.  The proof of \fullref{thm-ubiquity} is based on applying a
criterion of Steve Kerckhoff \cite{Kerckhoff1985}, which we discuss in
the next two subsections.

\subsection{Splitting interval exchanges as a dynamical system}

The bulk of the proof of \fullref{thm-ubiquity} will be
carried out in a slightly different setting.  This subsection is
devoted to describing this setting, and stating the analog therein of
\fullref{thm-ubiquity}.  Let $\cI$ be the set of complete genus 2
interval exchanges up to (abstract) isomorphism.  It is easy to see
that any $\tau \in \cI$ has 7 bands, so $\cI$ is finite.  Consider the set
\[
X = \coprod_{\tau \in \cI}\ML(\tau, \R).
\]
with the measure coming from Lebesgue measure on each $\ML(\tau, \R)$.
We can construct a transformation of $X$ by sending $(\tau, w)$ to its
splitting $(\tau', w')$.  This is well-defined provided $\tau'$ is also
complete; equivalently, we are not in the corner case where the
weights on the rightmost bands are equal, ie,~$w_t = w_b$ in the
notation of \fullref{fig:splitting}.  As this only concerns a
measure $0$ subset of $X$, we simply delete all $(\tau, w)$ in $X$ which
pass through the corner case somewhere in their splitting sequence.
(Although we have just removed all the integral points which are the
focus of \fullref{thm-ubiquity}, this will not inconvenience us
greatly.)  Thus given a $(\tau, w)$ in $X$, we get an infinite sequence
$\tau = \tau_1 \searrow \tau_2 \searrow \cdots $ of elements of $\cI$.  We are interested in
the behavior of this sequence for generic choices of the initial $(\tau,
w)$, where here generic means except for a measure~$0$ subset of $X$.
The goal is to show that $X$ with this transformation is
\emph{normal}; that is, a finite splitting sequence
\[
S \maps  \sigma_1 \searrow \sigma_2 \searrow \cdots \searrow \sigma_n
\]
which can happen must happen infinitely often for almost every $(\tau,
w)$.  Here, the phrase ``can happen'' must be interpreted correctly,
as we next discuss.

Make $\cI$ into a directed graph by putting an edge from $\tau$ to $\tau'$
if $\tau \searrow \tau'$.  The splitting sequence of some $(\tau, w)$ now
corresponds to an infinite directed path in this graph.  A complication
is that $\cI$ is not directedly connected; that is, there are
$\tau$ and $\tau'$ so that no directed path starting at $\tau$ ends at
$\tau'$.  In terms of our dynamical system $X$, this gives rise to
transient states which generically appear only finitely many times in
a splitting sequence.  If $\tau$ can be joined to $\tau'$ by a directed
path starting at $\tau$, we write $\tau \searrow \! \! \! \searrow \tau'$.  Let
\[
\cI_S = \setdef{ \tau \in \cI }{ \mbox{$\sigma \searrow\! \! \! \searrow \tau$ for every $\sigma \in \cI$}},
\]
which we refer to as the \emph{sink} of $\cI$.  A priori, $\cI$ could
be empty; for instance, this is the case if $\cI$ is not connected.

We now restrict attention to the subset of $X$ sitting over the sink
\[
X_S = \coprod_{\tau \in \cI_S} \ML(\tau, \R)
\]
which is closed under the splitting transformation.  Our precise
definition of normality is:

\begin{definition} 
  Suppose the sink $\cI_S$ is non-empty.  Then the transformation on
  $X_S$ is \emph{normal} if every finite splitting sequence
  \[
  \sigma_1 \searrow \sigma_2 \searrow \cdots \searrow \sigma_n \quad \mbox{which is contained in 
  $\cI_S$}
  \]  occurs infinitely often in the splitting sequence of almost
  all $(\tau, w) \in X_S$.  
\end{definition}

We will show the following, which will easily imply
\fullref{thm-ubiquity}: 

\begin{theorem}\label{thm-dyn-ubiq}
  As above, let $X_S$ be the set of weights on complete genus 2 interval
  exchanges lying in the sink $\cI_S$.  Then the splitting
  transformation on $X_S$ is normal.  Moreover, $\tau_0$ lies in the
  sink $\cI_S$.
\end{theorem}

Our proof of \fullref{thm-dyn-ubiq} is a direct application of a
normality criterion of Steve Kerckhoff \cite{Kerckhoff1985}, which he
used to prove the analogous result for \emph{classical} interval
exchanges, namely those without orientation reversing bands.  Readers
familiar with \cite{Kerckhoff1985} may wonder why we are working with
non-classical interval exchanges rather than just using train tracks,
as \cite{Kerckhoff1985} states the analog of
\fullref{thm-ubiquity} for complete train tracks in any genus.
There are two reasons for this.  The first is that using interval
exchanges makes it much easier to understand $\cI$ explicitly.  The
other reason is that the proof given in \cite{Kerckhoff1985} for the
train track case is incomplete.  The proof there involves two steps,
the first is to establish that a certain combinatorial criterion
implies normality, and the second to check that this criterion holds
for train tracks.  The criterion, which is the one we use here, is
certainly strong enough to ensure normality; the problem occurs in the
second step, as the criterion is violated for certain explicit train
tracks.  Kerckhoff informs us that he noticed this problem as well, and
that there should be a weaker combinatorial condition which still
ensures normality but also holds in the train track setting.
Kerckhoff is planning on publishing a correction along these lines.

We turn now to the combinatorics of $\cI$ and its sink.  It turns out
that there are 201 genus 2 interval exchanges in $\cI$, and 190 in the
sink $\cI_S$.  Also, as mentioned above, $\tau_0 \in \cI_S$.  We checked
this by brute force computer enumeration, and will not further justify
these facts here.  This is not a difficult calculation, requiring only
100 lines of code and a few minutes of computer time.  The part that's
actually used in \fullref{thm-ubiquity}, namely that $\tau_0$ is in
the sink of its connected component of $\cI$, is particularly
straightforward: First, start with $\tau_0$ and keeps splitting until
you don't generate any new exchanges.  Then for each of the 190
exchanges so generated, check that you can split them all back to
$\tau_0$.

Enumerating all of $\cI$, in particular to see that it is connected,
requires a little more thought to set up.  One approach is to think of
a complete genus 2 interval exchange $\tau$ as constructed by starting
with a disc and adding 7 bands.  The boundary of the disc is divided
into 28 segments, alternating between a place where we attach a band,
and a gap between bands.  Further, two of the gaps are distinguished
because they form the vertical sides of the thickened interval that is
the base of the exchange; we call these special gaps \emph{smoothed}.
If we ignore the smoothings, the complementary regions of $\tau$ are
ideal polygons of valence either $\{3,3,3,5\}$ or $\{3,3,4,4\}$.
Conversely, provided the complementary regions are of this form,
choosing smoothings of gaps in the larger complementary regions gives
a genus 2 interval exchange.  As there are only $135{,}135$ possible
gluings for the bands, one can simply try them all and thus calculate
$\cI$.

For train tracks, the structure of $\cI$ and its sink can also be quite
complicated.  For a 4--punctured sphere, for instance, the vast
majority of the complete train tracks do not lie in the sink.  It
would be quite interesting to answer the following question.

\begin{question}
  Find a type of train track like object where the elements of the
  sink can be characterized topologically, and for which the
  corresponding dynamical system is normal.
\end{question}

\subsection{Normality after  Kerckhoff}

This subsection is devoted to the proof of \fullref{thm-dyn-ubiq}.
We begin by describing Kerckhoff's normality criterion.  First, we
will need to work with \emph{labeled} interval exchanges, that is, interval
exchanges where the bands are labeled by integers from $1$ to the
number of bands.  If $\tau$ is a labeled exchange, and we split $\tau$ to
$\tau'$, the convention for the labels on $\tau'$ is as follows.  The
bands of $\tau'$ that come unchanged from $\tau$ retain their labels; in
the notation of \fullref{fig:splitting}, the two modified bands
$t'$ and $b'$ get the labels of $t$ and $b$ respectively.

The set of labeled complete genus 2 exchanges is of course finite, and
we focus on the subset $\cI'$ where the underlying unlabeled exchange
lies in the sink $\cI_S$.  Again, $\cI'$ has structure of a directed
graph with edges given by splittings.  The forgetful map $\cI' \to
\cI_S$ is a covering map.  Fix a connected component $\cI'_0$ of
$\cI'$; the covering map $\cI'_0 \to \cI_S$ is also surjective.  (It
appears that $\cI'$ consists of two components.)  We claim that
$\cI'_0$ is its own sink.  Call a directed graph \emph{strongly
  connected} if for all vertices $v_1$ and $v_2$, there is a directed
path from $v_1$ to $v_2$.  What we need is equivalent to:
\begin{lemma} 
  Let $G$ be a finite strongly connected directed graph, and $H \to G$
  a finite covering map.  If $H$ is connected then it is strongly
  connected.
\end{lemma}
\begin{proof}
  A \emph{cycle} in a directed graph is a union of edges which form a
  directed closed loop, where no vertex is visited more than once.
  Observe that a connected directed graph is strongly connected if and
  only if every edge is part of a cycle --- the point here is that
  going almost all the way around a cycle effectively allows us go backwards
  along a directed edge.  The result now follows by noting that the
  preimage of a cycle in $G$ is a disjoint union of cycles in $H$.
\end{proof}

From now on, we work to show that 
\[
Z = \coprod_{\tau \in \cI_0'} \ML(\tau , \R)
\]
with its splitting transformation is normal; this suffices to prove
\fullref{thm-dyn-ubiq}.  We now set up the terminology needed to
state Kerckhoff's normality criterion.  Let $\tau$ be a labeled interval
exchange.  The rightmost places where the bands are glued on the top
and bottom are called the \emph{critical positions}.  The bands in
those positions are called the \emph{critical bands}, and are denoted
$t$ and $b$ respectively.  During a splitting move, the band whose
width is reduced is said to be split by the other.  For example, in
\fullref{fig:splitting} where $w_t > w_b$, we say that $t$ is
split by $b$, or equivalently $b$ splits $t$.  A \emph{block} is a
cycle in $\cI_0'$, that is, a splitting sequence
\[
\tau_1 \searrow \tau_2 \searrow \cdots \searrow \tau_n \searrow \tau_1
\]
starting and ending at the same point.   The key definition is as follows:
\begin{definition}  
  A block is said to be \emph{isolating} if we can partition the
  labels into non-empty subsets $V = \{v_i\}$ and $W = \{w_j\}$ such that:

 \begin{enumerate}
  \item Every $v_i$ splits some $v_j$ and is split by some $v_k$.
  \item No $w$ splits a $v$.  
  \end{enumerate}
\end{definition}
Then we have:
\begin{theorem}[Kerckhoff]\label{thm-steve-crit}
  If there are no isolating blocks, then the splitting transformation
  on $Z$ is normal. 
\end{theorem}

The way we have stated things differs slightly from
\cite{Kerckhoff1985}, so we now say how to directly connect our
presentation to his work there (readers unfamiliar with
\cite{Kerckhoff1985} will want to skip ahead to the proof of
\fullref{thm-dyn-ubiq}).  First, the notion of an isolated block
is defined at the bottom of page 262 of \cite{Kerckhoff1985}.  It is
given there in terms of a partition of vertices of a simplex $\Sigma$.  We
have that $\ML(\tau, \R)$ is a convex cone which is the intersection of
$\R_{\geq 0}^7$ with a hyperplane.  Here the coordinates of $\R_{\geq
  0}^7$ correspond to the labeled bands.  Kerckhoff's simplex $\Sigma$ is
just the convex hull of the positive unit vectors along each of the
coordinate axes.  Splitting the interval exchange corresponds to
adding vertices of $\Sigma$, as detailed in the proof of Prop.~1.4 of
\cite{Kerckhoff1985}.  \fullref{thm-steve-crit} above is
essentially just Theorem~2.1 of \cite{Kerckhoff1985}, with a slight
modification because $\ML(\tau, \R)$ is not all of $\R_{\geq 0}^7$.  This
modification is justified in the 2nd paragraph of page 268 of
\cite{Kerckhoff1985}.  (The problem mentioned above with the proof of
normality for train tracks occurs later, namely in the proof of
Proposition~2.2 of \cite{Kerckhoff1985}.)  Finally, the notion of
normality is not precisely defined in \cite{Kerckhoff1985}, as it is a
standard concept in dynamical systems.  In particular, there is no
mention there of the graph structure of $\cI$ or the need to focus on
splitting sequences lying in the sink.  However, these notions are
implicit in \cite{Kerckhoff1985}, see in particular the second
sentence of the proof of Corollary 1.9.  We return now to the matter at
hand.

\begin{proof}[Proof of \fullref{thm-dyn-ubiq}]
  The proof of this theorem is a little involved, but is purely
  combinatorial, and essentially self-contained.  By
  \fullref{thm-steve-crit}, we just need to show that the
  splitting transformation on $Z$ has no isolating blocks.  Suppose to
  the contrary we have a block
  \[
  \tau_1 \searrow \tau_2 \searrow \cdots \searrow \tau_n \searrow \tau_1
  \]
  isolating band subsets $V = \{v_i\}$ and $W = \{w_j\}$ as above.
  Throughout, we will think of a $\tau_i$ as being specified by two
  lists of band labels, one each for the top and bottom interval of
  the exchange.  For instance, the standard exchange $\tau_0$ is given
  by
  \[
  \tau_0 = \frac{1 2 3 4 2 3 4}{ 5 67 5 67 1}.
  \]
  A statement like ``$v_i$ lies to the right of $w_j$'' means that one
  occurrence of $v_i$ lies to the right of $w_j$ in its list.  For
  $\tau_0$, both the statements ``2 lies to the right of 3'' and ``3
  lies to the right of 2'' are correct; they simply refer to different
  occurrences of the labels.  We begin with:
  \begin{lemma}
    No $w$ ever enters the critical position.  Moreover, on both the
    top and the bottom interval, every $v$ lies to the right of every
    $w$.
  \end{lemma}
  \begin{proof}
    First, we argue as in Proposition 1.4 of \cite{Kerckhoff1985} that no
    $w$ ever enters one of the critical positions.  By axiom (1) of
    isolation, at some stage along the block both of the critical
    positions are occupied by $V$ bands.  Reindexing, we can assume
    that this is the case for $\tau_1$.  Consider the largest $k$ so
    that $\tau_k$ has a band $w_j$ in the critical position.  As a
    single splitting only changes one of the labels in the critical
    positions, the other critical band for $\tau_k$ is some $v_i$.
    As the next stage $\tau_{k+1}$ has only $V$ bands in the critical
    positions, for $\tau_k \searrow \tau_{k+1}$ we must have $w_j$ splitting
    $v_i$, violating axiom (2).  So no $w_j$ ever enters a critical
    position.
    
    The rest of the proof of the lemma is based on considering the
    following quantity:
    \[
    C = \num{ \mbox{$v_i$ to the right of all $w_j$ on top} } + \num{ \mbox{$v_i$ to the right of all $w_j$ on bottom} }.
      \]
      How does $C$ change as we split $\tau_k$ to $\tau_{k+1}$?  For
      notation, suppose $v_i$ splits $v_j$.  Then the
      critical end of $v_i$ is removed, decreasing $C$ by $1$.  On the
      other hand, the non-critical end of $v_j$ is divided into two,
      and so overall either:
      \begin{enumerate}\renewcommand{\labelenumi}{(\alph{enumi})}
      \item The non-critical end of $v_j$ lies to the right of all
        $w_k$, in which case $C$ is unchanged.
      \item The non-critical end of $v_j$ lies to left of some $w_k$,
        and $C$ decreases by $1$.  
      \end{enumerate}
      Axiom (2) of isolation forces each $v_i$ to be split, and hence
      case (b) does occur somewhere along the block.  But as $C$ is
      non-increasing, this gives a contradiction as the block starts
      and ends at the same interval exchange, and so $C$ must be
      unchanged after running all the way through the block.
    \end{proof}
  
    The above lemma says that $\tau_k$ splits up, in a certain sense,
    into two distinct interval exchanges which have been stuck next to
    each other.  Here one exchange consists of the $W$ bands, and the
    other of the $V$ bands; we denote them by $W$ and $V_k$
    respectively ($W$ is unchanging and so needs no subscript).  For
    example, we might have:
    \begin{equation}\label{eq-almost-iso}
    \frac{ w_1 w_2 w_3 w_1 | v_1 v_2 v_3 v_1}{w_3 w_2 | v_4 v_3 v_2 v_4}.
    \end{equation}
    Note that $W$ and $V_k$ may be quite degenerate as interval
    exchanges, in particular they need not be recurrent. Moreover,
    this decomposition is purely at the combinatorial level; a measure
    $\mu \in \ML(\tau_k, \R)$ is typically not the result of taking $\mu_w
    \in \ML(W, \R)$ and $\mu_v \in \ML(V_k, \R)$ and amalgamating them.
    We now work to acquire more information about $W$ and the $V_k$,
    eventually deriving a contradiction.

    \begin{lemma}\label{lemma-rev-bands}
      Both $W$ and $V_k$ have orientation reversing bands.
      Moreover, $V_k$ has such bands on both top and bottom.
    \end{lemma}
    \begin{proof}
      Suppose $W$ has only orientation preserving bands.  Then for all
      measure laminations $\mu \in \ML(\tau_k,\R)$ the gaps on the top and
      bottom between the $W$ bands and the $V$ bands line up exactly.
      Thus we can slice through $\tau_k$ at that point to get a
      (generalized) train track carrying $\mu$ which is the disjoint
      union of $W$ and $V$.  Since the complementary regions to $\tau_k$
      are ideal triangles, it follows that one of the complementary
      regions to $\mu$ is \emph{not} an ideal triangle.  But as this is
      true for \emph{every} $\mu$, the exchange $\tau_k$ cannot be
      complete as laminations with triangle complementary regions are
      dense in $\ML(\Sigma)$.  So $W$ has an orientation reversing band.
      
      The same argument shows that $V_k$ must have an orientation
      reversing band.  Suppose it has such a band only on one side,
      say the top.  As it is not possible to create a reversing band
      without a reversing band on the same side, it follows that none
      of the $V_k$ have a reversing band on the bottom.  Since there
      are no reversing bands on the bottom, we cannot ever decrease
      the number of such bands on the top.  As the block starts and
      ends at the same exchange, it follows that the number of
      reversing bands is constant throughout.  However, by axiom (1),
      at some point one of the reversing bands is split, necessarily
      by an orientation preserving band; as this increases the number
      of reversing bands, we have a contradiction.  So each $V_k$ must
      have orientation reversing bands on both sides.
    \end{proof}
    
    \begin{figure}
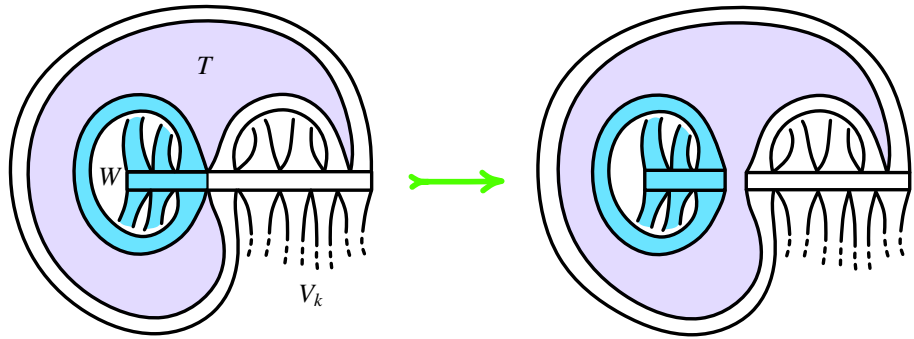

      \centering
      {\small
      \centering\leavevmode
      \begin{xyoverpic*}{(401,148)}{scale=0.85}{\figdir/outside}
        ,(87,120)*{T}
        ,(45,71)*{W}
        ,(134,18)*{V_k}
      \end{xyoverpic*}
      }
      \caption{
        Splitting apart $W$ and $V_k$ results in amalgamating two
        vertices of the triangle $T$ into a punctured monogon.
      }\label{fig-cut-tri}
    \end{figure}
    What are the complementary regions of $W$ and $V_k$, thought of as
    abstract interval exchange?  If we slice across to separate $\tau_k$
    into $W$ and $V_k$, we amalgamate two ideal triangles at a pair of
    vertices.  As this cutting separates $\tau_k$ into two pieces, in
    fact we must be amalgamating two ideal vertices of the \emph{same}
    ideal triangle, as shown in \fullref{fig-cut-tri}.  Thus the
    complementary regions to $W$ or $V_k$ are all ideal triangles with
    one exceptional region, which we call the \emph{outside} region.
    One of the outside regions is an ideal monogon, and the other is
    smooth.  Next, we use this to show:

    \begin{lemma}  
      The exchange $W$ has at least 3 bands, and $V_k$ has at least 4 bands.  
    \end{lemma}
    \begin{proof}
      Let's begin with $W$, which we know has an orientation reversing
      band.  If $W$ has only that one band, then the outside complementary
      region is smooth, but the interior complementary region is an
      ideal monogon.  So we must add a second band with one end glued
      inside the monogon to ``break'' it.  There are two possibilities
      depending on whether the second band reverses or preserves
      orientation, but these have outside regions a triangle and digon
      respectively.  So $W$ must have at least three bands.  
      
      Turning to $V_K$, we know it has orientation reversing bands
      on both sides.  The only way to get rid of the associated
      interior monogons by adding a single additional band is
      \[
      \frac{v_1 v_2 v_1}{v_3 v_2 v_3}
      \]
      which has a digon complementary region.  So $V_k$ needs at
      least 4 bands.
    \end{proof}
    
    Now, a complete genus 2 exchange has $7$ bands, so it follows from
    the lemma that $V_k$ consists of exactly 4 bands.  We still have
    work to do, as the example of \eqref{eq-almost-iso} is consistent
    with what we know so far.  To conclude the proof, we consider the
    possible choices for the number $r$ of reversing bands of $V_k$.
    \begin{enumerate}\leftskip 10pt
    \item[$r = 4$:] To avoid interior monogons, we must have two bands
      on each side, interlocking to form 
      \[
      \frac{v_1 v_2 v_1 v_2}{v_3 v_4 v_3 v_4}.
      \]
      But then the only complementary region is a hexagon.  So no
      $V_k$ has $r = 4$.
      
    \item[$r = 3$:] Since we have reversing bands on both sides, we
      can assume we have two such bands on the top and one on the
      bottom.  The top reversing bands must interlock since we have
      only a single additional orientation preserving band to break any
      interior monogons and digons.   Thus we must have
      \[
      \frac{v_1 v_2 v_1 v_2}{v_3 v_4 v_3}
      \]
      where the top end of $v_4$ has not yet been attached.  There are
      5 possibilities for the placement of $v_4$, and an easy check
      shows that all of them have the wrong complementary regions.
      Thus no $V_k$ has $r = 3.$
      
    \item[$r=2$:] As the example of \eqref{eq-almost-iso} shows, this
      is possible as a stand-alone exchange; we argue instead that
      such an exchange cannot lie in a isolating block --- that is, to
      return to where we start in a block we must involve some of the
      $W$ bands.  By \fullref{lemma-rev-bands} and the above cases,
      we must have that every $V_k$ has exactly two reversing bands,
      one on each side.  But by axiom (1) at some point a reversing
      band is split, which either creates or destroys an orientation
      reversing band, a contradiction.
      \end{enumerate}
      
      \noindent
      Thus we have ruled out all possibilities for an isolating
      block, proving the theorem.
\end{proof}

\begin{remark}
  For surfaces of higher genus, it seems very likely that there are
  isolating blocks.  Indeed, in genus 3 take
  \[
  W = \frac{w_1 w_2 w_3 w_1}{w_4 w_3 w_2 w_4}
  \]
  and let $V_1$ be a complete exchange for a once punctured genus 2
  surface.  If the monogon complementary region is in the right place,
  then $\tau_1 = W \cup V_1$ is a complete genus 3 exchange.  If $V_1$ is
  in some closed strongly connected subgraph of such exchanges, then
  completeness will allow us to construct a splitting sequence of
  $V_1$ back to itself which satisfies axiom (1).  What is not entirely
  clear is whether this isolating block is in the sink for
  genus 3 exchanges.

\end{remark}

\subsection[Proof of \ref{thm-ubiquity}]{Proof of \fullref{thm-ubiquity}}

We end this section by deriving \fullref{thm-ubiquity} from
\fullref{thm-dyn-ubiq}.

\begin{proof}[Proof of \fullref{thm-ubiquity}]
  As in the statement of the theorem, let $\tau_0$ be our particular
  complete genus 2 exchange, and $S$ some splitting sequence of $\tau_0$
  consisting of complete exchanges.  Set
  \[
  \cS = \setdef{ \gamma \in \ML(\tau_0, \Z)}{ \mbox{$\gamma$ exhibits $S$}}.
  \]
  Let $U$ be a bounded open set in $\ML(\tau_0, \R)$.  We need to show
  \[
  \frac{  \num{ \cS \cap t U} } { \num{ \ML(\tau_0, \Z) \cap t U}} %
    \to 1 \mtext{as} t \to \infty.
  \]
  Let $V_m$ be the set of $\mu \in \ML(\tau_0, \R)$ such that the
  splitting sequence of $\mu$ passes through $S$ ending at the $m^{\mathrm{th}}$
  stage.  The set $V_m$ is defined by a sequence of strict
  inequalities in the weights of the bands in the successive critical
  positions; in particular, it is open.  If we set $V = \bigcup V_m$, then
  by \fullref{thm-dyn-ubiq}, the complement of $V$ has measure $0$
  in $\ML(\tau_0, \R)$.  Moreover, as $\cS = \ML(\tau_0, \Z) \cap V$ and $V$ is
  invariant under positive scaling
  \[
   \cS \cap t U =  \ML(\tau_0, \Z) \cap V \cap t U = \ML(\tau_0 , \Z) \cap t (V \cap U).
  \]
  \begin{align*}
\tag*{\hbox{Therefore}}
  \frac{  \num{ \cS \cap t U} } { \num{ \ML(\tau_0, \Z) \cap t U}}  &= %
\frac{  \num{ (V \cap U ) \cap t^{-1} \ML(\tau_0,\Z)} } { \num{U \cap t^{-1} \ML(\tau_0, \Z)}} \\%
&= \frac{ t^{-6} \cdot  \num{ (V \cap U ) \cap t^{-1} \ML(\tau_0,\Z)} } {t^{-6} \cdot \num{U \cap t^{-1} \ML(\tau_0, \Z)}}
  \end{align*}
  As $t \to \infty$, the top and bottom of the right-hand fraction converge
  to the Lebesgue measures of $V \cap U$ and $U$ respectively, since
  these sets are open.  As the complement of $V$ has measure $0$, the
  sets $V \cap U$ and $U$ have the same measure.  Hence 
  the fraction limits to $1$ as $t \to \infty$, completing the proof.
\end{proof}

\section{Proof of the main theorems}\label{sec-main-thm-pf}

Here, we complete the proofs of Theorems~\ref{thm-main-nonsep} and
\ref{thm-main-sep}.  In addition to the results of
Sections~\ref{sec-magic-sequence} and \ref{sec-ubiquity}, we need one
additional ingredient, which says that specific types of
multicurves, eg,~connected non-separating curves, have positive
density in the set of all multicurves:
\begin{theorem}[Mirzakhani]\label{thm-maryam}
  Let $\Sigma$ be a closed surface of genus $g \geq 2$.  Fix a multicurve
  $\gamma \in \ML(\Sigma, \Z)$.  Consider the set $\cC = \MCG(\Sigma) \cdot \gamma$
  of all multicurves in $\Sigma$ of the same topological type.  Then for
  every bounded open subset $U$ of $\ML(\Sigma, \R)$ we have:
  \[
  \frac{ \num{ \cC \cap t U}}{\num{\ML(\Sigma, \Z) \cap t U}} \to d_{\gamma} \in
    \frac{1}{\pi^{6g - 6}} \Q,
  \]
  where $d_\gamma$ is positive and  independent of $U$.  
\end{theorem}
For us, the exact value of $d_\gamma$ will not be important, merely the
fact that it is positive; however, Mirzakhani does provide a recursive
procedure for computing it.  The above theorem is a slight
restatement of Theorem~6.4 of \cite{Mirzakhani2004} where we have set
$d_\gamma = c_\gamma/b_{g,0}$ in the notation there, and rewritten the measure
$\mu_{t,\gamma}$ in terms of $\num{\ML(\Sigma, \Z) \cap t U}$ rather than
$t^{6g - 6}$ via the proof of Theorem~3.1 of \cite{Mirzakhani2004}.
 
Let us now recall the statement of \fullref{thm-main-nonsep} and
the notation of \fullref{subsec:measured-lam-prob}.  Let $H$ be
our genus 2 handlebody, and fix Dehn--Thurston coordinates 
\[
(w_\alpha,w_\beta, w_\delta, \theta_\alpha, \theta_\beta, \theta_\delta)
\]
 for $\partial H$ compatible with $H$
as discussed in \fullref{subsec:measured-lam-prob} and shown in
\fullref{fig-DT-coor}.  We are interested in the set $\T$ of
attaching curves $\gamma \subset \partial H$ parameterizing tunnel number one
3--manifolds with one boundary component, where $\gamma$ satisfies the
restrictions:
\begin{enumerate}
  \item $\gamma$ is a non-separating simple closed curve.
  \item The weights $w_\alpha, w_\delta, w_\beta$ are $> 0$. 
  \item Each twist satisfies $0 \leq \theta < w$. 
  \item $w_\delta \leq \min( 2 w_\alpha , 2 w_\beta )$.
\end{enumerate}
We then consider the finite set $\T(r)$ of $\gamma \in \T$ where $w_\alpha +
w_\beta < r$.  \fullref{thm-main-nonsep} is that the
probability that $M_\gamma \in \T(r)$ fibers over $S^1$ goes to $0$ as $r
\to \infty$.

To prove this, we first reinterpret the setup in the context of
$\ML(\partial H, \R)$.  The Dehn--Thurston coordinates on $\ML(\partial H, \Z)$ have
natural extensions to coordinates on all of $\ML(\partial H, \R)$ \cite[Theorem
3.11]{PennerHarerBook}, which we denote in the same way.  Let $W \subset
\ML(\partial H, \R)$ consist of measured laminations satisfying conditions
(2--4) above.  The subset $W$ is a polyhedral cone in $\ML(\partial H, \R) \cong
\R^6$ with some faces removed.   Set 
\[
\cG = \setdef{\gamma \in \ML(\Sigma, \Z)}{\mbox{$\gamma$ is connected and non-separating }}.
\]
We have $\T = W \cap \cG$, and if we let $U$ be the open subset of
$\ML(\partial H, \R)$ defined by $w_\alpha + w_\beta < r$ then
\[
\T(r) = \T \cap r U = W \cap \G \cap r U.
\]
A slightly more general result immediately implying
\fullref{thm-main-nonsep} is:
\begin{theorem}\label{thm-main-nonsep-restated}
  Let $U \subset \ML(\partial H, \R)$ be an open set such that $W \cap U$ is bounded.
  Then the probability that $M_\gamma$ fibers over the circle for $\gamma
  \in \T \cap rU$ goes to $0$ as $r \to \infty$.
\end{theorem}

The rest of this section is devoted to the proof of
\fullref{thm-main-nonsep-restated}.  The proof of
\fullref{thm-main-sep} which concerns tunnel number one
3--manifolds with two boundary components is identical if one simply
replaces all occurrences of ``non-separating'' with ``separating''.  

\begin{proof}[Proof of \fullref{thm-main-nonsep-restated}]
  Let $S$ be the magic splitting sequence of \fullref{lemma-magic}.
  Let $\cS$ consist of those laminations $\mu \in \ML(\partial H, \R)$ which
  can be carried by some tightly labeled interval exchange $\tau$
  where the splitting sequence of $(\tau, \mu)$ exhibits $S$.  By
  \fullref{lemma-magic}, for $\gamma \in \cS \cap \T$ the manifold $M_\gamma$
  does not fiber.  So if $\cS^c$ denotes the complement of $\cS$, it
  suffices to show
  \begin{equation}\label{eqn-goal-1}
    \frac{\num{\cS^c \cap \T \cap r U}}{\num{\T \cap rU}} \to 0 \quad \mbox{as} \quad r \to \infty.
  \end{equation}
  We next show that \fullref{thm-maryam} allows us
  to replace $\T$ with $\ML(\partial H, \Z) \cap V$ in the above limit, where $V =
  W \cap U$.  In particular
  \begin{align*}
  \frac{\num{\cS^c \cap \T \cap r U}}{\num{\T \cap rU}} &\leq \frac{\num{\cS^c
      \cap \ML(\Sigma, \Z) \cap r V}}{\num{\G \cap rV}} \\
  &= \frac{\num{\cS^c \cap
      \ML(\Sigma, \Z) \cap r V}}{\num{\ML(\Sigma, \Z) \cap rV}} \cdot
  \frac{\num{\ML(\Sigma, \Z) \cap rV}}{\num{\G \cap rV} }
  \end{align*}
  By \fullref{thm-maryam}, the second factor in the final
  expression converges to a positive number as $r \to \infty$.  Thus to
  show \eqref{eqn-goal-1} it suffices to prove
    \begin{equation}\label{eqn-goal-2}
    \frac{\num{\cS^c \cap \ML(\partial H, \Z) \cap r V}}{\num{\ML(\partial H, \Z) \cap rV}} \to 0 \quad \mbox{as} \quad r \to \infty.
  \end{equation}
  
  Now by \fullref{lemma-std-track}, $W$ is covered by a countable
  collection of charts $\ML(\tau,\R)$ where $\tau$ is one of the
  standard tightly labeled interval exchanges.  In each chart, the
  proof of \fullref{thm-ubiquity} shows that
  \fullref{thm-dyn-ubiq} gives a \emph{homogeneous open} set $Y_\tau \subset \ML(\tau,
  \R) \cap \cS$ whose complement in $\ML(\tau, \R)$ has measure $0$.   Then 
  \[
  Y = \bigcup_{\tau} \mathrm{int}(Y_\tau) \subset \cS
  \]
  is an open subset of $W$ whose complement has measure $0$.  As in
  the proof of \fullref{thm-ubiquity}, it is easy to show that 
  \begin{equation}
    \frac{\num{Y^c \cap \ML(\partial H, \Z) \cap r V}}{\num{\ML(\partial H, \Z) \cap rV}} \to 0 \quad \mbox{as} \quad r \to \infty, 
  \end{equation}
  since the top and bottom converge to the Lebesgue measure of $Y^c \cap
  V$ and $V$ respectively.  This implies \eqref{eqn-goal-2} and hence the theorem.  
\end{proof}

\bibliographystyle{gtart} 
\bibliography{link}

\end{document}